\documentstyle[12pt,bezier]{book}
\newcounter{nonsense}[section]
\newcommand{\Phrase}[2]{\begin{trivlist}
   \def\thenonsense{\thesection.\arabic{nonsense}}
   \refstepcounter{nonsense}
   \item[\hskip\labelsep\bf\thenonsense\ #1] #2
   \end{trivlist}}
\newcommand{\Lemma}[1]{\Phrase{Lemma}{\em #1}}
\newcommand{\Prop}[1]{\Phrase{Proposition}{\em #1}}

\newcommand{\Thm}[1]{\Phrase{Theorem}{\em #1}}
\newcommand{\Cor}[1]{\Phrase{Corollary}{\em #1}}
\newcommand{\Def}[1]{\Phrase{Definition}{#1}}

\newfont{\sfb}{cmssbx10 scaled 1200}
\expandafter\chardef\csname pre amssym.def at\endcsname=\the\catcode`\@
\catcode`\@=11

\def\undefine#1{\let#1\undefined}
\def\newsymbol#1#2#3#4#5{\let\next@\relax
 \ifnum#2=\@ne\let\next@\msafam@\else
 \ifnum#2=\tw@\let\next@\msbfam@\fi\fi
 \mathchardef#1="#3\next@#4#5}
\def\mathhexbox@#1#2#3{\relax
 \ifmmode\mathpalette{}{\m@th\mathchar"#1#2#3}%
 \else\leavevmode\hbox{$\m@th\mathchar"#1#2#3$}\fi}
\def\hexnumber@#1{\ifcase#1 0\or 1\or 2\or 3\or 4\or 5\or 6\or 7\or 8\or
 9\or A\or B\or C\or D\or E\or F\fi}

\font\twlmsb=msbm10 scaled \magstep1
\font\egtmsb=msbm8 
\font\sixmsb=msbm6
\newfam\msbfam
\textfont\msbfam=\twlmsb
\scriptfont\msbfam=\egtmsb
\scriptscriptfont\msbfam=\sixmsb
\edef\msbfam@{\hexnumber@\msbfam}
\def\Bbb#1{\fam\msbfam\relax#1}

\font\twleusb=eusb10 scaled \magstep1
\font\egteusb=eusb8 
\font\sixeusb=eusb6
\newfam\eusbfam
\textfont\eusbfam=\twleusb
\scriptfont\eusbfam=\egteusb
\scriptscriptfont\eusbfam=\sixeusb

\font\twleufm=eufm10 scaled \magstep1
\font\egteufm=eufm8
\font\sixeufm=eufm6
\newfam\eufmfam
\textfont\eufmfam=\twleufm
\scriptfont\eufmfam=\egteufm
\scriptscriptfont\eufmfam=\sixeufm
\def\frak#1{{\fam\eufmfam\relax#1}}

\let\ps@plain\ps@empty

\font\gross=eufm10 scaled \magstep5
\newbox\initiale

\def\boxit#1{\vbox{\hrule\hbox{\vrule\kern3pt
\vbox{\kern3pt#1\kern3pt}\kern3pt\vrule}\hrule}}

\def\Initiale#1{\setbox\initiale=%
\hbox{\hbox{\vrule%
\vbox{\hrule\kern3pt%
\hbox{\kern3pt{\gross #1}\kern3pt}%
\kern3pt\hrule}\vrule}\kern2pt}%
\hangindent=\wd\initiale\hangafter=-2%
\ \unskip\hskip-\wd\initiale%
\dp\initiale=-\baselineskip\lower\baselineskip\box\initiale\kern-1.5pt}

\catcode`\@=\csname pre amssym.def at\endcsname

\newcommand{\imb}{\hookrightarrow}
\newcommand{\dto}{\downarrow}
\newcommand{\uto}{\uparrow}
\newcommand{\DTO}[1]{\downarrow\rlap{$\scriptstyle#1$}}

\newcommand{\bra}[1]{\langle#1\rangle}
\newcommand{\eps}{\varepsilon}
\newcommand{\RR}{{\Bbb{R}}}
\newcommand{\SS}{{\Bbb{S}}}
\newcommand{\ZZ}{{\Bbb{Z}}}
\newcommand{\NN}{{\Bbb{N}}}

\newcommand{\CC}{{\Bbb{C}}}
\newcommand{\HH}{{\Bbb{H}}}
\newcommand{\II}{{\Bbb{I}}}
\newcommand{\QQ}{{\Bbb{Q}}}

\newcommand{\OO}{{\Bbb{O}}}

\newcommand{\ot}{\leftarrow}

\newcommand{\too}{\longrightarrow}
\newcommand{\TO}[1]{\stackrel{#1}{\to}}

\newcommand{\OT}[1]{\stackrel{#1}{\ot}}

\newcommand{\SUB}{\subseteq}
\newcommand{\proof}{\par\medskip\rm{\em Proof. }}
\newcommand{\qed}{\ \hglue 0pt plus 1filll $\Box$}

\renewcommand{\P}{{\cal P}}
\renewcommand{\L}{{\cal L}}
\newcommand{\F}{{\cal F}}
\newcommand{\V}{{\cal V}}
\newcommand{\PLF}{(\P,\L,\F)}
\newcommand{\Gal}{{\rm Gall}}
\newcommand{\StamGal}{{\rm StamGall}}
\newcommand{\PropGal}{{\rm PropGall}}
\newcommand{\opp}{{\rm opp\,}}
\newcommand{\flag}{{\rm fl}}
\newcommand{\Cl}{{\rm Cl}}
\newcommand{\id}{{\rm id}}
\renewcommand{\H}{{\bf H}}
\newcommand{\Hbar}{{\bf\bar H}}
\newcommand{\p}{{\bf p}}
\newcommand{\e}{{\bf e}}
\newcommand{\pr}{{\rm pr}}

\renewcommand{\aa}{{\frak a}}
\newcommand{\bb}{{\frak b}}
\newcommand{\cc}{{\frak c}}
\newcommand{\dd}{{\frak d}}

\newcommand{\ff}{{\frak f}}
\renewcommand{\gg}{{\frak g}}
\newcommand{\Sq}{{\rm Sq}}

\let\perp\bot

\begin{document}

\begin{titlepage}

\
\vspace{6cm}
\center{\Huge\bf Compact Polygons}

\vspace{4cm}
\begin{center}
{\Large\sc
Dissertation der Mathematischen Fakult\"at \\
der Eberhard-Karls-Universit\"at T\"ubingen \\
zur Erlangung des Grades eines Doktors \\
der Naturwissenschaften \\
\ }
\end{center}

\vspace{2cm}
\begin{center}
{\Large\sc
vorgelegt von \\
Linus Kramer \\
aus Marburg an der Lahn \\
\ }
\end{center}

\vfill
\center{\Large 1994}

\end{titlepage}

\pagenumbering{roman}

\chapter*{Einleitung}
\addcontentsline{toc}{chapter}{\protect\numberline{}{Einleitung}}

Verallgemeinerte Polygone sind nat\"urliche Verallgemeinerungen
projektiver Ebenen; an die Stelle des Axioms, da{\ss} zwei Punkte
sich durch eine eindeutig bestimmte Ge\-ra\-de verbinden lassen, tritt
die Forderung, da{\ss} zwei Elemente (Punkte oder Geraden) sich
durch einen eindeutig bestimmten Polygonzug einer gewissen L\"ange
verbinden lassen. Projektive Ebenen sind verallgemeinerte Dreiecke.
'Erfunden' wurden ver\-all\-ge\-mei\-ner\-te Polygone (oder
sph\"arische Geb\"aude vom Rang 2) von Tits, um eine geometrische Interpretation
f\"ur gewisse einfache (Ausnahme-) Liegruppen zu
finden; sph\"arische Geb\"aude liefern eine einheitliche geometrische Interpretation
f\"ur alle einfachen algebraischen Gruppen. Nun tritt bei sph\"arischen
Geb\"auden das gleiche Ph\"anomen wie in der projektiven Geometrie auf:
die sph\"a\-ri\-schen Ge\-b\"au\-de vom Rang $\geq 3$ sind in einem gewissen Sinne
klassisch, das hei{\ss}t, sie stammen von klassischen oder von einfachen algebraischen Gruppen.
Verallgemeinerte Polygone sind dagegen (wie projektive Ebenen) im
allgemeinen nicht klassisch, und gerade deswegen interessant.

Will man substantielle Aussagen \"uber verallgemeinerte Polygone machen, 
so sind zus\"atzliche Strukturannahmen wie etwa Endlichkeit
oder die Existenz gewisser Automorphismen notwendig. 
Eine andere M\"oglichkeit ist,
das Polygon zu\-s\"atz\-lich mit einer topologischen Struktur zu versehen
und zu verlangen, da{\ss} der oben er\-w\"ahn\-te Polygonzug stetig von
seinen Enden abh\"angt. Im wesentlichen ist das die Definition eines
topologischen  Polygons (in dieser Dissertation wird etwas weniger Stetigkeit
verlangt). Bereits dieser recht allgemeine Ansatz liefert
Einschr\"ankungen f\"ur die Topologie: der unterliegende topologische Raum wird
regul\"ar und ist entweder zusammenh\"angend oder total
unzusammenh\"angend. Ist das Polygon lokalkompakt, so wird die Topologie
des Polygons metrisierbar und abz\"ahlbar. 

Wesentlich st\"arkere Aussagen \"uber die globale topologische
Struktur ergeben sich mit Hilfe der algebraischen Topologie, wenn man
voraussetzt, da{\ss} die Topologie auf dem Punktraum des verallgemeinerten $n$-Eckes lokalkompakt und
endlichdimensional ist. In diesem Falle ist $n\in\{3,4,6\}$, das hei{\ss}t, es gibt unter diesen topologischen Vor\-aus\-setzungen neben den projektiven Ebenen \"uberhaupt nur verallgemeinerte Vier- und Sechsecke (rein inzidenzgeometrisch
lassen sich mit Hilfe freier Konstruktionen verallgemeinerte $n$-Ecke f\"ur
jedes $n$ konstruieren). 

Eine unmittelbare Anwendung dieser Ergebnisse \"uber die algebraische
Topologie verallgemeinerter Polygone ist die Klassifikation
punkthomogener kompakter Polygone mit gleichen
topologischen Parametern. Dieses Resultat beinhaltet als Spezialfall
die Klassifikation der punkthomogenen kompakten projektiven Ebenen, die
in den siebziger und achtziger Jahren von Salzmann und L\"owen abgeschlossen
wurde (im Prinzip l\"a{\ss}t es sich auch als Klassifikation der
homogenen Fokalmannigfaltigkeiten isoparametrischer Hyperfl\"achen mit
gleichen Multiplizit\"aten deuten). Im Folgenden sollen die Methoden etwas
genauer beschrieben werden.

\medskip
Um die Topologie auf dem Punktraum $\P$ eines topologischen $n$-Eckes
zu verstehen, wird zun\"achst die induzierte Topologie auf den Punktreihen
betrachtet. Sticht man einen beliebigen Punkt $\infty_L$ aus der Punktreihe $L$ heraus, dann
l\"a{\ss}t sich auf dem Komplement $L-\{\infty_L\}$ eine Addition einf\"uhren, die
es erlaubt, Gleichungen (von einer Seite) stetig zu l\"osen. Zusammen
mit der zweifachen Homogenit\"at von $L$ unter der Projektivit\"atengruppe
folgt daraus (ganz \"ahnlich wie bei topologischen Gruppen)
die Regularit\"at von $L$. 

Weiter kann man eine (sehr schwache) Multiplikation einf\"uhren; wenn
$L$ weg\-zu\-sam\-men\-h\"an\-gend ist, dann l\"a{\ss}t sich die affine Punktreihe
$L-\{\infty_L\}$ mit Hilfe dieser Multiplikation gleichm\"a{\ss}ig auf das Neutralelement der
Addition zusammenziehen.

W\"ahlt man eine Fahne $(p,\ell)$, so erh\"alt man eine
Filtrierung des Punkt\-raumes
\[
\{p\}\SUB L\SUB p^\perp \SUB  \ldots \SUB \P
\]
durch $n$ Mengen von Punkten, die sich ungef\"ahr wie folgt beschreiben lassen:
die Ele\-men\-te der $k$-ten Menge lassen sich durch einen Polygonzug
der L\"ange $\leq k$ mit der Fahne $(p,\ell)$ verbinden.
Die Mengen dieser Filtrierung sind die {\em Schubertvariet\"aten} des
$n$-Eckes. Jede Schubertvariet\"at ist eine abgeschlossene Teilmenge
des Punkt\-raumes, und die mengentheoretische Differenz zweier aufeinanderfolgender 
Schubertvariet\"aten (die zugeh\"orige {\em Schubertzelle}) ist
ein Produkt affiner Punktreihen und affiner Gera\-den\-b\"u\-schel.
Insbesondere ist der Punktraum $\P$ lokal ein Produkt aus Punktreihen
und Ge\-ra\-den\-b\"u\-scheln. Alle Eigenschaften der mengentheoretischen
Topologie von $\P$ ergeben sich aus diesem Resultat.

F\"ur die algebraische Topologie von $\P$ ist das folgende Theorem von L\"owen
von zentraler Bedeutung: {\em ist $X$ ein endlichdimensionaler, kompakter
ANR (absoluter Umgebungsretrakt) mit der Eigenschaft, da{\ss} das Komplement jedes Elementes $x\in X$
kontrahierbar ist, dann ist $X$ eine verallgemeinerte Mannigfaltigkeit und
homotopie\"aquivalent zu einer Sph\"are.}

Dieses Kriterium l\"a{\ss}t sich auf die Punktreihen und Geradenb\"uschel
eines ver\-all\-gemeinerten $n$-Eckes anwenden. Damit wird zun\"achst jede Schubertzelle eine ver\-all\-ge\-mei\-ner\-te Mannigfaltigkeit, deren
Ein-Punkt-Kompaktifizierung eine Ho\-mo\-to\-pie\-sph\"a\-re 
ist. Wenn $n>3$ ist, dann sind
die Schubertvariet\"aten im allgemeinen keine verallgemeinerten
Mannigfaltigkeiten mehr; der Punktstern $p^\perp$ etwa, also die Menge aller
Punkte, die mit $p$ durch eine Gerade verbindbar sind,  hat in seinem
Mittelpunkt $p$ im allgemeinen eine Singularit\"at. Diese Singularit\"aten lassen
sich "gl\"atten", indem man die Schubertvariet\"at $p^\perp$ durch den Raum der
zugeh\"origen Galerien ersetzt; 
\[
\Gal_2(\ell,p)\cong\{(h,q)|\ h\mbox{ geht durch }p\mbox{ und }
q\mbox{ liegt auf }h\}
\]
ist eine verallgemeinerte Mannigfaltigkeit, und $p^\perp\cong
\Gal_2(\ell,p)/(\L_p\times\{p\})$ ist ein Quotient diese Raumes.
Diese Beobachtung, verbunden mit der Tatsache, da{\ss} die Menge der
$k$-Galerien ein $k$-fach iteriertes B\"undel ist, erlaubt die
Berechnung der Homologiegruppen der Schubertvariet\"aten; es ergibt
sich, da{\ss} jede Schubertvariet\"at im Punktraum einen Erzeuger
der Homologiegruppen repr\"asentiert, ein Resultat, das f\"ur
projektive R\"aume (oder Gra{\ss}mann-Mannigfaltigkeiten)
wohlbekannt ist.

Sch\"arfere Aussagen \"uber die algebraische Topologie erh\"alt man
mit Knarrs to\-po\-lo\-gi\-scher Veroneseeinbettung: der doppelte Abbildungszylinder \"uber dem Fahnenraum ist eine Homotopiesph\"are.
Mit einem Theorem von M\"unzner folgt, da{\ss} {\em die Coxetergruppe
eines lokalkompakten, endlichdimensionalen $n$-Eckes kristallographisch
ist (die Coxetergruppe eines verallgemeinerten $n$-Eckes ist eine Diedergruppe
der Ord\-nung $2n$).} Es gibt unter diesen to\-po\-lo\-gi\-schen Voraussetzungen also
nur projektive Ebenen, Vierecke und Sechsecke. Dieses Resultat
(das topologische Pendant zu den S\"atzen von
Tits-Weiss und Feit-Higman)
wurde zuerst von Knarr unter der Zusatzannahme bewiesen, da{\ss} das
$n$-Eck eine topologische Mannigfaltigkeit ist.

Setzt man voraus, da{\ss} die Punktreihen und Geradenb\"uschel
lokal euklidisch sind, dann wird die Schubertzellzerlegung
eine CW-Zerlegung. F\"ur projektive Ebenen wurde das erstmals
von Breitsprecher bewiesen; in kleinen Dimensionen lassen sich
damit gewisse Schubertvariet\"aten topologisch klassifizieren.
F\"ur projektive Ebenen wurde das bereits von Salzmann, Breitsprecher
und Buchanan gezeigt.

Eine Anwendung dieser Ergebnisse ist die geometrische Klassifikation
gewisser homogener $n$-Ecke. Gegeben sei ein lokalkompaktes,
zusammenh\"angendes verallgemeinertes $n$-Eck mit punkttransitiver 
Automorphismengruppe. Zun\"achst folgt, da{\ss} die
Automorphismengruppe eine Liegruppe ist, und da{\ss} $n=3,4,6$ ist.
F\"ur $n=3,6$ ist der Punktraum entweder $(n-1)$-dimensional, oder
er hat positive Eulercha\-rak\-teristik $n$. F\"ur Vierecke ist
diese Bedingung nicht notwendig erf\"ullt, weswegen wir sie
f\"ur $n=4$ zus\"atzlich fordern. {\em Unter diesen Voraussetzungen
ist das $n$-Eck klassisch (das hei{\ss}t Mou\-fangsch) und -- via Cartans Klassifikation der
einfachen Liegruppen -- explizit bekannt.} Der Beweis zerf\"allt in zwei Teile: wenn die Dimension des
Punkt\-rau\-mes $n-1$ ist, so folgt, da{\ss} die auftretenden
Liegruppen klein und dadurch eindeutig bestimmt sind; der andere
Fall ergibt sich mit Hilfe der
Borel-De Siebenthal Klassifikation maximaler Untergruppen von
maximalem Rang in kompakten Liegruppen.

Verallgemeinerte Vierecke nehmen hier eine Sonderstellung ein;
es gibt nicht-Moufangsche, punkthomogene verallgemeinerte
Vierecke (mit Eulercharakteristik $0$). Auch die Erfahrung,
da{\ss} die nichtklassischen kompakten Ebenen durch die
klas\-si\-schen Ebenen dominiert werden, hat f\"ur Vierecke keinerlei
Entsprechung. Kompakte Sechsecke zeigen wieder ein \"ahnliches Verhalten
wie kompakte projektive Ebenen; allerdings kennt man hier
gegenw\"artig \"uberhaupt keine nicht-Moufangschen Beispiele.

\"Uber die kompakten Vierecke mit Eulercharakteristik $0$ ist
gegenw\"artig noch wenig bekannt. Ferus-Karcher-M\"unzner und
Thorbergsson haben
gezeigt, da{\ss} jede Darstellung einer reellen Cliffordalgebra
ein solches verallgemeinertes Viereck liefert. Diese Konstruktion umfa{\ss}t
alle kompakten zusammenh\"angenden hermiteschen und reell-orthogonalen Mou\-fang-Vier\-ecke, dar\"uber hinaus aber sehr viele nicht-Mou\-fang\-sche
Beispiele. Einige dieser Vierecke haben punkttransitive
Automorphismengruppen. In dieser Richtung gibt es, auch in der Verbindung
zur Differentialgeometrie, noch interessante offene Fragen.

\medskip
Im ersten Kapitel werden alle ben\"otigten Fakten \"uber verallgemeinerte
$n$-Ecke zusammengestellt. Das zweite Kapitel besch\"aftigt sich mit der
mengentheoretischen Topologie verallgemeinerter Polygone. Einige der
Resultate sind inzwischen Folklore, aber eine einheitliche Darstellung
(f\"ur beliebiges $n$) lag bisher nicht vor. 

Das dritte Kapitel
ist das Zentrum dieser Dissertation. Es behandelt die al\-ge\-bra\-ische
Topologie endlichdimensionaler Polygone. Als 'Richtschnur' f\"ur
dieses Kapitel haben vor allem die Arbeiten von Breitsprecher und
L\"owen \"uber die algebraische Topologie projektiver Ebenen gedient.
Im vierten Kapitel werden lokal eu\-kli\-di\-sche Polygone untersucht.
Diese Polygone sind in nat\"urlicher Weise CW-Komplexe. In Verbindung
mit einem solchen Polygon treten verschiedene Vektorb\"undel
auf, die sich (in kleinen Dimensionen) zur Klassifikation bestimmter
Schubertvariet\"aten verwenden lassen.
Im f\"unften Kapitel werden homogene Polygone untersucht. 
Mit Hilfe der Ergebnisse aus Kapitel 3 werden die punkt\-ho\-mo\-ge\-nen
kompakten zusammenh\"angenden Polygone mit gleichen topologischen
Parametern vollst\"andig klassifiziert.
Im letzten Kapitel werden verschiedene topologische Hilfsmittel
be\-reit\-ge\-stellt.

\bigskip
Mathematik lebt vom Austausch und von Diskussionen; f\"ur beides
m\"ochte ich
mich bei meinen MitdoktorandInnen (und Habilitanden) Richard B\"odi,
Martina J\"ager, Michael Joswig, Bernhard M\"uhlherr und
Markus Stroppel bedanken. Norbert Knarr und
Stephan Stolz haben mir viele gute Anregungen gegeben, ohne die die
Kapitel 3, 4 und 6 sicher k\"urzer geworden w\"aren. 
F\"ur das sch\"one Thema und die vorz\"ugliche Betreuung dieser
Dissertation m\"ochte ich mich bei Theo Grundh\"ofer sehr herz\-lich
bedanken. Das Arbeiten in der AG Geometrie und Topologie war
f\"ur mich sowohl mathematisch als auch pers\"onlich immer
sehr angenehm; daf\"ur geb\"uhrt nicht zuletzt Herrn Salzmann
mein Dank.

W\"ahrend der beiden letzten Jahre hat die Studienstiftung des
deutschen Volkes mich dankenswerter Weise mit einem Promotionsstipendium unterst\"utzt.

\tableofcontents

\subsection*{Notation:}
\addcontentsline{toc}{chapter}{\protect\numberline{}{Notation}}

\begin{tabular}{ll}
$\SUB$      & subset \\
$A-B$       & set-theoretic difference $\{a\in A|\ a\not\in B\}$ \\
$\NN$       & natural numbers $\{0,1,2,3,\ldots\}$ \\
$\ZZ$       & integers \\
$\ZZ_k$     & integers mod $k$ \\
$\QQ$       & rational numbers \\
$\RR$       & real numbers \\
$\CC$       & complex numbers \\
$\HH$       & quaternions \\
$\OO$       & octonions (Cayley numbers) \\
$\II$       & unit interval \\
$\SS^k$     & $k$-dimensional unit sphere \\
$A\imb B$   & topological imbedding (homeomorphic onto its image) \\
$X/A$       & quotient space obtained by collapsing $A\SUB X$ to a point \\
$\cong$     & isomorphism in the appropriate category \\
            & (e.g. homeomorphism of topological spaces) \\
$\simeq$    & homotopy equivalence \\
$[X;Y]$     & set of all homotopy classes of maps from $X$ to $Y$ \\
$[X;Y]^0$   & set of all base-point preserving homotopy classes of
              maps \\
            & (between the pointed spaces $X$, $Y$)
\end{tabular}

\medskip\noindent
We use singular homology and cohomology, unless indicated otherwise.

\clearpage
\pagenumbering{arabic}
\setcounter{page}{1}

\chapter{Generalized polygons}

The first six sections of this chapter contain the basic
definitions concerning incidence structures and generalized polygons.
The most important
notions are the space of galleries $\Gal_k(u,v)$ based on a flag,
and the corresponding Schubert cells and Schubert varieties. Our
approach to the topology of generalized polygons is wholly based
on these spaces. 

In Section 7 it is shown that each Schubert cell is a
product of punctured pencils of lines and punctured point rows;
this decomposition is natural with respect to automorphisms, once
an ordinary $n$-gon is chosen (the Schubert cells are precisely
the preimages of the vertices of this ordinary $n$-gon with
respect to the associated retraction). 

In the last section we define the addition and multiplication
on the point rows and pencils of lines. It is well-known that
in the case of a projective plane, these algebraic operations
belong to a ternary field, and that there are strong analogies to
ordinary fields. For $n>3$ the algebraic properties of these
binary operations are much weaker; the addition yields only
a right loop, i.e. equations may be solved only from the right-hand
side. The multiplication has even weaker properties. However, it
turns out that these algebraic operations still have all the properties
which are needed to carry over most of the concepts used in the
theory of topological projective planes to generalized polygons.

\medskip
Generalized polygons were introduced by Tits in \cite{Ti0}. See also
\cite{Ti1}, \cite{Bro89}, \cite{Ron89}.

\section{Incidence structures}

\label{GenPoly1_1}
An {\em incidence structure} is a triple $\frak P=\PLF$, consisting of a
{\em point space} $\P$, a {\em line space} $\L$, and a {\em flag space}
$\F\SUB\P\times\L$. We assume that $\P$ and $\L$ are nonempty, disjoint
sets. The union $\V=\P\cup\L$ is called the set of {\em vertices}.
Two vertices $x,y\in\V$ are called {\em incident}, if $(x,y)$ or $(y,x)$ is
a flag; we denote the corresponding flag by $\flag(x,y)$. If $(p,\ell)$ is a 
flag, we say that the point $p$ lies on the line $\ell$, and that the line 
$\ell$ passes through the point $p$.

Given a set of vertices $A\SUB\V$, we let $\V_A$ denote the set of all 
vertices that are incident with some member of $A$; we put 
$\P_A=\V_A\cap\P$, $\L_A=\V_A\cap\L$, and 
$\F_A=(\V_A\times\V\cup\V\times\V_A)\cap\F$. In the special case that $A$ 
consists of a single point $p$, we call $\L_p$ the {\em pencil of lines} 
through $p$; similarly, if $A$ consists of a single line $\ell$, then 
$L=\P_\ell$ is called the {\em point row} consisting of all points 
lying on the line $\ell$. It is customary to put $x^\perp=\{y\in\V|\ 
y\mbox{ is incident with some }z\in\V_x\}$; 
this set is sometimes called the {\em star \em or \em perp} of $x$.

The incidence structure $\frak P$ is called {\em thick}, if every
point row and every pencil of lines contains at least $3$ elements.

A {\em $k$-chain} joining $x_0$ and $x_k$ is a sequence ${\bf x}=
(x_0,x_1,\ldots,x_k)$ of vertices with the property that $x_{i-1}$ is 
incident with $x_i$ for $1\leq i\leq k$. A $k$-chain {\em stammers}, 
if $x_i=x_{i-2}$ for some $1<i\leq k$. An {\em ordinary $k$-gon} is a 
$2k$-chain $(x_0,x_1,x_2,\ldots,x_{2k-1},x_0)$ with the property that 
$x_i\neq x_j$ for $0\leq i < j <2k$. If two
vertices $x,y$ can be joined by a $k$-chain, but not by any $k'$-chain for
$k'<k$, we say that the distance $d(x,y)$ is $k$. If there is no chain
joining $x$ and $y$, then we put $d(x,y)=\infty$. If the distance between
any pair of vertices if finite, the incidence structure $\frak P$ is called
{\em connected} (in the graph-theoretic sense). The {\em diameter} of an incidence structure is the supremum
of all the distances between vertices in $\V$. If $\P$ has finite diameter
$n$, then we call two vertices $x,y$ {\em opposite}, if $d(x,y)=n$.
We put $\opp x=\{ y\in \V|\ d(x,y)=n\}$. We also put
$\V^{(d\geq k)}=\{(x,y)\in\V\times\V|\ d(x,y)\geq k\}$, and
$\V^{(d=k)}=\{(x,y)\in\V\times\V|\ d(x,y)= k\}$.

\section{Galleries}

\label{GenPoly1_2}
A $(k+1)$-chain is called a {\em gallery} of length $k$. Given a flag 
$(p,\ell)$, let $\Gal_k(p,\ell)$ denote the set of all $(k+1)$-chains
of the form $(x_0=p,x_1=\ell,x_2,\ldots,x_{k+1})$, and let $\Gal_k(\ell,p)$ 
denote the set of all $(k+1)$-chains of the form 
$(x_0=\ell,x_1=p,x_2,\ldots,x_{k+1})$.

Put $\{u,v\}=\{p,\ell\}$. We let $\StamGal_k(u,v)$ denote the set of all 
stammering galleries in $\Gal_k(u,v)$. The non-stammering galleries 
$\PropGal_k(u,v)=\Gal_k(u,v)-\StamGal_k(u,v)$ are called
{\em proper galleries}. Note that $\PropGal_k(u,v)$ is nonempty, provided
that there are at least two points on every line and two lines through
every point.

The gallery ${\bf x}=(u,v,x_2,\ldots,x_{k+1})$ 
{\em ends at the flag $\flag(x_k,x_{k+1})\in\F$}.
We say also that the gallery $\bf x$ {\em ends at the vertex $x_{k+1}$.}

There is a canonical projection 
\begin{eqnarray*}
\pr : \Gal_{k+1}(u,v) & \to & \Gal_k(u,v) \\
(u,v,\ldots,x_{k+1},x_{k+2}) & \mapsto &
(u,v,\ldots,x_{k+1}),
\end{eqnarray*}
and a section 
\begin{eqnarray*}
s: \Gal_k(u,v) & \to & \StamGal_{k+1}(u,v) \\
(u,v,\ldots,x_k,x_{k+1}) & \mapsto & (u,v,\ldots,x_k,x_{k+1},x_k).
\end{eqnarray*}
Clearly $\pr\circ s=\id$.

There is also an injection
\begin{eqnarray*}
\Gal_k(v,u) & \to & \StamGal_{k+1}(u,v) \\
(v,u,x_2,\ldots,x_{k+1}) & \mapsto & (u,v,u,x_2,\ldots,x_{k+1}).
\end{eqnarray*}

\section{Automorphisms}

\label{GenPoly1_3}
An {\em automorphism} of an incidence structure is a bijection of
the set of vertices that maps points to points, lines to lines,
and flags to flags. It is clear that automorphisms map 
$k$-chains to $k$-chains, and that they preserve the distance $d$.

Homomorphisms between generalized polygons are investigated
in \cite{BK2}.

\section{Coset geometries}

\label{GenPoly1_4}
Suppose we are given a group $G$ with subgroups $A,B$. Then we may
form the {\em coset geometry} $(G/A,G/B,G/(A\cap B))$. The points are
the cosets of $A$, the lines are the cosets of $B$, and two cosets
are incident, if their intersection is nonempty (there is a slight
problem if $A=B$, but this case is not interesting anyway). The
flag space of this incidence structure can be identified with
$G/(A\cap B)$. The coset $g(A\cap B)$  represents the flag $(gA,gB)$.

Note that if $G$ acts as a flag-transitive automorphism group
on an incidence structure $\frak P$, and if $(p,\ell)$ is a flag,
then the coset geometry $(G/G_p,G/G_\ell,G/G_{p,\ell})$ is
isomorphic to $\frak P$. The existence of ordinary $k$-gons in 
$\frak P$ can be translated into a group-theoretic property, 
see \cite[3.1]{GKK}.

\section{Generalized polygons}

\label{GenPoly1_5}
A thick incidence structure $\frak P=\PLF$ is called a {\em generalized
$n$-gon}, if it satisfies the following two conditions:
\begin{enumerate}
\item $\frak P$ contains no ordinary $k$-gons for $2\leq k<n$.
\item Any two vertices $x,y$ are contained in an ordinary $n$-gon.
\end{enumerate}
It is readily verified that the generalized digons are exactly the
trivial thick incidence structures $(\P,\L,\P\times\L)$. So we assume
from now on that $n>2$.

An incidence structure is called a {\em partial $n$-gon}, if it satisfies
condition (i). This condition guarantees that for any two vertices
$x,y$ with $d(x,y)=k<n$, the $k$-chain $(x,x_1,\ldots,x_{k-1},y)$ joining
$x$ and $y$ is unique. So we may define a map $f_k(x,y)=x_{k-1}$ on the
set $\V_k^{(d=k)}$
of pairs of vertices at distance $k$. We are particularly interested 
in the map $f_{n-1}$.

Note that for every non-stammering $k$-chain $(x_0,\ldots,x_k)$, we
have $d(x_0,x_k)=k$, provided that $k\leq n$.

The following lemma will be required later.

\Lemma{\label{L1_5_1} (cp. \cite[3.30]{Ti1})
Let $x,y$ be vertices in a generalized $n$-gon. If $x$ and $y$ have
the same type, then there exists a vertex $z\in\opp x\cap\opp y$.
If $x$ and $y$ have different types, then there exists a vertex 
$z\in\opp x$ with $d(z,y)=n-1$.

\proof Let $z$ be a vertex in $\opp x$ that has maximal distance $k$ to $y$,
and suppose that $k\neq n,n-1$. Choose a vertex $b\in\V_z-\{f_k(y,z)\}$. 
We have $d(b,x)=n-1$ and $d(b,y)=k+1$, hence there is an element 
$z'\in\V_b-\{f_{n-1}(x,b),f_{k+1}(y,b)\}$. But now $d(z',y)=k+2$ and 
$d(z',x)=n$, a contradiction. (Here we needed the fact that the 
incidence structure $\frak P$ is thick.)\qed}

%
%

\section{Projectivities}

\label{GenPoly1_6}
Let $\frak P=\PLF$ be a generalized $n$-gon, and let $x,y$ be vertices
at distance $n$. We may define a map $[y,x]$ from $\V_x$ to $\V_y$
by $z\mapsto f_{n-1}(z,y)$, with inverse $[x,y]$. These maps are
called {\em perspectivities}. A concatenation of perspectivities is
called a {\em projectivity}, and we write $[z,y][y,x]=[z,y,x]$ etc.

\Lemma{Given any two vertices $x,y$ of the same type,
there exists a projectivity between $\V_x$ and $\V_y$; if $n$ is odd, then
such a projectivity exists also if $x$ and $y$ have different types.

\proof Let $x,y\in\V$ be vertices. By Lemma \ref{L1_5_1}, there exists
a vertex $z\in\opp x\cap\opp y$. Now $[y,z,x]$ is a projectivity
as required.

The last statement is clear, because opposite vertices have different 
types if $n$ is odd.\qed}

The set of projectivities of $\frak P$ forms a groupoid; the group $\Pi(x)$
of all projectivities that start and end at a given vertex $x$ is called 
the {\em group of projectivities} of $\frak P$ at $x$; by the lemma above, 
its equivalence type as a permutation group depends only on the type of $x$. It will be shown below in
\ref{D1_8_2} that $\V_x$ is 
doubly homogeneous under the action of $\Pi(x)$.

\section{Schubert cells}

\label{GenPoly1_7}
Let $(p,\ell)$ be a flag in a projective plane, and let $L$ be the
point row corresponding to $\ell$. Then the point space $\P$ can be
decomposed as $\{p\}\cup (L-\{p\}) \cup (\P-L)$. Now choose points
$q\in L-\{p\}$, $o\in \P-L$, and let $H$ be the point row corresponding
to $h=o\vee q$. We may introduce coordinates as follows:
\[
\begin{array}{lclclcl}
L-\{p\} & \to     & \L_o-\{o\vee p\} & {:} 
    & x & \mapsto & x\vee o \\
\P-L    & \to     & (H-\{q\})\times(\L_q-\{p\vee q\}) & {:} 
    & x & \mapsto & (x\vee q,(p\vee x)\wedge h).
\end{array}
\]
There is a similar decomposition of the flag space and of the line space.

We want to generalize these concepts to arbitrary polygons.

\Def{Let $\frak P$ be a generalized $n$-gon. Let $(p,\ell)$ be a
flag, and let $\{u,v\}=\{p,\ell\}$. For $0\leq k\leq n$ we put
\[
\F_k(u,v)=\{ (q,h)\in\F|\ 
\mbox{there is an }{\bf x}\in
\PropGal_k(u,v)
\mbox{ that ends at }(q,h)\}.
\]
The set $\F_k(u,v)$ is called a {\em Schubert cell}.
For $0<k<n$, Schubert cells of different types have no flags in
common, whereas $\F_k(u,v)=\F_k(v,u)$ for $k=0,n$. Hence there are $2n$
Schubert cells that cover the flag space. If we choose any ordinary
$n$-gon $\bf x$ containing $(p,\ell)$, then the Schubert cells are
precisely the preimages of the {\em retraction}
$(\F,(p,\ell)) \to ({\bf x},(p,\ell))$, see \cite[3.3]{Ti1}.

The sets 
\[
\Cl\F_k(u,v)=\F_k(u,v)\cup\bigcup_{j<k}(\F_j(u,v)\cup\F_j(v,u))
\]
are called the 
{\em closed Schubert cells} or {\em Schubert varieties}. It is immediate 
that the Schubert variety $\Cl\F_k(u,v)$ consists precisely of all the flags 
that occur as ends of galleries in $\Gal_k(u,v)$.

Similarly, we put
\[
\begin{array}{lcl}
\P_{2k}(\ell,p) & = & \{ q\in\P|\ 
\mbox{there is an }{\bf x}\in
\PropGal_{2k}(\ell,p)
\mbox{ that ends at }q\} \\
\P_{2k+1}(p,\ell) & = & \{ q\in\P|\ 
\mbox{there is an }{\bf x}\in
\PropGal_{2k+1}(p,\ell)
\mbox{ that ends at }q\} 
\end{array}
\] 
for $0\leq 2k<n$ and $1\leq 2k+1<n$, respectively.
In the same way we define
\[
\begin{array}{lcll}
\L_{2k}(p,\ell) & = & \{ h\in\L|\ 
\mbox{there is an }{\bf x}\in
\PropGal_{2k}(p,\ell)
\mbox{ that ends at }h\} \\
\L_{2k+1}(\ell,p) & = & \{ h\in\L|\ 
\mbox{there is an }{\bf x}\in
\PropGal_{2k+1}(\ell,p)
\mbox{ that ends at }h\} 
\end{array}
\]
for $0\leq 2k <n$ and $ 1\leq 2k+1 <n$, respectively.
The $2n$ sets of this type are also called Schubert cells; they form a 
partition of $\P$ and $\L$, respectively. We also introduce Schubert
varieties by 
\[
\Cl\P_k(u,v)=\P_0(\ell,p)\cup\P_1(p,\ell)\cup\P_2(\ell,p)\cup\ldots
\cup\P_k(u,v),
\]
and 
\[
\Cl\L_k(v,u)=\L_0(p,\ell)\cup\L_1(\ell,p)\cup\L_2(p,\ell)\cup\ldots
\cup\L_k(v,u);
\]
these are the vertices that can be reached by 
(possibly stammering) galleries of the corresponding type.
Note also that 
\[
\Cl\P_k(u,v)=\{p\in\P|\ d(v,p)\leq k\}
\]
is precisely the set 
of all points that have distance $\leq k$ to $v$ and thus does not depend on
$u$. Similarly,
\[
\Cl\L_k(v,u)=\{\ell\in\L|\ d(u,\ell)\leq k\}.
\]

This decomposition is called the {\em Schubert cell decomposition} 
of the $n$-gon $\frak P$ (with respect to the flag $(p,\ell)$).
The sets $\P_{n-1}(u,v)$, $\L_{n-1}(v,u)$, and
$\F_n(u,v)$ are called the {\em big cells}.}

Returning to the example of the projective plane, we have
\begin{eqnarray*}
\P_0(\ell,p) & = & \{p\}   \\
\P_1(p,\ell) & = & L-\{p\} \\
\P_2(\ell,p) & = & \P-L.
\end{eqnarray*}

In order to introduce coordinates in the generalized $n$-gon $\frak P$,
we choose an ordinary $n$-gon 
${\bf x}=(x_0=p,x_1=\ell,x_2,\ldots,x_{2n-1},x_0)$.

Let $q\in\P_{2k+1}(p,\ell)$, and let ${\bf y}=(p,\ell,y_2,\ldots,y_{2k+2}=q)$ 
be the corresponding gallery. Note that $d(y_j,x_{n+j})=n$, because the
$n$-chain $(y_j,y_{j-1},\ldots,y_0,x_{2n-1},\ldots,
\linebreak
x_{n+j})$ does not 
stammer. Hence we have the relation
\[
y_{j-1}=f_{n-1}(x_{n+j+1},y_j),
\]
and therefore, the bijection $\P_{2k+1}(p,\ell)\to\PropGal_{2k+1}(p,\ell)$ 
can be expressed in terms of the map $f_{n-1}$ and the ordinary $n$-gon
$\bf x$. 

Now we attach to $q$ the $2k+1$ coordinates
\[
q_{j-1}=f_{n-1}(y_j,x_{n+j-1}) \in\V_{x_{n+j-1}}-\{x_{n+j}\},
\]
where $2\leq j\leq 2k+2$. The point $q$ can be recovered from these 
coordinates in terms of the function $f_{n-1}$, because
\[
y_j=f_{n-1}(q_{j-1},y_{j-1}).
\]
\newpage

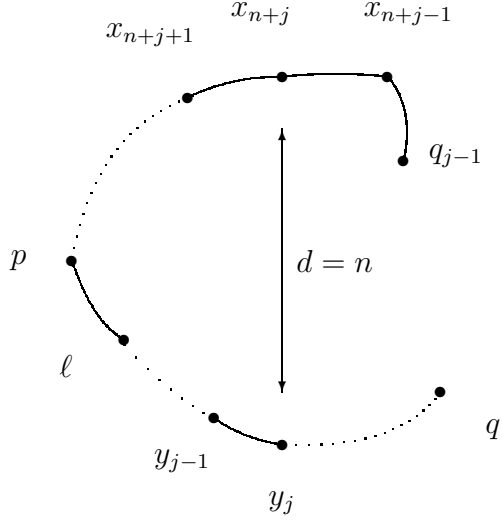
\[
\unitlength=0.70mm
\linethickness{0.4pt}
\begin{picture}(110.00,95.00)(0,15)
\put(70.00,30.00){\circle*{2.00}}
\put(70.00,100.00){\circle*{2.00}}
\put(30.00,65.00){\circle*{2.00}}
\put(40.00,50.00){\circle*{2.00}}
\put(100.00,40.00){\circle*{2.00}}
\put(90.00,100.00){\circle*{2.00}}
\put(70.00,60.00){\vector(0,1){30.00}}
\put(70.00,60.00){\vector(0,-1){20.00}}
\put(80.00,65.00){\makebox(0,0)[cc]{$d=n$}}
\put(20.00,65.00){\makebox(0,0)[cc]{$p$}}
\put(29.00,45.00){\makebox(0,0)[cc]{$\ell$}}
\put(70.00,19.00){\makebox(0,0)[cc]{$y_j$}}
\put(110.00,33.00){\makebox(0,0)[cc]{$q$}}
\put(66.00,112.00){\makebox(0,0)[cc]{$x_{n+j}$}}
\put(94.00,112.00){\makebox(0,0)[cc]{$x_{n+j-1}$}}
\put(57.00,35.00){\circle*{2.00}}
\put(51.00,27.00){\makebox(0,0)[cc]{$y_{j-1}$}}
\put(93.00,84.00){\circle*{2.00}}
\put(103.00,85.00){\makebox(0,0)[cc]{$q_{j-1}$}}
\put(52.00,96.00){\circle*{2.00}}
\put(45.00,108.00){\makebox(0,0)[cc]{$x_{n+j+1}$}}
\bezier{76}(30.00,65.00)(34.00,53.00)(40.00,50.00)
\bezier{56}(57.00,35.00)(63.00,31.00)(70.00,30.00)
\bezier{76}(52.00,96.00)(60.00,100.00)(70.00,100.00)
\bezier{80}(70.00,100.00)(80.00,101.00)(90.00,100.00)
\bezier{72}(90.00,100.00)(95.00,94.00)(93.00,84.00)
\bezier{20}(30.00,65.00)(34.00,89.00)(52.00,96.00)
\bezier{8}(40.00,49.00)(49.00,39.00)(57.00,35.00)
\bezier{17}(70.00,30.00)(92.00,29.00)(100.00,40.00)
\end{picture}
\]

Thus there is a bijection 
\[
\P_{2k+1}(p,\ell)\to
(\V_{x_{n+1}}-\{x_{n+2}\})\times \cdots \times
(\V_{x_{n+2k+1}}-\{x_{n+2k+2}\})
\]
that can in both directions be expressed in terms of the function $f_{n-1}$.

Observe also that $\pr_1$ maps $\F_{2k+1}(p,\ell)$ bijectively onto
$\P_{2k+1}(p,\ell)$ for $2k+1<n$, with inverse
$q\mapsto (q,f_{n-1}(x_{n+2k+3},q))$.

The other Schubert cells can be treated in a similar way, and we get the
following  result:

\Prop{\label{P1_7_2}
Let $(p,\ell)$ be a flag in the generalized $n$-gon $\frak P$, and let 
${\bf x}=(x_0=p,x_1=\ell,x_2,\ldots,x_{2n-1},x_0)$ be an ordinary $n$-gon.
The maps 
\[
\begin{array}{lclcll}
\PropGal_{2k}(\ell,p)   & \to        & \F_{2k}(\ell,p) 
                        & \TO{\pr_1} & \P_{2k}(\ell,p) 
& (2k<n) \\
\PropGal_{2k+1}(p,\ell) & \to        & \F_{2k+1}(p,\ell) 
                        & \TO{\pr_1} & \P_{2k+1}(p,\ell) 
& (2k+1<n) \\
\PropGal_{2k}(p,\ell)   & \to        & \F_{2k}(p,\ell) 
                        & \TO{\pr_2} & \L_{2k}(p,\ell) 
& (2k<n) \\
\PropGal_{2k+1}(\ell,p) & \to        & \F_{2k+1}(\ell,p) 
                        & \TO{\pr_2} & \L_{2k+1}(\ell,p) 
& (2k+1<n) \\
\PropGal_n(p,\ell)      & \to        & \F_n(p,\ell)  \\
\PropGal_n(\ell,p)      & \to        & \F_n(\ell,p) 
\end{array}
\]
are bijections, and their inverses can be expressed in terms of
the function $f_{n-1}$ and the ordinary $n$-gon $\bf x$.

Moreover, there are bijections
\[
\begin{array}{lcl}
\P_{2k}(\ell,p) & \to & 
  (\V_{x_n}-\{x_{n-1}\})\times \cdots \times 
  (\V_{x_{n+2k+1}}-\{x_{n+2k}\}) \\
\P_{2k+1}(p,\ell)   & \to & 
  (\V_{x_{n+1}}-\{x_{n+2}\})\times \cdots \times 
  (\V_{x_{n+2k+1}}-\{x_{n+2k+2}\}) \\
\L_{2k}(p,\ell)   & \to & 
  (\V_{x_{n+1}}-\{x_{n+2}\})\times \cdots \times 
  (\V_{x_{n+2k}}-\{x_{n+2k+1}\}) \\
\L_{2k+1}(\ell,p) & \to & 
  (\V_{x_n}-\{x_{n-1}\})\times \cdots \times
  (\V_{x_{n-2k}}-\{x_{n-2k-1}\})  \\
\F_k(p,\ell)      & \to & 
  (\V_{x_{n+1}}-\{x_{n+2}\})\times \cdots \times 
  (\V_{x_{n+k-1}}-\{x_{n+k}\}) \\
\F_k(\ell,p)      & \to & 
  (\V_{x_n}-\{x_{n-1}\})\times \cdots \times
  (\V_{x_{n-k}}-\{x_{n-k-1}\})    .
\end{array}
\]
that can (in both directions) be expressed in terms of the function
$f_{n-1}$ and the ordinary $n$-gon $\bf x$.

All of these maps are {\em natural with respect to automorphisms}, i.e.
it does not matter if the automorphism is applied to the coordinates,
or if the coordinates are taken with respect to the image of the ordinary 
$n$-gon $\bf x$ under the automorphism. In particular, an automorphism is
completely determined by its action on the sets 
$\V_{x_0},\ldots,\V_{x_{2n-1}}$.\qed}

The following fact will be required later.

\Prop{\label{Separate}
Let $q\in\P_{n-1}(u,v)$. Let $C\SUB\P$ be a Schubert cell
(with respect to $(u,v)$) different from $\P_{n-1}(u,v)$. 

There exists an
incident pair $(u',v')$ with the following properties:

\begin{enumerate}
\item $q\in\P_{n-1}(u',v')$
\item $C\SUB\P_{n-1}(u',v')$
\item The map $x\mapsto f_{n-1}(x,v')$ has constant value $z$ on 
$C$, and $f_{n-1}(q,v')\neq z$.
\end{enumerate}

\proof We treat only case that $n-k$ is odd (the case that $n-k$ is even is similar); so $C=\P_{k}(u,v)$
for some $k<n-2$.

Inductively, we can find a non-stammering chain 
$(x_0=v,x_1=u,x_2,\ldots,x_{n-k})$ with the property that 
$d(x_j,q)\geq n-1$ for $0\leq j\leq n-k$. Thus for
$o\in\P_k(u,v)$, we have $d(o,x_{n-k-1})=d(q,x_{n-k-1})=n-1$, and 
$f_{n-1}(o,x_{n-k-1})=x_{n-k-2}$.
On the other hand, $f_{n-1}(q,x_{n-k-1})\neq x_{n-k-2}$, because
$d(x_{n-k-2},q)=n$. Thus we may put $(u',v')=(x_{n-k},x_{n-k-1})$.\qed}

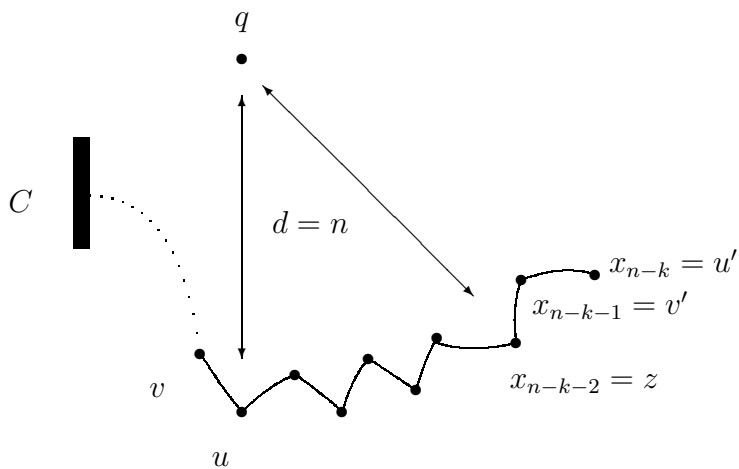
\[
\unitlength=0.70mm
\linethickness{0.4pt}
\begin{picture}(125.00,80.00)(0,40)
\put(60.00,117.00){\circle*{2.00}}
\put(60.00,50.00){\circle*{2.00}}
\put(52.00,61.00){\circle*{2.00}}
\put(60.00,80.00){\vector(0,1){30.00}}
\put(60.00,80.00){\vector(0,-1){20.00}}
\put(84.00,92.00){\vector(1,-1){20.00}}
\put(84.00,92.00){\vector(-1,1){20.00}}
\put(112.00,63.00){\circle*{2.00}}
\put(97.00,64.00){\circle*{2.00}}
\put(93.00,54.00){\circle*{2.00}}
\put(84.00,60.00){\circle*{2.00}}
\put(79.00,50.00){\circle*{2.00}}
\put(70.00,57.00){\circle*{2.00}}
\put(28.00,81.00){\rule{3.00\unitlength}{21.00\unitlength}}
\put(18.00,90.00){\makebox(0,0)[cc]{$C$}}
\put(60.00,124.00){\makebox(0,0)[cc]{$q$}}
\put(44.00,54.00){\makebox(0,0)[cc]{$v$}}
\put(56.00,41.00){\makebox(0,0)[cc]{$u$}}
\put(125.00,55.00){\makebox(0,0)[cc]{$x_{n-k-2}=z$}}
\put(113.00,75.00){\circle*{2.00}}
\put(130.00,70.00){\makebox(0,0)[cc]{$x_{n-k-1}=v'$}}
\put(127.00,76.00){\circle*{2.00}}
\bezier{60}(113.00,75.00)(121.00,78.00)(127.00,76.00)
\put(142.00,78.00){\makebox(0,0)[cc]{$x_{n-k}=u'$}}
\bezier{15}(31.00,91.00)(48.00,91.00)(52.00,61.00)
\bezier{56}(52.00,61.00)(57.00,53.00)(60.00,50.00)
\bezier{48}(60.00,50.00)(65.00,55.00)(70.00,57.00)
\bezier{44}(79.00,50.00)(73.00,55.00)(70.00,57.00)
\bezier{48}(79.00,50.00)(81.00,57.00)(84.00,60.00)
\bezier{44}(84.00,60.00)(90.00,56.00)(93.00,54.00)
\bezier{44}(93.00,54.00)(95.00,61.00)(97.00,64.00)
\bezier{64}(97.00,63.00)(103.00,61.00)(112.00,63.00)
\bezier{48}(112.00,63.00)(111.00,70.00)(113.00,75.00)
\put(73.00,86.00){\makebox(0,0)[cc]{$d=n$}}
\end{picture}
\]

\section{The algebraic operations}

Let $\frak P=\PLF$ be a generalized $n$-gon, and let 
$(x_0,x_1,\ldots,x_{2n-1},x_0)$ be an ordinary $n$-gon in $\frak P$.
Let $a\in\V_{x_{n-1}}-\{x_{n-2},x_n\}$. Since $d(a,x_0)=d(a,x_{2n-2})=n$,
we may define a projectivity $\mu_a=[x_{2n-2},a,x_0]$. This projectivity
maps $x_1$ to $x_{2n-3}$ and fixes $x_{2n-1}$. Choose an element
$1_L\in\V_{x_{n-1}}-\{x_{n-2},x_n\}$, and consider the map
\[
\mu_a^{-1}\mu_{1_L}(x)=f_{n-1}(f_{n-1}(f_{n-1}(f_{n-1}(x,1_L),x_{2n-2}),a),
x_0).
\]
The right-hand side makes sense also for $a=x_{n-2}$, provided that 
$x\neq x_{2n-1}$, and 
$f_{n-1}(f_{n-1}(f_{n-1}(f_{n-1}(x,1_L),x_{2n-2}),x_{n-2}),x_0)=x_1$.
Similarly, if $x\neq x_1$, then
$f_{n-1}(f_{n-1}(f_{n-1}(f_{n-1}(x,1_L),x_{2n-2}),x_{n-2}),x_n)=x_{2n-1}$.

This leads to the following definition: 

\Def{\label{D1_8_1}(cp. \cite[1.1]{GKK})
Let $(x_0,x_1,\ldots,x_{2n-1},x_0)$ be an ordinary $n$-gon
in the generalized $n$-gon $\frak P$. Put $0_K=x_1$,
$\infty_K=x_{2n-1}$, $0_L=x_{n-2}$, and $\infty_L=x_n$. Set
$K=\V_{x_0}-\{\infty_K\}$, $L=\V_{x_{n-1}}-\{\infty_L\}$, and choose
an element $1_L\in L-\{0_L\}$. We may define maps 
\begin{eqnarray*}
K\times L & \to     & K \\
(x,y)     & \mapsto & x\bullet y=f_{n-1}(f_{n-1}(f_{n-1}(f_{n-1}(x,1_L),
x_{2n-2}),y),x_0) \\ 
K\times (L-\{0_L\}) & \to     & K \\
(x,y)             & \mapsto & x/y=f_{n-1}(f_{n-1}(f_{n-1}(f_{n-1}(x,y),
x_{2n-2}),1_L),x_0)
\end{eqnarray*}
These maps have the following properties:
\begin{enumerate}
\item The relations $x\bullet 0_L=0_K\bullet y=0_K$ and
$x\bullet 1_L=x$ hold for all $(x,y)\in K\times L$.
\item The relation $(x\bullet y)/y=(x/y)\bullet y=x$ holds for all
$(x,y)\in K\times L-\{0_L\}$.
\end{enumerate}
The map $\bullet$ is called the {\em multiplication} with respect to
$(x_0,x_1,\ldots,x_{2n-1},1_L)$.}

We shall also need the map 
\[
\begin{array}{ccc}
\rho_b & : & L\cup\{\infty_L\} \to K\cup\{\infty_K\}, \\
  & & y\mapsto f_{n-1}(f_{n-1}(f_{n-1}(f_{n-1}(b,1_L),x_{2n-2}),y),x_0),
\end{array}
\]
which is defined for $b\in K-\{0_k\}$. Note that $\rho_b(0_L)=0_K$,
and $\rho_b(1_L)=b$.

\medskip
Next, choose an element $e\in\V_{x_{n+1}}-\{x_n,x_{n+2}\}$. For $y\in K$,
consider the projectivity $\pi_y=[x_0,x_n,f_{n-1}(e,y),x_{n+2},
f_{n-1}(e,x_1),x_n,x_0]\in\Pi(x_0)$. Note that $\pi_y(\infty_K)=\infty_K$,
$\pi_y(0_K)=y$, and $\pi_{0_K}=\id_{\V_{x_0}}$.
This leads to the following definition.

\Def{\label{D1_8_2}(cp. \cite[1.4]{GKK})
For $x,y\in K$, we define maps $\pm:K\times K\to K$ by
$x+y=\pi_y(x)$ and $x-y=\pi_y^{-1}(x)$. These maps have the
following properties:
\begin{enumerate}
\item $x+0_K=0_K+x=x$ for all $x\in K$.
\item $(x+y)-y=(x-y)+y=x$ for all $x,y\in K$.
\end{enumerate}
Thus $(K,0_K,+,-)$ is a {\em right loop} in the sense of \cite{BK1}.

In particular, the group of all projectivities is doubly-transitive on
$K\cup\{\infty_K\}$.}

\chapter{Topological polygons}

This chapter deals with the set-theoretic topology of
topological generalized polygons. The definition of a topological
polygon is rather general; we require only that the geometric operation
$f_{n-1}$ be continuous on its domain, and that there exists at least
{\em one} non-trivial open set in $\P\cup\L$ (\ref{TopPolyDef}).

The algebraic operations defined on the point rows and the pencils
of lines imply that every Schubert cell is homogeneous and regular 
(\ref{RegularCell}).
Moreover, the big cells are open (\ref{BigCellsOpen}),
hence the point space, the line
space, and the flag space are Hausdorff spaces (in fact even
regular, as J\"ager proved (\ref{Regular})).
The flag space $\F$ is closed in the product $\P\times\L$ (\ref{FlagsClosed}),
and the maps $\F\to\P$, $\F\to\L$ are locally trivial bundles
(\ref{F2PBundle}).

A topological polygon is either connected or totally disconnected
(\ref{P2_2_3}).
If it is path connected, then it is locally path-connected and locally
contractible (\ref{LocContract}).
The point rows are discrete if and only if the pencils
of lines are discrete (\ref{P2_2_6}).
The Schubert varieties are the topological closures of the corresponding
Schubert cells, provided that the topology is not discrete
(\ref{SchubertClosed}). In particular,
the point rows and the sets $p^\perp$ are closed.

If the topology of the polygon is locally compact, then it is second
countable and metrizable \cite{GKK} (\ref{T2_4_3}).
If the polygon is in addition connected, then it is path-connected (and thus
locally contractible) and compact. 

\medskip
The topological triangles are precisely the topological projective
planes, see \cite{Sal57,Sal67}. For topological quadrangles see
also \cite{For81}, \cite{GK}, \cite{Sch92}.

\setcounter{nonsense}{0}
\bigskip
\Def{\label{TopPolyDef}
A generalized $n$-gon $\frak P=\PLF$ is called a {\em topological
$n$-gon}, if the point space $\P$ and the line space $\L$ carry topologies
such that the map $f_{n-1}$ is continuous on its domain
$\V^{(d=n-1)}=\{(x,y)\in\V\times\V|\ d(x,y)=n-1\}$ (we
endow the set of vertices $\V=\P\cup\L$ with the sum topology). 
 
{\em In order to avoid trivialities, we shall always assume that there exists 
an open set in $\V$ besides $\emptyset$, $\P$, $\L$, and $\V$.}}

The following observations are simple, but important:

\Prop{Ever projectivity is a homeomorphism. All point rows are homeomorphic,
doubly homogeneous spaces, and similarly, all pencils of lines are 
homeomorphic, doubly homogeneous spaces. If $n$ is odd, then all point rows 
and all pencils of lines are homeomorphic via projectivities.

The coordinate functions defined in \ref{P1_7_2} are homeomorphisms. The
algebraic operations $\bullet$, $/$, and $\pm$ defined in \ref{D1_8_1},
\ref{D1_8_2} are continuous.\qed}

\Prop{(\cite{Hofm})
Let $(G,0,\pm)$ be a right loop in the sense of \cite{BK1}, see
\ref{D1_8_2}. Suppose that $G$ is a topological
${\bf T}_1$-space, and that the maps $\pm$ are continuous. Then $G$ is
a regular (i.e. a ${\bf T}_3$-) space.

\proof Let $U$ be a neighborhood of $0$. Choose a
neigborhood $V$ of $0$
with $V+V\SUB U$. We claim that $\overline V\SUB V+V$. For every element
$x\in\overline V$, there is a net $(x_\nu)\SUB V$ converging to $x$. Hence 
the net $(x-x_\nu)$ converges to $0$. Thus there is an element $y\in V$ with
$x-y\in V$, and therefore we have $x=(x-y)+y\in V+V$. 
 
Since $G$ is homogeneous, it is a regular space.\qed}

\section{The topology of Schubert cells}

Let $\frak P=\PLF$ be a topological $n$-gon.

\Lemma{Every point row and every pencil of lines is a ${\bf T}_1$-space.

\proof We may assume that the point space $\P$ contains a proper,
open subset $U$, with $\emptyset\neq U\neq\P$. Thus there exists a point row
$L$ with $L-U\neq L\neq L\cap U$. Since $L$ is doubly homogeneous, 
every point is closed in $L$. 

If $n$ is odd, then all point rows and all pencils of lines are
homeomorphic via projectivities, so we are done.

If $n$ is even, then the map $\rho_b$ of \ref{D1_8_1} provides a nonconstant
map of some pencil of lines $\L_p$ into some point row. Thus $\L_p$
contains proper, nonempty open subsets, and the claim follows
for the pencils of lines by duality.\qed}

\Cor{\label{RegularCell}
Every Schubert cell is homogeneous and regular.\qed}

\Thm{\label{BigCellsOpen}
Put $(u,v)=(p,\ell)$ if $n$ is even, and $(u,v)=(\ell,p)$ if
$n$ is odd. The big cells $\P_{n-1}(u,v)$, $\L_{n-1}(v,u)$, and
$\F_n(u,v)$ are open.

\proof Let $q\in\P_{n-1}(u,v)$, and let $C\SUB\P$ be a Schubert cell
different from $\P_{n-1}(u,v)$. By \ref{Separate} there exists
an incident pair $(u',v')$ with $C\cup\{q\}\SUB\P_{n-1}(u',v')$, 
and such that the map $x\mapsto f_{n-1}(x,v')$ has constant value 
$z\neq f_{n-1}(q,v')$ on $C$. Therefore there exists an open
neighborhood $U$ of $q$ that does not meet $C$. We may do this
for each of the $n-1$ Schubert cells in $\P$
different from $\P_{n-1}(u,v)$;
taking the intersection of the corresponding $n-1$ open neighborhoods
of $q$, we obtain an open neighborhood of $q$ contained in
$\P_{n-1}(u,v)$.

The same applies to $\L_{n-1}(v,u)$, and finally 
\[
\F_n(u,v)=\pr_1^{-1}(\P_{n-1}(u,v))\cap\pr_2^{-1}(\L_{n-1}(v,u)).
\]
\qed}

\Cor{The point space $\P$, the line space $\L$, and the flag space
$\F$ are Hausdorff spaces.\qed}

\Lemma{Let $x$ be a vertex, and let $\V_x^{(k)}$ denote the set
of all vertices that have distance $k$ to $x$. Then for $k<n$, the
minimal $k$-chain $(x=x_0,x_1,\ldots,x_{k-1},x_k=y)$ that joins 
$y\in\V_x^{(k)}$ and $x$ depends continuously on $y$.

\proof Clearly, it suffices to show that $x_{k-1}$ depends continuously
on $y$. So let $y\in\V_x^{(k)}$. We may choose a non-stammering chain
$(x_{-n+k},\ldots,x_0,\ldots,x_k)$. Thus $y$ is contained in the open
cell determined by the pair $(x_{-n+k},x_{-n+k+1})$, and 
$x_{k-1}=f_{n-1}(x_{-n+k+1},x_k)$.\qed}

\Lemma{\label{L2_1_6}
Let $p\in\P$ be a point and consider the map 
$p^\perp-\{p\}\to \L_p$, $q\mapsto f_2(p,q)$. This defines a locally 
trivial bundle with the punctured point rows through $p$
as fibers. Let $x$ be
a vertex opposite to $p$. Then the map $\ell\mapsto f_{n-1}(x,\ell)$
is a section of this bundle. In particular, there is an imbedding
$\L_p\imb p^\perp$. 

\proof Let $\ell\in\L_p$ be a line, and put $y=f_{n-1}(\ell,x)$,
$z=f_{n-1}(p,u)$.
The map $\L_p-\{\ell\}\times\V_y-\{z\}\to p^\perp-\{p\}$,
$(h,w)\mapsto f_{n-1}(w,h)$ is a trivialization of
$p^\perp-\{p\}\to\L_p$ over the open set $\L_p-\{\ell\}$.\qed}

This may be put in a slightly more general form:

\Lemma{\label{SchubertBundle}
Let $p$ be a point, and let $\V_p^{(k)}$ denote the set of all
vertices $x$ with $d(x,p)=k$. For $k<n$, let $x_j$ be the $j$-th 
vertex in the minimal $k$-chain from $x_k\in\V_p^{(k)}$ to $p$.
The map $x_k\mapsto x_j$ defines for $0<j<k<n$
a locally trivial bundle whose fibers are Schubert cells.

For example, we get bundles $p^\perp-\{p\}\to\L_p$ (for $n\geq 3$),
and $\P-L\to L$ (for $n=4$).

\proof This is clear form the coordinatization with respect to
$(p,\ell)$, where $\ell$ is a line through $p$.\qed}

\Prop{\label{F2PBundle}
The maps $\pr_1:\F\to\P$ and $\pr_2:\F\to\L$ are locally trivial
bundles, and hence open.

\proof Let $p\in\P$ be a point. Choose a big cell $\P_{n-1}(u,v)$
containing $p$. For every $q\in\P_{n-1}(u,v)$, we have
$d(u,q)=n$, thus we may consider the map 
$\V_u\times\P_{n-1}(u,v)\to\F$,
$(q,z)\mapsto (q,f_{n-1}(z,q))$.
This is clearly a trivialization
of $\pr_1$ over the open set $\P_{n-1}(u,v)$.\qed}

\Cor{If $U$ is an open set of points, then the set $\L_U$ of all lines that
meet $U$ is open (and dually).
\qed}

\Cor{The map $\pr:\Gal_k(u,v)\to\Gal_{k-1}(u,v)$ is a locally trivial
bundle, with a section $s$ as defined in \ref{GenPoly1_2}.

\proof This bundle is just the pullback of one of the bundles
$\F\to\P$ or $\F\to\L$ by the map that sends a gallery to the vertex
where it ends.\qed}

\Prop{The Schubert varieties are closed subspaces in $\P$, $\L$, and
$\F$, respectively. In particular, every point row, every pencil of lines,
and every star $x^\perp$ is closed.

\proof The Schubert variety $\Cl\P_{n-2}(v,u)$ is closed, because it is the
complement of the big cell $\P_{n-1}(u,v)$. We proceed by induction on
$n-k$ and may assume that $\Cl\L_{k+1}(u,v)$ is closed. Put 
$U=\L-\Cl\L_{k+1}(u,v)$. Now $\P_U=\P-\Cl\P_k(u,v)$ is open.\qed}

\Prop{\label{FlagsClosed}
The flag space $\F$ is closed in $\P\times\L$.

\proof Let $(p_\nu,\ell_\nu)$ be a net of flags, converging to $(p,\ell)$.
By Lemma \ref{L1_5_1} there is a vertex $u\in \opp p$ with $d(u,\ell)=n-1$.
Choose a vertex $v\in\V_u$. Then the point $p$ is contained in the big cell
$\P_{n-1}(u,v)$, therefore we may assume that all points $p_\nu$ are
contained in $\P_{n-1}(u,v)$, and thus $d(u,\ell_\nu)=d(u,\ell)=n-1$.
Now we have the relation $\ell_\nu=f_{n-1}(f_{n-1}(\ell_\nu,u),p_\nu)$.
Passing to the limit we obtain $\ell=f_{n-1}(f_{n-1}(\ell,u),p)\in
\L_p$.\qed}

The following theorem is due to J\"ager.

\Thm{\label{Regular}{\bf(J\"ager)} \cite{Jag}
Every Schubert variety is regular.

\proof It suffices to show that $\P$ and $\L$ are regular.
Let $q\in\P_{n-1}(u,v)$, and let $A\SUB\P-\{x\}$ be closed.
For every Schubert cell $C\SUB\P$, there exists an open set
$U_C$ containing $A\cap C$, and an open set $W_C$ containing $q$ with
$U_C\cap W_C=\emptyset$ by
\ref{Separate} and \ref{RegularCell}. Thus, the union of the
open sets $U_C$ and the intersection of the open sets $W_C$ are
disjoint open neighborhoods of $A$ and $q$, respectively.\qed}

\section{Connectivity properties}

\Lemma{\label{L2_2_1}
If the (path-) component of some point $p$ in a point row $L$ 
is non-trivial, then the point row $L$ is (path-) connected.

\proof This is clear from the fact that $L$ is doubly homogeneous.\qed}

\Prop{\label{P2_2_2}
Let $\frak P$ be a topological $n$-gon. Then the following
assertions hold:

\begin{enumerate}
\item If the point rows are (path-) connected, then the pencils of lines are 
(path-) connected. 
\item If the point rows are discrete, then the pencils of lines are 
discrete.
\end{enumerate}

\proof This is clear if $n$ is odd.

If $n$ is even, and if the point rows
are (path-) connected, then the map $\rho_b$ of \ref{D1_8_1} provides
a nonconstant map from some point row into some pencil of lines.
Therefore the pencils of lines are (path-) connected by \ref{L2_2_1}.

Next, assume that the point rows are not discrete. 
Consider the multiplication $K\times L\to K$, where 
$H=L\cup \{ \infty_L \}$ is a point row, and $K\cup\{\infty_K\}$ is
a pencil of lines. Since $H$ is homogeneous
and not discrete, there is a net $(x_\nu)\SUB L-\{0_L\}$ converging
to $0_L$. Choose $b\in K-\{0_K\}$. Now $y_\nu=b\bullet x_\nu$ converges
to $b\bullet 0_L=0_K$, but $b\bullet x_\nu \neq 0_K$, hence $0_K$
is not isolated in $K$.\qed}

\Prop{\label{P2_2_3}
The following are equivalent:
\begin{enumerate}
\item Some point row is (path-) connected.
\item Some pencil of lines is (path-) connected.
\item The point space is (path-) connected.
\item The line space is (path-) connected.
\item The flag space is (path-) connected.
\end{enumerate}
If one of these equivalent conditions is violated, then the 
(path-) component of every point, line, or flag is trivial.

\proof The conditions (i), (ii) are equivalent by \ref{P2_2_2}.
If the point
rows and the pencils of lines are (path-) connected, then the flag space is
(path-) connected, because for any flag $(p,\ell)$, the space
$(\{p\}\times L)\cup(\L_p\times\{\ell\})$ of all flags that have a point
or a line in common with $(p,\ell)$ is (path-) connected. Thus $\P$ and
$\L$ are (path-) connected. 

Now suppose that the point rows and the pencils of lines are not (path-)
connected. Then the (path-) component of every point in a point row
is trivial, and thus the (path-) component of every member of a Schubert
cell is trivial. Thus, if the point rows are not connected, then
the big cells are totally disconnected. Since any two points (any two lines,
any two flags) are contained in a big cell, the spaces $\P$, $\L$ and $\F$
are totally disconnected. 

If the point rows and pencils of lines are not path-connected, then
every map of the unit interval into a Schubert cell is constant;
hence the path-component of every point in the point space is trivial.\qed}

\Def{We call a topological polygon {\em (path-) connected}, if it satisfies
one of the five equivalent conditions of \ref{P2_2_3}. If it is
not connected, we call it {\em totally disconnected.}}

\Prop{\label{LocContract}
Let $\frak P$ be a path-connected polygon. Then $\P$, $\L$, $\F$,
every point row, and every pencil of lines is locally contractible.
Every Schubert cell is pseudo-isotopically contractible.

\proof Consider the multiplication defined in \ref{D1_8_1}.
Let $\gamma :\II \to L$ be a path with $\gamma(0)=0_L$ and $\gamma(1)=1_L$.
Let $U$ be a neighborhood of $0_K$. Every element $(0_K,t)\in U\times\II$
has a neighborhood $V_t\times I_t\SUB U\times\II$ with 
$V_t\bullet\gamma(I_t)\SUB U$. Since the unit interval $\II$ is compact, 
there is a finite subcovering $I_{t_1},\ldots,I_{t_k}$ of $\II$. We put
$V=\bigcap_{i=1}^k V_{t_i}$. Now $V\bullet\gamma(\II)\SUB U$, hence
$V$ can be contracted to $0_K$ within $U$.

Thus every point row and every pencil of lines is locally contractible.
This implies that every Schubert cell is locally contractible, and
the claim follows.

The space $K$ is pseudo-isotopically contractible by means of the map
$(x,t)\mapsto x\bullet\gamma(t)$ (we may assume that $\gamma(t)\neq 0_K$
for $t>0$).\qed}

\Prop{\label{P2_2_6}
Let $\frak P$ be a topological polygon.
The following are equivalent:
\begin{enumerate}
\item Some point row is discrete.
\item Some pencil of lines is discrete.
\item The point space is discrete.
\item The line space is discrete.
\item The flag space is discrete.
\item The point space has an isolated element.
\item The line space has an isolated element.
\item The flag space has an isolated element.
\end{enumerate}

\proof Conditions (i) and (ii) are equivalent by \ref{P2_2_2}.
If the point rows and pencils of lines are discrete, then the
big cells are discrete, and hence $\P$, $\L$, and $\F$ are discrete.
The other implications are clear.\qed}

\Def{We call a topological polygon {\em discrete,} if it satisfies one of 
the eight equivalent conditions of \ref{P2_2_6}.}

\Cor{\label{SchubertClosed}
If the polygon $\frak P$ is not discrete, then the Schubert varieties
are the topological closures of the corresponding Schubert cells.

\proof If $\frak P$ is not discrete, then the point rows are not 
discrete by \ref{P2_2_6}. Because point rows are homogeneous, no point on 
a point row is isolated. Now let $\P_k(u,v)$ be a Schubert cell, and 
let $q\in\P_{k-1}(v,u)$. Choose $\ell\in\L_q-\{f_{k-1}(u,q)\}$. The
corresponding punctured point row $L-\{q\}$ is contained in $\P_k(u,v)$,
hence $q$ is contained in the closure of $\P_k(u,v)$. Therefore, the
closure of $\P_k(u,v)$ contains $\P_{k-1}(v,u)$, and the claim follows
by induction on $k$.

The proofs for the other kinds of Schubert cells are similar.\qed}

%
%
%

\section{Automorphisms}

\Prop{Let $\phi$ be an automorphism of a topological $n$-gon.
If there is a point row $L$ and a pencil of lines $\L_p$,
such that the restrictions $L\TO{\phi|L} \phi(L)$ and 
$\L_p\TO{\phi|\L_p}\phi(\L_p)$ are continuous,
then $\phi$ is continuous on $\P$, $\L$, and $\F$. If $n$ is odd,
then it suffices that the restriction $L\TO{\phi|L}\phi(L)$ is continuous 
for some point row $L$.

In particular, every root collineation is a homeomorphism.

\proof The restriction of $\phi$ to any point row or pencil of
lines is continuous. For if $h$ is a line with point row $H$, then
there exists a vertex $z\in\opp\ell\cap\opp h$. Now 
$\phi|H=[\phi(h),\phi(z),\phi(\ell)]\circ(\phi|L)\circ[\ell,z,h]$ is
continuous. 

By \ref{GenPoly1_7} the restriction of $\phi$ to any Schubert cell
is continuous, hence $\phi$ is continuous.\qed}

This proposition is generalized to homomorphisms of topological
polygons in \cite{BK2}.

\section{Locally compact polygons}

\Def{A topological polygon $\frak P$ is called {\em locally compact},
if the point rows and the pencils of lines are locally compact.}

\Prop{Let $\frak P$ be a topological polygon. The following are 
equivalent:
\begin{enumerate}
\item The polygon $\frak P$ is locally compact.
\item The point space $\P$ is locally compact.
\item The line space $\L$ is locally compact.
\item The flag space $\F$ is locally compact.
\end{enumerate}

\proof If $\frak P$ is locally compact, then every Schubert cell is locally
compact, hence the spaces $\P$, $\L$, and $\F$ are locally compact.
Conversely, if $\P$ is locally compact, then the big cell in $\P$ is
locally compact, hence every punctured point row and every punctured
pencil of lines is locally compact.\qed}

The following theorem is proved in \cite{GKK}. There, it is stated for 
compact polygons, but the compactness is never really used in
the proof.

\Thm{\label{T2_4_3} \cite[1.5,1.6]{GKK}
Let $\frak P$ be a locally compact, non-discrete polygon.
Then the point space, the line space, the flag space, every point row,
and every pencil of lines is second countable. In particular, each
of these spaces is metrizable (and thus paracompact) and separable.

If $\frak P$ is in addition connected, then it is path connected and
hence locally contractible.\qed}

\Cor{Let $\frak P$ be a locally compact, connected polygon.
Fix a point $p\in\P$ and consider the bundle map
$\V_p^{(k)}\to\V_p^{(j)}$, $x_k\mapsto x_j$ defined in
\ref{SchubertBundle} for $0<j<k<n$. This bundle map admits
a section and is a homotopy
equivalence. For example in a generalized quadrangle, we
get homotopy equivalences
$\P-L \simeq \ell^\perp-\{\ell\} \simeq L$.

\proof The bundle has contractible fibers, thus it admits
a section and is shrinkable, see \cite{Dold2}.\qed}

\section{Compact polygons}

\Def{A topological polygon $\frak P$ is called {\em compact}, if
each point row and each pencil of lines is compact.}

\Prop{Let $\frak P$ be a topological polygon.
The following are equivalent:
\begin{enumerate}
\item The polygon $\frak P$ is compact.
\item The point space and the line space are compact.
\item The flag space is compact.
\end{enumerate}

\proof Suppose $\frak P$ is compact. If the base and the fiber of 
a locally trivial bundle are compact Hausdorff spaces, then the 
total space is also a compact Hausdorff space.

Let $(p,\ell)$ be a flag. The iterated bundle $\Gal_n(p,\ell)$ is compact. 
Thus the flag space $\F$ is compact, because it is the image of 
$\Gal_n(p,\ell)$ under the map 
$(x_0,\ldots,x_{n+1}) \mapsto \flag(x_n,x_{n+1})$.

The other implications are trivial.\qed}

I do not know if the compactness of the point space implies the
compactness of the line space in case that $n$ is even. It does, 
provided that $\frak P$ is connected, see \ref{P2_5_4}.

\Cor{(cp. \cite[1.2]{GKK})
Let $\frak P$ be a compact polygon. If the point rows are
finite, then the pencils of lines are finite.\qed}

\Prop{\label{P2_5_3}
(cp. \cite[2.1(a)]{GvM})
Let $\frak P$ be a generalized polygon. Suppose that
$\P$ and $\L$ are compact Hausdorff spaces, and that the
subspace $\F\SUB\P\times\L$ is closed. Then $\frak P$ is a
compact polygon.

\proof Clearly, the space of all $(n-1)$-chains 
$C^{(n-1)}$ and the space of all stammering $(n-1)$-chains $S^{(n-1)}$ are 
compact. Let $A=\{(x,y)\in\V\times\V|\ d(x,y)<n-1\}$. Then $A$ is compact, 
because it is the image of $S^{(n-1)}$ under the map
$(x_0,\ldots,x_{n-1})\mapsto (x_0,x_{n-1})$.

We want to show that the graph $G$ of the function $f_{n-1}$ is closed in
$X=((\V\times\V)-A)\times\V$. It is given by
\[
G=\{(x_0,x_{n-1},x_{n-2})|\ (x_0,\ldots,x_{n-1})\in C^{(n-1)}-S^{(n-1)}\}.
\]
The set $Y=\{(x_0,x_{n-1},x_{n-2})\in\V\times\V\times\V|\ (x_0,\ldots
x_{n-1})\in C^{(n-1)}\}$ is compact and hence closed in $\V\times\V\times\V$, 
and thus $G=X\cap Y$ is closed in $X$. Therefore, the map $f_{n-1}$
is continuous.\qed}

\Prop{\label{P2_5_4}
Let $\frak P$ be a connected, locally compact polygon.
Then $\frak P$ is compact.

\proof By duality, it suffices to show that some (and hence every)
pencil of lines is compact. Pick a point $q\in\P_{n-1}(u,v)$, and let 
$U\SUB\P_{n-1}(u,v)$ be an open neighborhood of $q$ with compact 
closure. Since $\frak P$ is connected, the boundary $\partial U$ 
of $U$ is compact and nonempty. Every point row $L$ through $p$ 
meets $\P_{n-2}(v,u)$, thus $L\cap \partial U$ is nonempty. 
Now the map $\partial U\to\L_q$, $x\mapsto f_2(x,q)$ is a continuous 
surjection.\qed}

For compact, non-connected polygons, we get the following result:

\Prop{\label{P2_5_6}
Let $\frak P$ be a compact polygon. If $\frak P$ is not connected,
then the point space, the line space, the flag space, every point row,
and every pencil of lines is homeomorphic to the Cantor set.

\proof Being totally disconnected, these spaces are zero-dimensional,
and the result follows from \cite[2-98]{HY}.\qed}

The following problem seems to be open:
If the point rows are (locally) compact, does it follow that the
pencils of lines are (locally) compact? Note that it is an open question
if there exist generalized polygons with finite point rows and
infinite pencils of lines.

\chapter{Finite-dimensional polygons}

The main theme of this chapter is the algebraic topology of a
finite-dimensional $n$-gon $\frak P=\PLF$, that is, a locally
compact $n$-gon
with the property that the dimensions $m,m'$ of the point rows and
of the pencils of lines are finite and positive (\ref{FiniteDDef}). 

Let $C$ be an $r$-dimensional
Schubert cell (in $\P$, say), let $X=\bar C$ be the corresponding Schubert 
variety, and put $A=X-C$. Thus $A$ is a smaller Schubert variety. In the
first section it is shown that $A\SUB X$ is a cofibration,
that $X,C$, and
$A$ are ANRs, and that $C$ is a generalized manifold. It follows that
$\P,\L$, and $\F$ are generalized manifolds.

The question is then how the homology of $X$ is related to the 
homology of $A$. It turns out that the pair $(X,A)$ has the
same homology as an $r$-sphere; in fact, the one-point 
compactification $X/A$ of the Schubert cell $C$ is a 
homotopy $r$-sphere (\ref{L3_4_2}). Hence
$X$ and $A$ have the same homology, except possibly in dimensions
$r,r-1$. We show that for a suitable choice of the coefficient domain $R$, the exact sequence
\[
0\to\H_r(A;R)\to\H_r(X;R)\to\H_r(X,A;R) \to \H_{r-1}(A;R) \to \H_{r-1}(X;R) \to 0
\]
breaks up to
\[
\begin{array}{rcl}
0 \to \H_r(A;R) \to \H_r(A;R)\oplus R \TO{\scriptstyle\ot} R & \!\to\! & 0 \\
         0 & \!\to\! & \H_{r-1}(A;R) \TO\cong \H_{r-1}(X;R) \to 0.
\end{array}
\]
Thus $\H_\bullet(X;R)=R^k$, where $k$ is the number of the Schubert 
cells contained in $X$. In particular, $\H_\bullet(\P;R)=R^n$.
The main idea is to use the gallery space to construct a linking cycle, 
a method developed by Palais, Terng, and Thorbergsson for 
isoparametric submanifolds. Similar results hold 
for the cohomolgy of $(X,A)$ (\ref{LinkCyc}).

Now suppose that $X$ is a Schubert variety in the flag space.
In the next section, we consider the double mapping cylinder
$DX$ of 
\[
\pr_1(X)\ot X \to \pr_2(X).
\]
It turns out that $DX$ is
contractible if $X\neq\F$, and that $D\F$ is a homotopy 
$(\dim\F+1)$-sphere. It follows from a theorem of M\"unzner
that $n\in\{3,4,6\}$, i.e. a finite-dimensional polygon is
a projective plane, a quadrangle, or a hexagon (\ref{346-Thm}).
This is one of our main results; it was first proved by Knarr
\cite{Kna2} under the additional assumption that the Schubert
cells are locally euclidean.
For polygons with equal parameters $m=m'$ we may calculate the
Steenrod squares; as an application, we show that the point
space does not factor as a product (\ref{NoProduct}). 
It is also shown that for a quadrangle with parameters $(m,1)$ 
the flag space $\F$ is a product $\F=\P\times\SS^1$, provided
that $\P$ is orientable (\ref{F2PTrivial}).

The connectivity properties of $DX$ are useful
in order to determine the homotopy properties of the
inclusion $A\SUB X$. Essentially, it turns out that the pair $(X,A)$ 
is $(\dim X-1)$-connected (\ref{RConn}), 
and thus the Schubert cell decomposition has
similar properties as a CW decomposition (there is a difficulty 
with the point and the flag space of
$(1,m')$-quadrangles for $m'>2$, though, due to the fact that
the fundamental groups of these spaces are non-trivial). 
The section closes with a table of the fundamental groups
of $\P,\L,$ and $\F$ (\ref{P3_2_6}).

In the last section we prove that the bundle map $\F\to\P$
does not admit a section, provided that $n\neq 4$ or that
$m=m'$ (\ref{NoSection}, with additional assumptions for the case of
$(1,1)$-quadrangles). This fact is useful in connection 
with automorphism groups, see Chapter \ref{PHomPoly}, and it leads
also to a non-existence result for closed ovoids in certain
polygons. It is proved in \ref{F2PTrivial} that for $n=4$,
the bundle $\F\to\P$ is sometimes trivial,
hence quadrangles may behave quite differently.

This seems to be a general phenomenon: There are strong analogies
between projective planes, quadrangles with equal parameters $m=m'$,
and hexagons, and one might conjecture
that for these geometries the dimension determines the homotopy
type (or even the homeomorphism type) of $\P,\L,$ and $\F$,
although a proof of this conjecture even for $n=3$ seems to run 
into some subtle and difficult
questions about characterizations of manifolds and on the
homotopy fiber ${\bf TOP}(k)/{\bf O}(k)$ for $k=4,8$ (but see \cite{Kr2}
for the smooth case).

In any case, this conjecture fails for quadrangles with $m\neq m'$
\cite{Wng,GTH91}.
So far,  little is known about these quadrangles, their
topology, and their automorphism groups.

\section{Dimension of Schubert cells}

We will use the following facts from dimension theory.

\Prop{\label{P3_0_7}
Let $A$ be a subspace of a second countable, metrizable space
$X$. Then the
{\em covering dimension $\dim A$}, the 
{\em large inductive dimension ${\rm Ind}\,A$}, and the
{\em small inductive dimension ${\rm ind}\,A$}
of $A$ coincide \cite[7.3.3]{Eng}, \cite[4 5.9]{Pears}.
Moreover, we have the inequality 
$$
\dim A\leq \dim X
\eqno(1)
$$
cp. \cite[7.3.4]{Eng}, \cite[4 6.4]{Pears}.
If $A$ is closed in $X$, then 
$$
\dim X = \max\{\dim A,\dim(X-A)\},
\eqno(2)
$$
cp. \cite[3 5.8, 4 4.2, 4 4.8]{Pears}, combined with (1). 
If $Z$ is a second countable, metrizable space, then the
inequality 
$$
\dim(X\times Z)\leq\dim X+\dim Z
\eqno(3)
$$
holds \cite[7.3.17]{Eng}, \cite[9 3.3]{Pears}.

Let $R$ be a principal ideal domain.
We denote the {\em sheaf-theoretic dimension} by $\dim_R$, cp. \cite{Bredon65}.
If $X$ is paracompact and if the covering dimension of $X$ is finite, then
$$
\dim_\ZZ X=\dim X,
\eqno(4)
$$
see \cite[4.1]{Loew}.
If $X$ is a locally compact Hausdorff space, then it follows readily
from the universal coefficient theorem that $\dim_RX\leq\dim_\ZZ X$
\cite[II 15.5, 18.3]{Bredon65}.
\qed}

\Def{\label{FiniteDDef}
A locally compact $n$-gon is called {\em finite-dimensional},
if the covering dimension $m=\dim L$ of the point rows and the covering
dimension $m'=\dim\L_p$ of the pencils of lines are finite and positive
(the case of zero-dimensional compact polygons is covered by \ref{P2_5_6}).
We call $(m,m')$ the {\em topological parameters} of the $n$-gon
$\frak P$, and we put $m_k(p,\ell)=m+m'+m+\ldots$ ($k$ summands), and
$m_k(\ell,p)=m'+m+m'+\ldots$ ($k$ summands). 

Note that a finite-dimensional polygon $\frak P$ is second countable, 
metrizable, path-connected, locally contractible, and compact by 
\ref{T2_4_3}, \ref{P2_5_4}.}

\Prop{\label{P3_0_9}
Let $\frak P=\PLF$ be a locally compact polygon.
The following are equivalent:
\begin{enumerate}
\item The polygon $\frak P$ is finite-dimensional.
\item The point space $\P$ has finite and positive covering dimension.
\item The line  space $\L$ has finite and positive covering dimension.
\item The flag  space $\F$ has finite and positive covering dimension.
\end{enumerate}

\proof If $\frak P$ is finite-dimensional, then every Schubert cell
is finite-dimensional (because of the product formula \ref{P3_0_7} (3)), 
and thus every Schubert variety is finite-dimensional by \ref{P3_0_7} (2)
(because it is the union of a Schubert cell and a smaller Schubert 
variety). Therefore $\P$, $\L$, and $\F$ are finite-dimensional.

Conversely, if $\P$ is finite-dimensional, then the point rows
and the pencils of lines are finite-dimensional by \ref{P3_0_7} (1), 
because they may be imbedded into $p^\perp\SUB\P$ by 
\ref{L2_1_6}.\qed}

{\em From now on, we assume that $\frak P$ is a finite-dimensional
$n$-gon with parameters $(m,m')$.}

\Prop{Every (punctured) point row and every (punctured) pencil of 
lines of $\frak P$ is an ANR (i.e. an absolute neighborhood retract for the class of all
metrizable spaces, see \ref{D6_1_7}).
Thus every Schubert cell is an ANR. The point space, 
the line space, and the flag space are ANRs.

\proof By \ref{T2_4_3}, every (punctured) point row is locally contractible. 
Being a finite-dimensional, locally contractible space,
it is an ANR \cite[V 7.1]{Hu}. A finite product of ANRs is again an ANR
\cite[III 7.6]{Hu}. A space which is locally an ANR is an ANR
\cite[III 8.1]{Hu}.\qed}

\Thm{\label{LinesMnf}
Every point row of $\frak P$ is a generalized $m$-manifold (see \ref{GenMnf}),
homotopy equivalent to an $m$-sphere. Similarly, every pencil of lines is a
generalized $m'$-manifold, homotopy equivalent to an $m'$-sphere.

If $m,m'\leq 2$, then the point rows and the pencils of lines are
homeomorphic to spheres, and every Schubert cell is homeomorphic to
some euclidean space. In particular, $\frak P$ is a manifold in the
sense of Chapter \ref{ManPoly}.
 
Every Schubert cell is a generalized manifold. Hence
\[
\begin{array}{ccccc}
\dim \P_k(u,v) & = & \dim \Cl\P_k(u,v) & = & m_k(u,v), \\
\dim \L_k(v,u) & = & \dim \Cl\L_k(v,u) & = & m_k(v,u), \\
\dim \F_k(u,v) & = & \dim \Cl\F_k(u,v) & = & m_k(u,v).
\end{array}
\]

The point space, the line space, and the flag space are 
generalized manifolds of dimension 
$\dim\P=n{m+m' \over 2}-m'$, $\dim\L=n{m+m' \over 2}-m$,
and $\dim\F=n{m+m' \over 2}$, respectively.

\proof The point rows and the pencils of lines are generalized
manifolds and homotopy spheres by L\"owen's Theorem \ref{Lowen}.
Being a generalized manifold is a local
property, hence every Schubert cell is a generalized manifold
by \ref{LocGenMnf}.\qed}

\Thm{\label{SchubertANR}
Let $\frak P$ be a finite-dimensional polygon. Then every Schubert
variety is an ANR (and thus locally contractible). 

For every $k\geq 0$ the gallery space $\Gal_k(u,v)$, as well as 
the subset of stammering galleries $\StamGal_k(u,v)$, is an ANR.
Moreover, $\Gal_k(u,v)$ is a generalized $m_k(u,v)$-manifold.
If $m,m'>1$, then $\Gal_k(u,v)$ is simply connected and hence
orientable.

Each of the following spaces is a compact ANR, and hence
the inclusions 
\[
\begin{array}{rcl}
\StamGal_k(u,v) & \SUB & \Gal_k(u,v) \\
\Cl\P_k(u,v)    & \SUB & \Cl\P_{k+1}(v,u) \\
\Cl\L_k(v,u)    & \SUB & \Cl\L_{k+1}(u,v) \\
\Cl\F_k(u,v)\cup\Cl\F_k(v,u) & \SUB & \Cl\F_{k+1}(u,v) \\
\Cl\F_{k+1}(u,v) & \SUB & \Cl\F_{k+1}(u,v)\cup\Cl\F_{k+1}(v,u)
\end{array}
\]
are cofibrations by \ref{D6_1_7}.

\proof Being an ANR is a local property \cite[III 8.1]{Hu}, hence
the iterated bundles $\Gal_k(u,v)$ are ANRs
by induction on $k$. By \ref{LocGenMnf},
the space $\Gal_k(u,v)$ is a generalized $m_k(u,v)$-manifold.
If $m,m'>1$, then it follows inductively from the homotopy exact
sequence for  bundles that the total
space of the bundle $\Gal_k(u,v)\to\Gal_{k-1}(u,v)$ is
simply connected.

Next, we want to show that the set $\StamGal_k(u,v)$ is an ANR. This is
trivial for $k=0$, so we proceed again by induction. Consider the maps
\[
\Gal_{k-1}(u,v)\TO{s}\Gal_k(u,v)\TO\pr\Gal_{k-1}(u,v).
\]
The set $A=\pr^{-1}(\StamGal_{k-1}(u,v))$ is an ANR, because it is the total
space of the bundle $\pr^{-1}(\StamGal_{k-1}(u,v))\to\StamGal_{k-1}(u,v)$,
and because $\StamGal_{k-1}(u,v)$ is an ANR. Clearly, the set
$B=s(\Gal_{k-1}(u,v))$, as well as the intersection 
$A\cap B=s(\StamGal_{k-1}(u,v))$ is an ANR. Thus the union 
$A\cup B=\StamGal_k(u,v)$ is an ANR by \ref{D6_1_7}.

Finally, we want to show that the Schubert varieties are ANRs. This
follows again by induction. Consider the relative homeomorphism
\[
(\Gal_k(u,v),\StamGal_k(u,v))\to(\Cl\P_k(u,v),\Cl\P_{k-1}(v,u)).
\]
By \ref{D6_1_7}, the space $\Cl\P_k(u,v)$ is an ANR. The proof for
the other types of Schubert varieties is similar. The fact that
$\Cl\F_k(u,v)\cup\Cl\F_k(v,u)$ is an ANR
follows again from \ref{D6_1_7}.\qed}

\clearpage

\section[Homological properties]%
{Homological properties of finite-dimensional pol\-y\-gons; first results}

\Def{\label{D3_1_1}
Throughout this section, we let $(X,A)$ and $(Y,B)$ 
denote one of the following pairs.
\[
\begin{array}{ll}
(X,A)   & (Y,B) \\   \hline
(\Cl\P_{k+1}(v,u),                     \Cl\P_k(u,v))                 &
(\Gal_{k+1}(v,u),                      \StamGal_{k+1}(v,u))          \\
(\Cl\L_{k+1}(u,v),                     \Cl\L_k(v,u))                 &
(\Gal_{k+1}(u,v),                      \StamGal_{k+1}(u,v))          \\
(\Cl\F_{k+1}(u,v),                     \Cl\F_k(u,v)\cup\Cl\F_k(v,u)) &
(\Gal_{k+1}(u,v),                      \StamGal_{k+1}(u,v)           \\
(\Cl\F_{k+1}(u,v)\cup\Cl\F_{k+1}(v,u), \Cl\F_{k+1}(u,v))             & 
(\Gal_{k+1}(v,u),                      \StamGal_{k+1}(v,u))
\end{array}
\]
Consider the canonical map $(Y,B)\to (X,A)$ that sends a gallery to
its end. The restriction
$Y-B \to X-A$ is a homeomorphism, and we get a homeomorphism
$(Y/B,*)\cong (X/A,*)$. We put $r=\dim (X-A)$.}

\Lemma{\label{LLinkCyc}
The pairs $(X,A)$ and $(Y,B)$ have the same (co-) homology
as $(\SS^r,*)$. The quotient space $Y/B$ is a generalized $r$-manifold and
a (co-) homology $r$-sphere. The collapsing map $Y\to Y/B$ has $R$-degree
$1$, i.e. it induces isomorphisms in $r$-dimensional (co-) homology
with coefficients in $R$, where $R=\ZZ$ for $m,m'>1$, and $R=\ZZ_2$ else.

\proof The space $Y/B=Y-B\cup\{*\}$ has finite covering
dimension by \ref{P3_0_7} (2),
thus its sheaf-theoretic dimension is finite as well. Now
$\H_\bullet(Y/B,Y-B)\cong{\bf\tilde H}_\bullet(Y/B)$, because $Y-B$ is
contractible. 
The quotient $Y/B$ is a wedge of point rows and pencils of lines, thus it 
has the same (co-) homology as an $r$-sphere by \ref{Wedge}.

Thus the space $Y/B$ has the 'right' local homology groups, and by
\ref{ANRMnf}, it is a generalized $r$-manifold. 

Now pick $y\in Y-B$, and consider the excision maps
\[
(Y,Y-\{y\})\ot(Y-B,Y-(B\cup\{y\}))\to(Y/B,Y/B-\{y\}).
\]
It follows that the collapsing map induces an ismomorphism
\[
\H_\bullet(Y,Y-\{y\};R)\to\H_\bullet(Y/B,Y/B-\{y\};R),
\]
and hence
\[
\H_r(Y;R)\TO\cong\H_r(X;R)
\]
is an isomorphism, since $X$ and $Y$ are $R$-orientable.\qed}

\Thm{\label{LinkCyc}

The maps $(A,\emptyset) \to (X,\emptyset) \to (X,A)$
induce short exact sequences
\[
\begin{array}{ccccccccc}
0 & \to & \H_\bullet(A;R) 
  & \to & \H_\bullet(X;R) 
  & \to & \H_\bullet(X,A;R) 
  & \to & 0 \\
0 & \ot & \H^\bullet(A;R) 
  & \ot & \H^\bullet(X;R) 
  & \ot & \H^\bullet(X,A;R) 
  & \ot & 0 
\end{array}
\]
which are split (note that
$\H_\bullet(X,A;R)\cong R\cong \H^\bullet(X,A;R)$ by \ref{LLinkCyc}). 
Again, the coefficient ring $R$ is $\ZZ_2$ if $1\in\{m,m'\}$, 
and $\ZZ$ else. Thus the attaching of the Schubert cell $X-A$ to the
Schubert variety $A$
corresponds to adding an $r$-dimensional (co-) homology class to
the (co-) homology of $A$.

In particular, every Schubert variety $X$ represents a (co-) homology class in
dimension $r=\dim (X-A)$, and the (co-) homology of a Schubert variety (and 
in particular of $\P$, $\L$ and $\F$) is (additively) a free $R$-module 
over the Schubert varieties contained in it.

The point space is orientable if $m>1$ (see also \ref{Stiefel-Whitney}).

\proof Consider the long exact sequence of the pair $(X,A)$. Since
$X/A$ has the same (co-) homology as an $r$-sphere, we have only
to show that $\H_r(X;R)\to \H_r(X,A;R)$ is surjective
(that $\H^r(X;R)\ot \H^r(X,A;R)$ is injective).

But this as well as the splitting follows from the diagram
\[
\begin{array}{ccccc}
Y    & \to & (Y,B) & \to & (Y/B,*) \\
\dto &     & \dto  &     & \DTO\cong \\
X    & \to & (X,A) & \to & (X/A,*) 
\end{array}
\]
and the fact that the composite $Y\to (Y,B) \to (X,A) $ induces an
isomorphism in $r$-dimensional (co-) homology by \ref{LLinkCyc},
\ref{Collapse}.\qed}

\Lemma{\label{L3_1_5}
If $m'>1$, then the inclusions $L\hookrightarrow p^\perp
\hookrightarrow\ldots\hookrightarrow\P$
induce isomorphisms in $1$-dimensional homology.
 
\proof By \ref{LinkCyc}, we have $\H_i(\P_{k+1}(v,u),\P_k(u,v))=0$ for
$k\geq 1$ and $i< m'+m$. Therefore we get isomorphisms
$\H_1(L)\TO\cong\H_1(p^\perp)\TO\cong\ldots\TO\cong\H_1(\P)$.\qed}

\section{The Veronese imbedding}

\Def{Let $A\SUB\F$ be a subset of the flag space. The double mapping 
cylinder $DA$ is defined as follows: on $A\times\II$ we introduce an 
equivalence relation by putting
\[
\begin{array}{ll}
(a,0)\sim (a',0) & \mbox{ if }\pr_1(a)=\pr_1(a') \\
(a,1)\sim (a',1) & \mbox{ if }\pr_2(a)=\pr_2(a') \\
(a,s)\sim (a,s)  & \mbox{ for }0\leq s\leq 1.
\end{array}
\]
The quotient $DA=(A\times\II)/\sim$ is called the {\em double mapping 
cylinder of $A$} 
\cite[1.29]{tDKP}. 
Note that the double mapping cylinder is obtained by
pasting the (unreduced) mapping cylinders of $A\to\pr_1(A)$ and 
$A\to\pr_2(A)$ together along $A$. Thus, if $A$ as well as $\pr_1(A)$,
$\pr_2(A)$ are ANRs, then $DA$ is also an ANR 
\cite[VI 1.2]{Hu}.

We may consider $\P$, $\L$, and $\F$ as subspaces of $D\F$.
These imbeddings are called the {\em topological Veronese imbeddings
of $\P$, $\L$, and $\F$}.

Note that up to homotopy, we may replace the maps $\pr_1(A)\ot A\to\pr_2(A)$
by the inclusions $DA-\L\hookleftarrow DA-(\P\cup\L) \hookrightarrow
DA-\P$.}

Our next aim is the following result.

\Thm{\label{DMC}
The double mapping cylinder $D\F$ of the flag space is homotopy 
equivalent to an $(n(m+m')/2+1)$-sphere.}

We fix a flag $(p,\ell)\in\F$, and we put $X_k=\Cl\F_k(p,\ell)$,
$Y_k=\Cl\F_k(\ell,p)$. The one-point compactification of a space
$Z$ is denoted by $Z^{+}$.

As an intermediate step, we prove the following result.

\Prop{For $0\leq k <n$, the double mapping cylinders $DX_k$, $DY_k$ and
$D(X_k\cup Y_k)$ are contractible.

\proof
This is certainly true for $DX_0=DY_0=\{(p,\ell)\}\times\II$.
Now we assume that the double mapping cylinders $DX_{k-1}$, $DY_{k-1}$, and
$D(X_{k-1}\cup Y_{k-1})$ are contractible.
Then the collapsing map 
\[
DX_k\mapsto DX_k/D(X_{k-1} \cup Y_{k-1})
\]
is a homotopy equivalence by \ref{CollapseA}.
Now 
\[
DX_k-D(X_{k-1} \cup Y_{k-1})=
\F_k(p,\ell)\times (0,1]
\]
for $k<n$.
The one-point compactification of this space is just the reduced cone 
$C_{\F_k(p,\ell)^+}$ of the one-point compactification of $\F_k(p,\ell)$ 
and thus contractible. Hence we find that 
\[
DX_k/D(X_{k-1} \cup Y_{k-1})=
C_{\F_k(p,\ell)^+}
\]
is contractible. Similarly, 
$DY_k/D(X_{k-1}\cup Y_{k-1})=C_{\F_k(\ell,p)^+}$ and
$D(X_k\cup Y_k)/D(X_{k-1}\cup Y_{k-1})=C_{\F_k(p,\ell)^+}\vee
C_{\F_k(\ell,p)^+}$ are contractible.\qed}

{\em Proof of the theorem.} Since $D(X_{n-1}\cup Y_{n-1})$ is contractible,
$D\F$ is homotopy equivalent to the one-point compactification
of $\F_n(p,\ell)\times (0,1)$, and this is precisely the reduced
suspension of $\F_n(p,\ell)^{+}$.

Now $\pi_1S(\F_n(p,\ell)^+)=0$, because $\F_n(p,\ell)^{+}=\F/(X_{n-1}\cup
Y_{n-1})$ is path connected. Thus $D\F$ is a simply connected homology 
$(n(m+m')/2+1)$-sphere by \ref{LLinkCyc}.
Being a compact ANR, the space $D\F$ is homotopy
equivalent to a CW complex, and hence to an $(n(m+m')/2+1)$-sphere by
\ref{CWSphere}.\qed

\Thm{The double mapping cylinder $D\F$ is a generalized
$(n(m+m')/2+1)$-manifold.

\proof See \ref{ManBundle}.\qed}

\Cor{\label{MVS}
For any coefficient domain $R$, we have the relations
\[
\H_{\dim\F}(\F;R) \cong R \cong \H^{\dim\F}(\F;R);
\]
in particular, the flag space is orientable.
The maps $\pr_1,\pr_2$ induce isomorphisms
$\H_j(\F;R)\TO\cong\H_j(\P;R)\oplus\H_j(\L;R)$ and
$\H^j(\F;R)\OT\cong\H^j(\P;R)\oplus\H^j(\L;R)$
in dimension $j$, for $0<j<\dim\F$.

\proof This follows from the Mayer-Vietoris sequence of
$(D\F;D\F-\P,D\F-\L)$.\qed}

The following theorem was first proved by Knarr \cite{Kna2}
under the assumption
that the polygon $\frak P$ is a manifold (see Chapter \ref{ManPoly}).

\Thm{\label{346-Thm}
Let $\frak P$ be a finite-dimensional $n$-gon with parameters
$(m,m')$. Then $n\in\{3,4,6\}$.
If $n=3$, then $m=m'\in\{1,2,4,8\}$.
If $n=6$, then $m=m'\in\{1,2,4\}$. 
If $n=4$, and if $m,m'>1$, then
either $m=m'\in\{2,4\}$, or $m+m'$ is odd.

The cohomology rings of $\P$, $\L$, and $\F$ are listed in section
\ref{MueThm}. 
Some of the Steenrod squares are calculated in \ref{SteenrodAlg}.

\proof This follows from M\"unzner's Theorem, see \ref{Muenzner}.\qed}

In case that $n\neq 4$ or that $m=m'$, it is easy to calculate the
Steenrod squares. We use the notation of \ref{MueThm}
for the $\ZZ_2$-cohomology.

\Lemma{\label{SteenrodAlg}
Let $\frak P$ be a finite-dimensional $n$-gon. If $n\neq 4$,
or if $m=m'$, then the Steenrod squares are as follows:

If $n=3,6$, then $\Sq\,x_m=x_m+x_m^2$, and $\Sq\,y_m=y_m+y_m^2$.
If $n=6$, then $\Sq\,x_{3m}=x_{3m}+x_mx_{3m}$, and $\Sq\,y_{3m}=y_{3m}+y_m
y_{3m}$.
If $n=4$, then $\Sq\,x_m=x_m+x_m^2$, $\Sq\,y_m=y_m$, and $\Sq\,y_{2m}=y_{2m}+
y_my_{2m}$.

\proof The operation of $\Sq$ on $x_m,y_m$ follows readily
from the properties of $\Sq$, see eg. \cite[VI 15.]{Bredon93}.
Now let $n=4$. Then $y_{2m}=x_m^2+x_my_m$,
hence $\Sq\,y_{2m}=x_m^2+x_my_m+x_m^2y_m=y_{2m}+y_my_{2m}$.
For $n=6$ consider $x_m^2y_m=x_{3m}+y_{3m}$. Applying $\Sq$ to both sides
yields $x_m^2y_m+x_m^2y_m^2=\Sq\,x_{3m}+\Sq\,y_{3m}$. Thus $\Sq\,x_{3m}=
x_{3m}+x_mx_{3m}$ and $\Sq\,y_{3m}=y_{3m}+y_my_{3m}$.\qed}

\Thm{\label{NoProduct}
Let $\frak P$ be a finite-dimensional $n$-gon.
If $n\neq 4$, or if $m=m'$, then neither $\P$ nor $\L$ is 
homeomorphic to a product of two topological spaces
(containing more then one point, of course).

\proof Recall that the $\ZZ_2$-Poincar\'e polynomial of a space $X$ is defined
by 
$P_X(t)=\sum_{i\geq 0}{\beta_i}t^i$, where $\beta_i=\dim_{\ZZ_2}\H^i(X;\ZZ_2)$.

A generalized $k$-manifold $M$ factors as a product $M=X\times Y$
only if $X$ and $Y$ are generalized $i$- and $(k-i)$-manifolds
\cite[V 15.8]{Bredon65}. We have the relation $p_X(t)p_Y(t)=p_M(t)$
for the $\ZZ_2$-Poincar\'e polynomials (taken with respect to
sheaf-theoretic $\ZZ_2$-cohomology). There is
a natural isomorphism $\H^\bullet(\P;\ZZ_2)\cong\Hbar^\bullet(\P;\ZZ_2)$,
because $\P$ is compact and $HLC$ (\ref{CompareHom}).

Next, note that the $\ZZ_2$-cohomology of a space with $\ZZ_2$-Poincar\'e
polynomial $1+t^k$ is necessarily of the form $\ZZ_2[x_k]/(x_k^2)$, and
the $\ZZ_2$-cohomology of a generalized manifold with with $1+t^k+t^{2k}$
as $\ZZ_2$-Poincar\'e polynomial is of the form $\ZZ_2[x_k]/(x_k^3)$
by \cite[V 10.6]{Bredon65} (since every generalized manifold is 
$\ZZ_2$-orientable).

\smallskip
If $n=3$, then the $\ZZ_2$-Poincar\'e polynomial of $\P$
is given by $1+t^m+t^{2m}$, with $m\in\{1,2,4,8\}$. These polynomials do
not factor over $\NN$.

If $n=4$, then the $\ZZ_2$-Poincar\'e 
polynomial is $1+t^m+t^{2m}+t^{3m}=
(1+t^m)(1+t^{2m})$, and $m\in\{1,2,4\}$. Now 
\[
\H^\bullet(\P;\ZZ_2)=
\ZZ_2[x_m]/(x_m^4),
\]
and
\[
\H^\bullet(\L;\ZZ_2)=\ZZ_2[y_m,y_{2m}]/(y_m^2,y_{2m}^2).
\]

Clearly, the $\ZZ_2$-cohomology of $\P$ is not isomorphic to 
\[
[\ZZ_2[u_m]/(u_m^2)]^\bullet\otimes_{\ZZ_2}
[\ZZ_2[v_{2m}]/(v_{2m}^2)]^\bullet
=\ZZ_2[u_m,v_{2m}]/(u_m^2,v_{2m}^2),
\]
hence $\P$ does not factor.
The $\ZZ_2$-cohomology ring of $\L$ is of this form; however,
we have the relation $\Sq^m y_{2m}=y_{2m}y_m\not\in\bra{1,y_{2m}}$.
Since the Steenrod squares commute with maps, we conclude that
$\L$ does not factor.

If $n=6$, then the $\ZZ_2$-Poincar\'e polynomial of $\P$ and $\L$ is
$1+t^m+t^{2m}+t^{3m}+t^{4m}+t^{5m}$, for $m\in\{1,2,4\}$. It
factors as $(1+t^m)(1+t^{2m}+t^{4m})=(1+t^m+t^{2m})(1+t^{3m})$
over $\NN$. The first factorization is excluded by the structure
of the $\ZZ_2$-cohomology ring, and the second possibility is excluded
again by the fact that $\Sq^mx_{3m}=x_{3m}x_m\not\in\bra{1,x_{3m}}$.\qed}

Note that for quadrangles with $m\neq m'$, the point space can be a
product of two spheres \cite{FKM,Wng,GTH91}.

\Thm{\label{F2PTrivial}
Let $\frak P$ be a finite-dimensional quadrangle with $m'=1$.
If $m=1$, assume in addition that $\P$ is orientable, see
\ref{Stiefel-Whitney}. Then the circle bundle $\F\to\P$ is trivial, 
i.e. $\F=\SS^1\times\P$.

\proof By Kneser's result \ref{Kneser}, there is an
orthogonal $2$-plane bundle $\xi$ over $\P$ that has $\F\to\P$ 
as its circle bundle, and by assumption $\xi$ is orientable. We may
identify the total space of $\xi$ with $D\F-\L$.
Consider the diagram
\[
\begin{array}{ccccccc}
(\RR^2,\RR^2-\{0\}) & \to & (D\F-\L,D\F-(\P\cup\L)) & \ot & D\F-\L & \to & \P \\
                    &     & \dto                    &     & \dto \\
                    &     & (D\F,D\F-\P)            & \ot & D\F.
\end{array}
\]
Since $\H^2(D\F)=0$, the Euler class of $\xi$, i.e. the image of the
orientation class $U_\xi\in\H^2(D\F-\L,D\F-(\P\cup\L))$ in $\H^2(\P)$,
vanishes. Thus the bundle  $\xi$ is trivial by \ref{LineBundle}.\qed}

\section{Homotopy properties}

We may use the double mapping cylinder to calculate some homotopy
groups of the Schubert cell decomposition.
We use the notation of the previous sections.

\Thm{\label{RConn}
Suppose that $m,m'>1$. Let $(X,A)$ denote one of the pairs
of \ref{D3_1_1}, and put $r=\dim (X-A)$. Then the pair $(X,A)$ is
$(r-1)$-connected, and $\pi_r(X,A)=\ZZ$. In particular, every Schubert
variety is $(\min\{m,m'\}-1)$-connected.

\proof By \ref{LinesMnf} the Schubert varieties $L=\Cl\P_1(p,\ell)$
and $\L_p=\Cl\L_1(\ell,p)$ are simply connected. 
Now suppose that $\Cl\P_k(u,v)$ is simply connected. It follows
from the homotopy sequence for bundles that the total space 
$\F_{\Cl\P_k(u,v)}=\Cl\F_{k+1}(u,v)$ is simply connected. Up to 
homotopy, we may replace the diagram 
\[
\Cl\P_k(u,v)\ot\Cl\F_{k+1}(u,v)\to\Cl\L_{k+1}(u,v)
\]
by the inclusions
\[
D\Cl\F_{k+1}(u,v)-\L \hookleftarrow D\Cl\F_{k+1}(u,v)-(\L\cup\P)
\hookrightarrow D\Cl\F_{k+1}(u,v)-\P.
\]
Since the double mapping cylinder $D\Cl\F_{k+1}(u,v)$ is
simply connected by \ref{DMC}, we conclude from the Seifert-Van Kampen
Theorem that $\Cl\L_{k+1}(u,v)$ is simply connected.
Thus every Schubert variety is simply connected.
By the Seifert-Van Kampen Theorem, the union $\Cl\F_k(u,v)\cup
\Cl\F_k(v,u)$ is also simply connected, see
\cite[II 2.5 and p.94]{Whitehead}. From the
relative Hurewicz isomorphism Theorem \cite[7 5.4]{Spanier}
we infer that
$\pi_s(X,A)\cong\H_s(X,A)=0$ for $s<r$, and that $\pi_r(X,A)\cong
\H_r(X,A)=\ZZ$.\qed}

\Cor{If $m,m'>1$, then $\P_k(u,v)$ is homotopy equivalent to a
CW complex $X$ with $k+1$ cells
\[
X=e^0\cup e^m\cup e^{m+m'}\cup \ldots \cup e^{m_k(u,v)}.
\]
A similar statement holds for the other kinds of Schubert varieties.

\proof The space $\P_k(u,v)$ is simply connected by \ref{RConn},
and $\H_\bullet(\P_k(u,v))\cong\ZZ^k$, with $\ZZ$-Poincar\'e polynomial
$1+t^m+t^{m+m'}+\cdots+t^{m_k(u,v)}$ by \ref{LinkCyc}.
Hence the claim follows from \cite[4.1]{Wall}.\qed}

\Lemma{\label{L3_4_2}
The quotient $X/A$ is homotopy equivalent to an $r$-sphere.

\proof By \ref{LinesMnf} and by induction, it suffices to consider
the following two cases:

\smallskip (1) $1\in\{m,m'\}$. Then the Schubert cell
$X-A$ is homeomorphic to
$\RR\times K$, where the one-point compactification $K^+$ is a homotopy 
sphere. The one-point compactification of $\RR\times K$ is the reduced
suspension of $K^+$, and hence simply connected. 

\smallskip (2) $m,m'>1$. Then the Schubert cell
$X-A$ is homeomorphic to $K\times L$,
where $K^+$ and $L^+$ are simply connected homotopy spheres. Thus
$\pi_1(K^+\times L^+)=0=\pi_1(K^+\vee L^+)$, cp. \cite[II 5.9]{tDck}.
From the relative Hurewicz
isomorphism, we find that $\pi_1(K^+\times L^+,K^+\vee L^+)\cong\H_1
(K^+\times L^+,K^+\vee L^+)=0$. It follows from \ref{HETC} that
$\pi_1(K\times L)^+\cong\pi_1(K^+\wedge L^+)=0$.

\smallskip Hence in each case the quotient $X/A$ is simply connected, and the
claim follows from \ref{LLinkCyc} and \ref{CWSphere}, since
every ANR is homotopy equivalent to a CW complex \cite[p.218]{Web}.\qed}

The following corollary may be used to calculate some homotopy
groups of the Schubert varieties.

\Cor{If $m,m'>1$, then $\pi_s(X,A)\cong\pi_s(\SS^r)$ for $0\leq s\leq
r+\min\{m,m'\}-2$, where $r=\dim(X-A)$.

\proof This follows from \ref{HETC}.\qed}

We will see in \ref{CWPoly}
that the Schubert cell decomposition is indeed a
CW decomposition, provided that $m,m'\leq 2$, or if the point rows
and pencils of lines are locally euclidean. Hence the pair
$(X,A)$ is also $(r-1)$-connected in these cases.

It remains to investigate quadrangles with parameters
$(1,m')$ and $m'>2$. In this case we do not get a complete
result. However, we may calculate the fundamental groups of the
Schubert varieties.

\Prop{Suppose that $\frak P$ is a finite-dimensional
quadrangle with parameters
$m=1$ and $m'>1$. Then $L$ is homeomorphic to
$\SS^1$, and $p^\perp$ is the one-point compactification of
$\L_p\times\RR$. The inclusions $L\SUB p^\perp\SUB\P$ induce
isomorphisms $\ZZ=\pi_1L\TO\cong\pi_1p^\perp\TO\cong\pi_1\P$.
The Schubert varieties $\L_p,\ell^\perp$ and $\L$ are
$(m'-1)$-connected.

\proof We use the same ideas as in the proof of \ref{RConn}.
The bundle map $\F_L\to L$ induces an isomorphism on the
fundamental groups, because the fiber is $(m'-1)$-connected.
Again, we conclude from the Seifert-Van Kampen Theorem, applied
to $D\F_L$, that $\ell^\perp$ is simply connected. Similarly,
it follows from the diagram $p^\perp\ot\F_{p^\perp}\to\L$
that $\L$ is simply connected. The line bundle
$p^\perp-\{p\}\to\L_p$ is trivial, because $\L_p$ is a
cohomology $m'$-sphere, see \ref{LineBundle} (using this fact,
it is not difficult to see that the pair $(p^\perp,L)$ is in fact
$m'$-connected). In \ref{P3_2_6} we prove
that $\pi_1\P\cong\ZZ$. Thus, the fundamental group of every
Schubert variety in $\P$ is abelian, and the claim about
the inclusions $L\SUB p^\perp\SUB\P$ follows from \ref{L3_1_5}.\qed}

\Prop{Let $\frak P$ be a finite-dimensional polygon. The fundamental group of $\F$ is generated by the images of $\pi_1(L\times\{\ell\})\cong\pi_1\SS^m$ and
$\pi_1(\{p\}\times\L_p)\cong\pi_1\SS^{m'}$ in $\pi_1\F$.

The fundamental group of the point space is generated by the image of
$\pi_1 L$ in $\pi_1\P$ (and dually, the fundamental group of the line
space is generated by the image of $\pi_1\L_p$).

\proof (cp. \cite[3.5]{GH})
It follows from the exact sequence of a fibration that the maps $\pi_1\F\to\pi_1\P$, $\pi_1\F\to\pi_1\L$ are
surjective, and the kernels of these maps are precisely the images
of $\pi_1(\{p\}\times\L_p)$ and $\pi_1(L\times\{\ell\})$, respectively.
Now we may replace $\P\ot\F\to\L$ by $D\F-\L\ot D\F-\{\L\cup\P\}\to
D\F-\P$, and the claim follows from the Seifert-Van Kampen Theorem.\qed}

\Cor{\cite[3.5]{GH} The fundamental group $\pi_1\F$ acts trivially on the
homotopy groups $\pi_s\F$ for $s>1$.

\proof Suppose $m'=1$.
It follows rom the homotopy exact sequence of the circle bundle
$\F\to\P$ that $\pi_s\F\to\pi_s\P$ is an injection for $s=2$, and
an isomorphism for $s>2$, and the projection $\F\to\P$ kills the image
of $\pi_1(\{p\}\times\L_p)$.
Thus $\pi_1(\{p\}\times\L_p)$ operates trivially on $\pi_s\F$.\qed}

\Cor{(cp. \cite[4.8]{GH})
If $\F\to\P$ is non-orientable, then the image of $\pi_{m'}(\{p\}\times
\L_p)$ in $\pi_{m'}\F$ is cyclic of order at most $2$.

\proof If $\F\to\P$ is non-orientable, then some $a\in\pi_1\F$ acts
(through the epimorphism $\pi_1\F\to\pi_1\P$) as $\bf-1$ on $\pi_{m'}(\{p\}\times\L_p)=\ZZ$. On the other hand, $a$
acts trivially on $\pi_{m'}\F$, hence every element in the image of $\pi_{m'}(\{p\}\times\L_p)\cong\ZZ$ is of order at most $2$.\qed}

\Cor{(cp. \cite[4.8]{GH})
Let $\frak P$ be a finite-dimensional quadrangle with parameters $(1,m' )$, for some $m'>1$. Then $\P$ and $\F\to\P$ are orientable if and only if
$m'$ is even.

\proof If $\F\to\P$ is orientable, then $1+m'$ is odd by
\cite[Satz 7]{Mue}.

Suppose that $\F\to\P$ is non-orientable. Then there is a
two-fold covering $\tilde\P\to\P$ such that the induced bundle
$\tilde\F\to\tilde\P$ is orientable. By the corollary above,
the image of $\pi_{m'}(\{p\}\times\L_p)$ in $\pi_{m'}\tilde\F$
is finite. From the diagram
\[
\begin{array}{ccc}
\pi_{m'}(\{p\}\times\L_p) & \to & \pi_{m'}\tilde\F \\
\DTO\cong                 &     & \dto \\
\H_{m'}(\{p\}\times\L_p)  & \to & \H_{m'}(\tilde\F) 
\end{array}
\]
we conclude that $\H_{m'}(\{p\}\times\L_p;\QQ)\to\H_{m'}(\F;\QQ)$
is trivial, and the claim follows from \ref{OddBundle}.\qed}

\Cor{\label{OddCor}
The integral cohomology of a quadrangle with parameters $(1,2k)$ is
as given in \ref{MueThm} $\bf4_3$.\qed}

\Thm{\label{P3_2_6}
Let $\frak P$ be a finite-dimensional $n$-gon with parameters
$(m,m')$. Up to duality, the fundamental groups of $\P$, $\L$, and
$\F$ are as follows:

\begin{center}
\begin{tabular}{l|ccc}
                  & $\pi_1\P$ & $\pi_1\L$      & $\pi_1\F$      \\ \hline 
$m,m'>1$          & $0$       & $0$            & $0$            \\
$m=1$, $m'>1$     & $\ZZ$     & $0$            & $\ZZ$          \\
$m=m'=1$, $n=3,6$ & $\ZZ_2$   & $\ZZ_2$        & $Q_8$          \\
$m=m'=1$, $n=4$   & $\ZZ_2$   & $\ZZ$          & $\ZZ+\ZZ_2$    
\end{tabular}
\end{center}

Here, $Q_8$ denotes the quaternion group of order $8$. In the case $n=4$,
$m=m'=1$, the point space $\P$ is orientable, and $\L$ is not,
see \ref{Stiefel-Whitney}.
(A simply connected generalized manifold is orientable.)

\proof Up to homotopy, we may replace the maps $\P\OT{\pr_1}\F\TO{\pr_2}\L$
by the inclusions $D\F-\L\hookleftarrow D\F-(\P\cup\L)
\hookrightarrow D\F-\P$. Assume that $m'>1$.
Since $\pr_1$ is a locally trivial fibration whose fiber is a homotopy 
$m'$-sphere, the induced map $\pi_1\F\to\pi_1\P$ is an isomorphism.
The Seifert-Van Kampen Theorem, applied to 
$(D\F-\L,D\F-\P,D\F-(\P\cup\L))$ yields $\pi_1\L=0$.
From the homotopy exact sequence of $\SS^m\to\F\to\L$, we deduce
that $\pi_1\F$ is cyclic, and hence isomorphic to $\H_1(\F)$.

The case $m=m'=1$ will be treated separately in \ref{11FundGrp}.\qed}

\section{Sections and ovoids}

From the knowledge of the cohomology of finite-dimensional polygons,
we get the following result which was first obtained by Breitsprecher 
\cite[2.3.2]{Brsp} for projective planes.

\Thm{\label{NoSection}
Let $\frak P$ be a finite-dimensional $n$-gon
with parameters $(m,m')$. If $n\neq 4$, or if $m=m'>1$, then
the bundles $\F\TO{\pr_1}\P$ and $\F\TO{\pr_2}\L$ admit no sections.

If $\frak P$ is a quadrangle with parameters $(1,1)$, and if the point
space is orientable (see \ref{Stiefel-Whitney}),
then the bundle $\F\to\L$ admits
no section (the bundle $\F\to\P$ is trivial by \ref{F2PTrivial}).

\proof Assume that there is a section $s:\P\to\F$. We put $R=\ZZ_2$ for
$n=3,6$, and $R=\QQ$ for $n=4$, $m=m'>1$.
Let $x_1,x_2,y_1,y_2$ be generators of $\H^m(\P;R)$, $\H^{2m}(\P;R)$,
$\H^m(\L;R)$, and $\H^{2m}(\L;R)$, respectively. By \ref{MVS} we may identify
$\H^m(\P;R)\oplus\H^m(\L;R)$ with $\H^m(\F;R)$, and 
$\H^{2m}(\P;R)\oplus\H^{2m}(\L;R)$ with $\H^{2m}(\F;R)$.
We have the relations $x_1^2=\alpha x_2$
$y_1^2=\beta y_2$, and $x_1y_1=x_2+y_2$, with
$\alpha,\beta\in R$. Note that $s^\bullet\pr_1^\bullet=\id$.

For $n=3,6$, we get $(\alpha,\beta)=(1,1)$, and for $n=4$ we get
$(\alpha,\beta)\in\{(1,2),(2,1)\}$, according to \ref{MueThm}. 
Now let
$ax_1+y_1$ be a generator of $\ker(s^\bullet:\H^m(\F;R)\to\H^m(\P;R))$.
Since $s^\bullet$ is a ring homomorphisms, we have
\begin{eqnarray*}
0 & = & s^\bullet((ax_1+y_1)(y_1-(\beta+a)x_1)) \\
  & = & s^\bullet(-(a^2\alpha+a\alpha\beta+\beta)x_2) \\
  & = & -(a^2\alpha+a\alpha\beta+\beta)x_2.
\end{eqnarray*}
But the polynomial $f(a)=a^2\alpha+a\alpha\beta+\beta$ has no roots in $R$,
hence we get a contradiction.

In the case of an $(1,1)$-quadrangle, let $ay_1+x_1$ be in the kernel
of 
\[
s^\bullet:\ZZ_2[x_1,y_1,y_2]/(x_1^4,y_1^2,y_2^2,x_1^2+y_2+x_1y_1)
\to\ZZ_2[y_1,y_2]/(y_1^2,y_2^2).
\]
Then $0=s^\bullet(((1+a)y_1+x_1)(ay_1+x_1))=s^\bullet(y_1x_1+x_1^2)=y_2$,
a contradiction.\qed}

\Cor{{\bf(Breitsprecher)}
Every continuous collineation $\phi$ of a finite-di\-men\-sion\-al projective
plane has a fixed point.

\proof Otherwise, the map $p\mapsto(p,p\vee\phi(p))$ would be a 
section.\qed}

Generalized quadrangles with $m\neq m'$ may have sections; in fact, 
the bundle $\F\to\P$ is sometimes trivial, see \ref{F2PTrivial}.
Every compact connected Moufang quadrangle
or hexagon admits an involution with no fixed point and no fixed line
(because every isoparametric foliation is invariant under the antipodal
map of the ambient sphere).
\smallskip

An {\em ovoid} in a generalized $n$-gon (for $n\geq 4$)
is a set of points $O\SUB\P$ with
the property that every point row meets $O$ in exactly one
point. Cameron has proved that the point space of 
every generalized $n$-gon with infinite point rows and infinite 
pencils of lines of the same
cardinality may be partitioned into ovoids \cite{Cam}.

\Cor{Let $\frak P$ be a finite-dimensional $n$-gon, for some 
$n\geq 4$. If $n=4$,
assume in addition that $m=m'$, and that $\P$ is orientable. Then
$\P$ contains no closed ovoid.

\proof Let $O\SUB\P$ be a closed and therefore compact
ovoid. Then $\F_O\to\L_O=\L$ is a 
continuous bijection and hence a homeomorphism. Thus the bundle
$\F\to\L$ has a section.\qed}

\chapter{Polygons which are manifolds}
\label{ManPoly}

In this section we investigate polygons which are locally
euclidean. In this case the Schubert cells are homeomorphic
to some euclidean spaces, and the point rows and the pencils
of lines are spheres (more generally the one-point
compactification of every Schubert cell is a sphere).
This fact is due to the contractiblity properties of the Schubert
cells.
The Schubert cell decomposition of $\P,\L,$ and $\F$ is a CW
decomposition (\ref{CWPoly}), a fact that was proved for projective planes
by Breitsprecher \cite{Brsp}. The proof becomes surprisingly simple, if
one constructs first a CW decomposition of the gallery spaces
(which is easy).

The double mapping cylinder is a manifold; it follows from the
proof of the generalized Poincar\'e conjecture that it is indeed
a sphere \cite{Kna2}. Therefore $\P,\L,$ and $\F$ may be
imbedded into $\SS^{\dim\F+1}$ (with normal disk bundles) (\ref{VeroneseRev}).
We calculate the Stiefel-Whitney classes of these spaces,
as well as the Stiefel-Whitney classes of the normal bundles,
provided that $n\neq 4$ or that $m=m'$ (\ref{Stiefel-Whitney}).

In the next section we obtain a complete topological classification of the
bundles $p^\perp-\{p\}\to\L_p$, provided that $m=m'\leq 2$ (\ref{T4_2_1}).
For $n=3$, this is of course Salzmann's and Breitsprecher's topological
classification of $\P$. For $n=3$ and $m=m'\leq 2$, Knarr's
Veronese imbedding leads to a simple proof of the topological
classification of the flag space (\ref{BrspBuc}).

The last section gives a topological criterion for a partial
$n$-gon $\frak P$ to be an $n$-gon: 'if the
dimension is right', then $\frak P$ is a compact $n$-gon.

\section{A CW decomposition for $\P,\L$, and $\F$}

\Def{A topological polygon $\frak P=\PLF$ {\em is a manifold}, if the 
point rows and the pencils of lines are locally homeomorphic to
$\RR^m$ and $\RR^{m'}$, respectively, for some numbers $m,m'>0$.
Thus the parameters of $\frak P$ are $(m,m')$.
Note that $\P$, $\L$, and $\F$ are second countable and
metrizable by \ref{T2_4_3}.

A finite-dimensional polygon with parameters $m,m'\leq 2$
is a manifold by \ref{LinesMnf}.}

\Prop{Suppose the topological $n$-gon $\frak P=\PLF$ is a manifold
with parameters $(m,m')$. Then every point row is homeomorphic
to an $m$-sphere, and every pencil of lines is homeomorphic to
an $m'$-sphere. The point space $\P$, the line space $\L$, and
the flag space $\F$ are compact connected manifolds of dimension
$n(m+m')/2-m'$, $n(m+m')/2-m$ and $n(m+m')/2$, respectively. The maps
$\F\to\P$ and $\F\to\L$ are locally trivial $m'$- and $m$-sphere bundles,
respectively.

\proof Every punctured point row and every punctured pencil of
lines is a pseudo-isotopically contractible manifold. By the result
of Harrold \cite{Harrold}, it is homeomorphic to some euclidean space.\qed}

\Thm{\label{CWPoly}
Suppose the topological $n$-gon $\frak P=\PLF$ is a manifold
with parameters $(m,m')$. Fix a flag $(p,\ell)$. The Schubert cell
decomposition of $\frak P$ with respect to $(p,\ell)$ is a CW decomposition. Thus, the point space
$\P$ consists of $n$ cells of dimension $0,m,m+m',m+m'+m,\ldots$. Similarly,
the line space $\L$ consists of $n$ cells, and the flag space consists
of $2n$ cells.}

In order to prove the theorem, we need the following lemma.

\Lemma{The space $\Gal_k(u,v)$ is a manifold.
It is obtained from the subspace $\StamGal_k(u,v)$ 
by attaching an $m_k(u,v)$-cell,
\[
\Gal_k(u,v)=\StamGal_k(u,v)\cup e^{m_k(u,v)}.
\]

\proof We proceed by induction on $k$. For $k=0$, there is nothing
to show. Now we consider the sphere bundle $\pr:\Gal_{k+1}(u,v)\to
\Gal_k(u,v)$, with the section $s$ defined in \ref{GenPoly1_2}. 
The proper galleries in $\Gal_k(u,v)$ are contained in an
$m_k(u,v)$-cell. Hence by \ref{CWBundle}, the proper galleries
$\PropGal_{k+1}(u,v)=\pr^{-1}(\PropGal_k(u,v))-s(\PropGal_k(u,v))$
are contained in an $m_{k+1}(u,v)$-cell 
\[
(e^{m_{k+1}(u,v)},\dot e^{m_{k+1}(u,v)})\to
(\Gal_{k+1}(u,v),\StamGal_{k+1}(u,v)).
\]
\qed}

{\em Proof of the theorem.} By the lemma above, there is a an attaching
map 
\[
\begin{array}{ccc}
(e^{m_k(u,v)},\dot e^{m_k(u,v)}) & \to & (\Gal_k(u,v),\StamGal_k(u,v)) \\
 & & \dto \\
 & & (\Cl\F_k(u,v),\Cl\F_{k-1}(u,v)\cup\Cl\F_{k-1}(v,u)).
\end{array}
\]
A similar construction works for the point space and the line space.\qed

\section{The Veronese imbedding, revisited}
\label{VeroneseRev}

The following fact was first proved by Knarr \cite{Kna2}. 

\Prop{If the $n$-gon $\frak P$ is a manifold (with parameters $(m,m')$), 
then the double mapping cylinder is a $(\dim\F+1)$-sphere.
Hence $\F$ can be imbedded in $\RR^{1+\dim\F}$ as a topological 
hypersurface with trivial normal bundle. The point space and the line
space can be imbedded in $\RR^{1+\dim\F}$ as topological submanifolds
with normal disk bundles of dimension $m'+1$, $m+1$, respectively.

\proof Put $E=D\F-\L$. Then $E\to\P$ is a locally trivial open disk 
bundle, and thus $E$ is a topological manifold. Being a homotopy
sphere of dimension $>3$, the double mapping cylinder is a sphere,
see \ref{Poincare}.\qed}

\Prop{\label{Stiefel-Whitney}
Let $\frak P$ be an $n$-gon which is a manifold. Suppose that
$n\neq 4$ or that $m=m'$. We get the following table for the total
Stiefel-Whitney classes $\bf w$ of $\P,\L,\F$, and of the sphere bundles
$\F\to\P$, $\F\to\L$ (cp. \cite{Brsp} for $n=3$).

\[
\begin{array}{l|ccccc}
      & {\bf w}(\P) & {\bf w}(\L) & {\bf w}(\F) & {\bf w}(\F\to\P) & {\bf w}(\F\to\L) \\ \hline
n=3,6 & 1+x_m+x_m^2 & 1+y_m+y_m^2 & 1           & 1+x_m            & 1+y_m \\
n=4   & 1           & 1+y_m       & 1           & 1                & 1+y_m 
\end{array}
\]

Note that the orientability of the manifolds and their normal disk bundles
is determined by the corresponding first Stiefel-Whitney classes ${\bf w}_1$.
Thus, the point space of an $(1,1)$-quadrangle is orientable, but
the line space is not (up to exchanging $\P$ and $\L$, of course).

In each case the number $m+1$ is the minimal codimension for an imbedding
of $\L$ into some euclidean space \cite[6 10.23]{Spanier}.

\proof The Stiefel-Whitney classes of $\P$ and $\L$ can be easily
calculated from \ref{SteenrodAlg} and
the Wu formula \cite[VI 17.]{Bredon93}. Since the
sums of the topological tangent bundles and the normal bundles of
$\P,\L,\F$ are trivial, the total Stiefel-Whitney classes of these
bundles are inverses of each other.\qed}

\section{Classification in small dimensions}

It was shown in \ref{LinesMnf} that the point rows and pencils of lines of
a finite-dimensional polygon with parameters $m,m'\leq 2$ are spheres.
Next, we want to classify the small Schubert varieties 
$\Cl\P_2(\ell,p)=p^\perp$ and $\Cl\L_2(p,\ell)=\ell^\perp$. This
generalizes results of Salzmann \cite{Sal67}
and Breitsprecher \cite{Brsp} about projective 
planes of small dimension.

\Thm{\label{T4_2_1}
Let $\frak P=\PLF$ be a compact, connected, finite dimensional
$n$-gon with equal parameters $m=m'\leq 2$. This implies by \ref{LinesMnf}
that $\frak P$ is a manifold. Let $p$ be a point and let $\ell$ be a 
line. There exist orthogonal $m$-plane bundles $\eta,\eta'$ over $L$ and
over $\L_p$, such that the sphere bundles $S(\eta\oplus\eps)$
and $S(\eta'\oplus\eps')$ are isomorphic to the sphere bundles
$\F_L\to L$, and $\F_{\L_p}\to\L_p$, respectively. The trivial line
bundles $\eps,\eps'$ correspond to the sections $q\to(q,\ell)$
and $h\to(p,h)$, respectively. Moreover, there are homeomorphisms
$E(\eta)\to p^\perp-\{p\}$, $E(\eta')\to\ell^\perp-\{\ell\}$ that
take the fibers of the vector bundles onto the punctured point rows
through $p$ and the punctured pencils of lines through $\ell$,
respectively. Note that $p^\perp$ and $\F_{\L_p}$ are
completely determined, once the vector bundle $\eta$ is known.

Let $\eta_\RR$ and $\eta_\CC$ denote the real and the complex Hopf
vector bundle over the $1$- and the $2$-sphere, respectively.
The Thom spaces of these bundles are of course the real and the
complex projective plane.

Up to duality, i.e. up to exchanging $\P$ and $\L$, 
precisely one of the following cases occurs:
\begin{enumerate}
\item $\frak P$ is a $2$-dimensional projective plane (thus $m=m'=1$), and
the bundles $\eta\cong\eta'\cong\eta_\RR$
are isomorphic to the real Hopf vector bundle. Thus
$\F_L\cong\F_{\L_p}$ is the Klein bottle.
\item $\frak P$ is a $4$-dimensional projective plane (thus $m=m'=2$), and 
the bundles $\eta\cong\eta'\cong\eta_\CC$ are isomorphic to the complex Hopf vector bundle.
\item $\frak P$ is a $3$-dimensional quadrangle  (thus $m=m'=1$). Then
$\eta\cong\eta_\RR$ is isomorphic to the real Hopf vector bundle,
and $\eta'$ is the trivial line
bundle over $\SS^1$. Hence $\F_{\L_p}$ is the Klein bottle, and
$\F_L$ is a $2$-torus.
\item $\frak P$ is a $6$-dimensional quadrangle (thus $m=m'=2$). Then
$\eta\cong\eta_\CC$ is isomorphic to the complex Hopf vector bundle,
and $\eta'$ is isomorphic to the tangent bundle of the $2$-sphere.
\item $\frak P$ is a $5$-dimensional hexagon (thus $m=m'=1$). Then
the bundles $\eta\cong\eta'\cong\eta_\RR$ are isomorphic to the real Hopf vector bundle.
Thus $\F_L\cong\F_{\L_p}$ is the Klein bottle.
\item $\frak P$ is a $10$-dimensional hexagon (thus $m=m'=2$). Then
$\eta\cong\eta_\CC$ is isomorphic to the complex Hopf vector bundle, and
$\eta'\cong\eta_\CC^{\otimes 3}$ is isomorphic to the unique
$2$-plane bundle over $\SS^2$ that is obtained by pulling back
the complex
Hopf vector bundle by a map $\SS^2\to\SS^2$ of degree $3$.
\end{enumerate}

\proof Consider the locally trivial $\SS^m$-bundle $\F_{\L_p}\to\L_p$.
By \ref{Kneser} the bundle $\F_{\L_p}\to\L_p$ is the sphere
bundle of an orthogonal vector bundle $\xi$. The map $\ell\mapsto (p,\ell)$ 
provides a section. Let $\eps$ denote the corresponding trivial line 
bundle, and put $\eta=\eps^\perp$. Now the total space of $\eta$ is 
homeomorphic to $\F_{\L_p}-(\{p\}\times\L_p)$ (by a stereographic
projection in the fibers), which in turn is homeomorphic to 
$p^\perp-\{p\}$. Thus the Thom space of $\eta$ has the 
same cohomology as $p^\perp$, which is known
by \ref{LinkCyc} and \ref{MueThm}.

In the case $m=m'=1$, there are just two line bundles over $\SS^1$,
namely the trivial bundle and the M\"obius bundle. For $n=3,6$, we have 
$\H^\bullet(\ell^\perp;\ZZ_2)\cong
 \H^\bullet(p^\perp;\ZZ_2)=\ZZ_2[x_1]/(x_1^3)$, whereas for $n=4$, we 
have $\H^\bullet(p^\perp;\ZZ_2)=\ZZ_2[x_1]/(x_1^3)$, and 
$\H^\bullet(\ell^\perp;\ZZ_2)=\ZZ_2[x_1,x_2]/(x_1^2,x_2^2,x_1x_2)$.
This completes the case $m=m'=1$.

The oriented $2$-plane bundles over $\SS^2$ are classified by their Euler
class, cp. \ref{LineBundle}. For $n=3$, we have 
$\H^\bullet(\L)\cong\H^\bullet(\P)=\ZZ[x_2]/(x_2^3)$.
For $n=4$, we have $\H^\bullet(p^\perp)=\ZZ[x_2]/(x_2^3)$, and 
$\H^\bullet(\ell^\perp;\ZZ)=\ZZ[x_2,x_4]/(x_4^2,x_4-2x_2^2,x_2x_4)$.
For $n=6$, we have $\H^\bullet(p^\perp)=\ZZ[x_2]/(x_2^3)$, and 
$\H^\bullet(\ell^\perp;\ZZ)=\ZZ[x_2,x_4]/(x_4^2,x_4-3x_2^2,x_2x_4)$.

Now we may apply \ref{EulerClass} to obtain the Euler class of the 
bundles $\eta$, $\eta'$.\qed}

The flag spaces of low-dimensional projective planes have been
classified by Breitsprecher \cite{Brs2} and Buchanan \cite{Buc}.
With the aid of Knarr's imbedding theorem, the proof becomes
considerably simpler.

\Thm{\label{BrspBuc}
Let $\frak P=\PLF$ be a compact connected
$2m$-dimensional projective plane. If $m\in\{1,2\}$, then the $m$-sphere
bundle $\F\to\P$ is isomorphic to its classical counterpart.

\proof By \ref{Kneser}, there exists an orthogonal $(m+1)$-plane
bundle $\xi$ over $\P$ that has $\F\to\P$ as its
sphere bundle.

If $m=1$, the total Stiefel-Whitney class $\bf w$ of $\xi$ is given by 
$1+x_1$, cp. \ref{Stiefel-Whitney}.
The bundle $\xi$ has no section by \ref{NoSection}, and $\xi$ is
not a sum of two copies of the tautological bundle $\eta$ over
$\P\cong\P_2\RR$, because ${\bf w}(\eta\oplus\eta)=1+x_1^2$,
cp. \cite[p.43]{MS}. Thus $\xi$ does not split, see \ref{VectP2R}.
Consider the two-fold
covering $f:\SS^2\to\P$. It follows from \ref{VectP2R} that $\xi$ is
uniquely determined by its pullback $f^{*}\xi$, and hence by
the fundamental group of the sphere bundle $S(f^{*}\xi)$, see 
\cite[26.2]{Steenrod}.
Since $\pi_1\F=Q_8$, we conclude that $\pi_1S(f^{*}\xi)=\ZZ_4$.

Now suppose that $m=2$.
According to \cite[3.]{DW}, we have only to 
calculate the first Pontrjagin class and the second Stiefel-Whitney 
class of $\xi$ in order to classify $\F$. The total Pontrjagin class
of $\P_2\CC$ is $1+3x_m^2$ \cite[p.178]{MS}. By \ref{TopTrivial},
we get  $\p(\xi)=1-3x_m^2$. The total Stiefel-Whitney class
of $\xi$ is $1+x_m$ by \ref{Stiefel-Whitney}.\qed}

Finally, we want to calculate the fundamental groups for the
case $m=m'=1$.

\Prop{\label{11FundGrp}
Let $\frak P$ be a finite-dimensional polygon with parameters
$(1,1)$. Then the fundamental groups of $\P$, $\L$, and $\F$ are 
as given in \ref{P3_2_6}.

\proof Since $m,m'\leq 2$, the polygon is a manifold, and we may
use the CW decomposition of \ref{CWPoly}. The $2$-skeleton of the point space
is $p^\perp$, and the $2$-skeleton of the line space is $\ell^\perp$.
Thus we get the fundamental groups of $\P$ and $\L$ by \ref{T4_2_1}.

The $1$-skeleton of the flag space is $\{p\}\times\L_p\cup
L\times\{\ell\}\cong\SS^1\vee\SS^1$. The $2$-skeleton of $\F$
is $\F_{\L_p}\cup\F_L$, and $\F_{\L_p}\cap\F_L=
\{p\}\times\L_p\cup L\times\{\ell\}$. Consider the diagram
\[
\begin{array}{ccccc}
 & & \{p\}\times\L_p\cup L\times\{\ell\} \\
 & \swarrow & & \searrow \\
\F_{\L_p} & & & & \F_L. 
\end{array}
\]
Passing to the fundamental groups, we get a diagram
\[
\begin{array}{ccccc}
 & & \bra{a,b} \\
 & \swarrow & & \searrow \\
\bra{a,b|\ R_1} & & & & \bra{a,b|\ R_2}
\end{array}
\]
where $R_1=aba^{-1}b$ and $R_2=bab^{-1}a$ for $n=3,6$, and
$R_1=aba^{-1}b$ and $R_2=bab^{-1}a^{-1}$ for $n=4$ by \ref{T4_2_1}.

By the Seifert-Van Kampen Theorem
we get the result stated in \ref{P3_2_6}.\qed}

\section{A Lemma on partial polygons}

Suppose we are given an incidence structure $\frak P=\PLF$ which
is a partial $n$-gon. We want to show that if $\P$ and $\L$ carry nice 
topologies such that $\PLF$ 'looks like a compact $n$-gon', then it is 
indeed a compact $n$-gon. The next lemma makes this more precise.

\Lemma{Let $\PLF$ be a thick partial $n$-gon. Suppose that $\P,\L,$ and
$\F$ are compact connected manifolds, and that $\pr_1$ and $\pr_2$ are 
locally trivial bundle maps, with $m'$- and $m$-manifolds as fibers. If 
the dimension of $\frak P$ is right, that is, if $\dim\F=n(m+m')/2=
m_n(p,\ell)=m_n(\ell,p)$,
then $\frak P$ is a compact connected $n$-gon, and in
fact a manifold.

\proof Let $(p,\ell)\in\F$ be an arbitray flag. Since
$\pr_1$, $\pr_2$ are locally trivial bundles, the iterated
bundle $\Gal_n(p,\ell)$ is a compact $n(m+m')/2$-dimensional manifold.
The subset $\PropGal_n(p,\ell)$ is nonempty and open, so we may choose
a connected component $M$ of $\Gal_n(p,\ell)$ that contains an open 
nonempty set $U$ of non-stammering galleries. Since $\frak P$ is a
partial $n$-gon, the map $(x_0,\ldots,x_{n+1})\mapsto \flag(x_n,x_{n+1})$ 
maps $U$ injectively onto its image in $\F$, and the preimage of the image
of $U$ is precisely $U$. By \ref{LocSur}, the map is surjective.
Thus $d(\ell,z)\leq n$ for every vertex $z\in\V$. Similarly, $d(p,z)\leq n$ 
for every $z\in\V$. Since $(p,\ell)$ was an arbitrary flag, the diameter 
of $\frak P$ is $n$, and thus $\frak P$ is a generalized $n$-gon. By 
\ref{P2_5_3}, it is a compact $n$-gon, and hence a manifold.\qed}

The topological assumptions of this lemma are certainly satisfied, if
$\frak P$ is a compact coset geometry of the right dimension of a
connected Lie group, or if 
$\frak P$ comes from an isoparametric hypersurface. Thus, the main
property that has to be checked for this kind of geometries in order
to show that they are compact $n$-gons is the {\em non-existence} of 
ordinary $k$-gons for $2\leq k<n$.

\chapter{Point homogeneous polygons}
\label{PHomPoly}

In this chapter we consider compact polygons with point transitive
automorphism groups.  Point homogeneous compact 
projective planes of dimension
at most 16 have been classified by Salzmann; it turns out that these
planes are precisely the four classical compact connected
Moufang planes \cite{Sal75}.
Later, L\"owen showed that the dimension of a point homogenous (and
thus finite-dimensional) compact projective plane is at most 16 \cite{Loew}.
L\"owen gave also a second proof for the classification (for $m>1$)
which made strong use of the Borel-De Siebenthal classification of maximal
subgroups of maximal rank in compact simple Lie groups and the fact that the Euler characteristic of a projective plane is $3$ \cite{Loew1}.
His proof uses some facts about involutive automorphisms of
projective planes; however, it turns out that the differential-geometric
properties of the group action are already 'good enough' to reconstruct
the plane uniquely from the group.

Generalized hexagons with $m>1$ have Euler characteristic $6$, and
generalized quadrangles with $m=m'>1$ have Euler characteristic $4$
(generalized quadrangles with $m\neq m' $ have Euler characteristic $0$). Thus
its seems reasonable to try to classify the point homogenous polygons of
positive Euler characteristic. 

On the other hand, if both parameters $m=m'=1$, then $n$ is an upper
bound for the dimension of compact automorphism groups, and this fact
leads to a classification of the point homogeneous polygons with
parameters $(1,1)$.

The main result of this chapter is Theorem \ref{T5_0_7}: If a compact
connected $n$-gon $\frak P$ admits a point transitive group of automorphisms
$\Gamma$, and if $n\neq 4$ or if $m=m'>1$, then $\Gamma$ is transitive
on the flag space $\F$. By the main result of \cite{GKK}, this implies that $\frak P$
is classical. In the case of quadrangles with parameters $m=m' =1$,
a point transitive group need not be flag transitive; however, the
existence of a point transitive group of automorphisms implies that
the quadrangle is the real symplectic quadrangle (\ref{4(11)Homogen}).
Thus we obtain a complete
classification of point homogeneous $n$-gons, provided that $n\neq 4$ or
that $m=m'$ (\ref{PHomClass}).

The notation for compact Lie groups is adapted from \cite{GKK}; thus
${\bf SO}_k\RR$ is the same group as ${\bf SO}(k)$ (which is used in the
other chapters).

\section{Compact transformation groups}

\Lemma{\label{L5_0_1}
Suppose that ${\bf U}_k\CC$ acts on $\SS^2$. If the action is
not transitive, then there is a fixed point. If $k>2$, then there
is always a fixed point.

\proof If the action is not transitive, then the orbits are circles
or points. Thus, ${\bf SU}_k\CC<{\bf U}_k\CC$ is contained in the kernel
of the action, and we have in fact an $\SS^1$-action on $\SS^2$. By 
\cite[4.7.12]{Spanier},
this action has a fixed point. For $k>2$, there is no
subgroup of codimension $2$ and rank $k-1$ in ${\bf U}_k\CC$.\qed}

\Lemma{\label{L5_0_2}
Every action of ${\bf SO}_5\RR\times{\bf SO}_2\RR$ on a
$2$-sphere has a fixed point.

\proof The action cannot be transitive, since otherwise the isotropy
group would have the same rank as ${\bf SO}_5\RR\times{\bf SO}_2\RR$. But the largest
subgroup of rank $3$ is clearly ${\bf SO}_4\RR\times{\bf SO}_2\RR$. Hence ${\bf SO}_5\RR$
is contained in the kernel of the action, and we get an ${\bf SO}_2\RR$
action on $\SS^2$. Thus there is a fixed point by \cite[4.7.12]{Spanier}.\qed}

\Thm{\label{T5_0_3}{\bf(Szenthe)} \cite{Szenthe} \cite[2.2]{GKK}
If a locally compact, connected,
second countable group $\Gamma$ acts as an effective
and transitive tranformation group on a connected, locally contractible
space $X$, then $\Gamma$ is a Lie transformation group, and 
$X\cong\Gamma/\Gamma_x$ is a (smooth) manifold.\qed}

\Thm{\label{T5_0_4} \cite{Mon}
Let $X=\Gamma/\Gamma_x$ be a compact connected homogeneous space of a Lie group
$\Gamma$. If the fundamental group of $X$ is finite, then every maximal
compact connected subgroup of $\Gamma$ acts transitively on $X$.\qed}

\Thm{\label{T5_0_5}
Let $X$ be a simply connected, compact manifold, with positive Euler
characteristic $\chi(X)=n\leq 6$. Let $\Delta$ be a compact, connected, effective and transitive Lie transformation group on $X$, with isotropy
group $\Delta_x$.
The fact that $\chi(X)>0$ implies that $\Delta_x$ has the same rank as
$\Delta$, and that $\Delta_x$ is connected. Moreover, the group
$\Delta$ is a product of centerless and simple compact Lie groups,
and $\Delta_x$ factors accordingly \cite{Wng49}.

There are precisely the following possibilities (we list only the Lie algebras
of the groups):

\medskip
$\bf(n=2)$ The group $\Delta$ is centerless and simple, and $\Delta_x$ is a
maximal connected subgroup. There are only the following possibilities:
\[
\begin{array}{c|c|cl}
\Delta & \Delta_x & \dim X & \\ \hline
\bb_k  & \dd_k    & 2k & X=\SS^{2k} \\
\gg_2  & \aa_2    & 6  & X=\SS^6 
\end{array}
\]

\medskip
$\bf(n=3)$ The group $\Delta$ is centerless and simple, and $\Delta_x$ is a
maximal connected subgroup. There are only the following possibilities:
\[
\begin{array}{c|c|cl}
\Delta & \Delta_x & \dim X & \\ \hline
\aa_2  & \aa_1+\RR & 4 & X=\P_2\CC \\
\cc_3  & \cc_1+\cc_2 & 8 & X=\P_2\HH \\
\ff_4  & \bb_4   & 16 & X=\P_2\OO \\
\gg_2  & \aa_1+\aa_1 & 8  
\end{array}
\]

\medskip
$\bf(n=4)$ If $\Delta$ is not simple, then we have a product of two
transformation groups listed in $(n=2)$.

Otherwise, $\Delta$ is centerless and simple. If $\Delta_x$ is a
maximal connected subgroup, then there are only the following possibilities:
\[
\begin{array}{c|c|cl}
\Delta & \Delta_x & \dim X & \\ \hline
\aa_3  & \aa_2+\RR & 6 & X=\P_3\CC \\
\cc_2  & \aa_1+\RR & 6  \\
\cc_4  & \cc_1+\cc_3 & 12 & X=\P_3\HH 
\end{array}
\]

Otherwise, there is a maximal connected subgroup $\Phi$
with $\Delta>\Phi>\Delta_x$, and $\Delta/\Phi$, $\Phi/\Delta_x$ are
even-dimensional spheres. There is only the following possibility:
\[
\begin{array}{c|c|c|c}
\Delta & \Phi & \Delta_x & \dim X \\ \hline
\cc_2  & \aa_1+\aa_1 & \aa_1+\RR & 6
\end{array}
\]

\medskip
$\bf(n=5)$ The group $\Delta$ is centerless and simple, and $\Delta_x$ is a
maximal connected subgroup. There are only the following possibilities:
\[
\begin{array}{c|c|cl}
\Delta & \Delta_x & \dim X & \\ \hline
\aa_4  & \aa_3+\RR & 8 & X=\P_4\CC \\
\cc_5  & \cc_1+\cc_4 & 16 & X=\P_4\HH 
\end{array}
\]

\medskip
$\bf(n=6)$ If $\Delta$ is not simple, then we have a product of a
transformation group listed in $(n=2)$ with a transformation group
listed in $(n=3)$.

Otherwise, $\Delta$ is centerless and simple. If $\Delta_x$ is a
maximal connected subgroup, then there are only the following possibilities:
\[
\begin{array}{c|c|cl}
\Delta & \Delta_x & \dim X & \\ \hline
\aa_3  & \aa_1+\aa_1+\RR & 8 \\
\aa_5  & \aa_4+\RR & 10 & X=\P_5\CC \\
\bb_3  & \bb_1+\dd_2 & 12 \\
\bb_3  & \bb_2+\RR   & 10 \\
\cc_4  & \cc_2+\cc_2 & 24 \\
\cc_6  & \cc_1+\cc_5 & 20 & X=\P_5\HH 
\end{array}
\]

Otherwise, there is a maximal connected subgroup $\Phi$
with $\Delta>\Phi>\Delta_x$, and $\Delta/\Phi$, $\Phi/\Delta_x$ are
among the spaces listed in $(n=2)$ and $(n=3)$ (note, however, that
$\Phi$ need not be effective on $\Phi/\Delta_x$).
There are only the following possibilities:
\[
\begin{array}{c|c|c|c}
\Delta & \Phi & \Delta_x & \dim X \\ \hline
\aa_2  & \aa_1+\RR & \RR+\RR & 6 \\
\cc_3  & \cc_1+\cc_2 & \cc_1+\cc_1+\cc_1 & 12 \\
\cc_3  & \cc_1+\cc_2 & \RR+\cc_2 & 10 \\
\ff_4  & \bb_4 & \dd_4 & 24 \\
\gg_2  & \aa_2 & \aa_1+\RR & 10 \\
\gg_2  & \aa_1+\aa_1 & \aa_1+\RR & 10 \\
\gg_2  & \aa_1+\aa_1 & \RR+\aa_1 & 10
\end{array}
\]

\proof This follows from the Borel-De Siebenthal classification of
maximal subgroups of maximal rank. See \cite{Wng49},
\cite{BSie}, \cite[8.10]{Wo}, and
\cite[2.4 and Table 1]{GKK}.\qed}

\section{The classification}

\Thm{\label{T5_0_6}{\bf(Burns-Spatzier)} \cite[2.1]{BS}
Let $\frak P$ be a compact connected polygon. The group of all
continuous automorphisms of $\frak P$, endowed with the compact-open
topology, is a locally compact, second countable topological transformation
group.\qed}

\Lemma{\label{L5_0_7} \cite[3.2]{GKK}
Let $\frak P$ be a finite-dimensional polygon with parameters
$m=m'=1$. Suppose a compact Lie group $\Delta$ acts as an effective
automorphism group on $\frak P$. Then the isotropy group $\Delta_{p,\ell}$
of every flag $(p,\ell)$ is finite; hence $\dim\Delta\leq n$, and if 
$\dim\Delta=n$, then $\Delta$ is flag transitive.

It follows that the isotropy group $\Delta_p$ of a point (or a line) 
is at most one-dimensional, and if $\Delta_p$ is one-dimensional, then it is
transitive on the pencil of lines $\L_p$.

\proof Let $(p,\ell)$ be a flag. The connected component 
$(\Delta_{p,\ell})^1$ has a fixed point on $L\times\{\ell\}\cong\SS^1$ and
on $\{p\}\times\L_p\cong\SS^1$; thus it fixes every flag on
$\{p\}\times\L_p\cup L\times\{\ell\}$. But this implies that 
$(\Delta_{p,\ell})^1$ fixes every flag.\qed}

\Thm{\label{T5_0_7}
Let $\frak P$ be a compact, connected $n$-gon. Suppose that a closed
subgroup $\Gamma$ of the automorphism group of $\frak P$ acts transitively
on the point space $\P$. This implies that $\frak P$ is finite 
dimensional, so the topological parameters $(m,m')$ of $\frak P$ 
are defined.

If $n\neq 4$, or if $m=m'>1$, then the group $\Gamma$ is transitive on 
$\F$. By the main result (Theorem 3.8) of \cite{GKK}, the polygon
$\frak P$ is a compact connected Moufang polygon. 
(The point transitive groups are listed explicitly in \cite{GKK}.)

\proof (1) By \ref{T5_0_3}, the group $\Gamma$ is a Lie group, and $\P$ is a 
homogeneous space. Hence $\frak P$ is finite-dimensional by \ref{P3_0_9}. 

By \ref{P3_2_6}, the point space $\P$
has a finite fundamental group in the cases that we consider.
Thus, a compact connected subgroup $\Delta$ 
of $\Gamma$ is still transitive on $\P$ by \ref{T5_0_4}. 

\smallskip
(2) Let $p$ be a point. The isotropy group $\Delta_p$ fixes no line
$\ell\in\L_p$, because otherwise the map 
$\delta p\mapsto (\delta p,\delta\ell)$, with $\delta\in\Delta$,
would provide a section of the bundle $\F\to\P$, contradicting \ref{NoSection}.

\smallskip
(3) Now suppose that $m=m'=1$. Then $n=3,6$, and $\dim\P=n-1\leq \dim\Delta$.
If $\dim\Delta=n-1\in\{2,5\}$, then $\Delta$ has an infinite fundamental
group, and $\Delta\to\Delta/\Delta_p$ is a covering.
But the fundamental group of $\P$ is $\ZZ_2$. Thus $\dim\Delta=n$,
and by \ref{L5_0_7}, $\Delta$ is transitive on $\F$.

\smallskip
(4) The remaining cases are $n=3,4,6$, and $m=m'>1$.
Thus $\P$ is simply connected, and we may apply \ref{T5_0_5}, since
the Euler characteristic of $\P$ is $n=3,4,6$. Let $p$ be a point, 
and let $\Delta_p$ be the stabilizer of $p$.

\smallskip
If $n=3$, then $(\Delta,\Delta_p)\in\{
({\bf PSU}_3\CC,{\bf U}_2\CC),
({\bf PU}_3\HH,{\bf PU}_1\HH\times{\bf U}_2\HH),
({\bf F}_4,{\bf Spin}_9),
\linebreak
({\bf G}_2,{\bf SO}_4\RR)
\}$.
Since $\pi_2{\bf G}_2=0=\pi_1\P$, it follows
that $\pi_2({\bf G}_2/{\bf SO}_4\RR)\cong\pi_1{\bf SO}_4\RR=\ZZ_2$, 
and this excludes this space as a
candidate for an $8$-dimensional projective plane.
Thus the point space $\P$ is homeomorphic to the point space of one 
of the four classical compact connected planes. It is readily seen that in
each of the three cases, $\Delta$
contains no proper closed subgroups of codimension $\leq2m$
besides the conjugates of $\Delta_p$.
Thus, $\Delta$ is also transitive on the line space, and $\Delta_p$
fixes some line $h$. But $\Delta_p$ has only one orbit $X$ of dimension 
$\leq m$ in $\P-\{p\}$, namely the cut locus of $p$ (because $\P$ is
a compact symmetric space of rank $1$ under the action of $\Delta$),
and this set is precisely the polar line of $p$  with respect to
the elliptic polarity . Hence $H=X$ is a classical point row, and
therefore each point row of $\frak P$ is a point row of the 
classical plane. It is well known that the elliptic motion groups
of the four classical compact connected Moufang planes are flag transitive.

\smallskip
Now let $n=4$. It follows from the structure of $\H^\bullet(\P)$ that
$\P$ is not a product of two even-dimensional spheres (see also \ref{NoProduct}).
Thus we have to consider the simple groups listed in \ref{T5_0_5}.
The group $\Delta_p$ acts without fixed points on the homology $m$-sphere
$\L_p$ by (2). The point space $\P$ and the line space $\L$ are
$3m$-dimensional, and $m\in\{2,4\}$. 

The pair 
$({\bf PSU}_4\CC,{\bf U}_3\CC)$ is excluded by \ref{L5_0_1} and (2).

The group ${\bf PU}_4\HH$ has (up to conjugation) only one subgroup
of codimension $\leq 12$, as is easily seen. Thus, if ${\bf PU}_4\HH$
is point transitive, then it is line transitive as well, and $\P$
and $\L$ are homeomorphic. But the cohomology rings of $\P$ and $\L$
are not isomorphic, and this excludes the pair 
$({\bf PU}_4\HH,({\bf U}_3\HH\times{\bf U}_1\HH)/\pm\bf1)$.

In the two remaining cases, $\Delta_p={\bf U}_2\CC$. It follows
from (2) and \ref{L5_0_1} that $\Delta_p$ is transitive on $\L_p$,
and hence $\Delta={\bf U}_2\HH$ is transitive on $\F$. (There are two
possibilities for the imbedding $\Delta_p\SUB\Delta$, corresponding to
the complex symplectic quadrangle and its dual.)

\smallskip
Finally, let $n=6$. Then $\dim\P=\dim\L=5m\in\{10,20\}$.
It is clear from the structure of $\H^\bullet(\P)$ that
the action of $\Delta$ does not factor as a product
(see also \ref{NoProduct}),
so we have again only to consider the simple groups listed in 
\ref{T5_0_6}. 

The groups ${\bf U}_6\CC$ and ${\bf SO}_5\RR\times{\bf SO}_2\RR$
are excluded as canditates for $\Delta_p$ by (2) and \ref{L5_0_1}.

The groups ${\bf PU}_3\HH$ and ${\bf PU}_6\HH$ have only one conjugacy
class of subgroups of codimension $\leq 10$ and $\leq 20$, respectively.
Hence, if one of these groups is point transitive, then it is line transitive,
and $\P$ and $\L$ are homeomorphic. But the cohomology rings of
$\P$ and $\L$ are not isomorphic.

Thus, $\Delta={\bf G}_2$ and $\Delta_p={\bf U}_2\CC$. It follows
as in the case $(n=4)$ from (2) and \ref{L5_0_1} that $\Delta$
is flag transitive (there are again two imbeddings
${\bf U}_2\CC\SUB{\bf G}_2$, due to the fact that the complex
hexagon is not self-dual).\qed}

The group ${\bf SO}_3\RR$ acts as a sharply 
point transitive group on the real symplectic quadrangle. Thus we
have to modify our assumptions for $(1,1)$-quadrangles.

\Thm{Let $\frak P$ be a finite-dimensional quadrangle with parameters 
$(1,1)$. If its automorphism group contains a  
semi-simple compact subgroup $\Delta$ of positive dimension, 
then $\frak P$ is (up to duality) the real symplectic quadrangle.

\proof (cp. the proof of \cite[3.4]{GKK} due to Knarr)
By \ref{L5_0_7} we have $\dim\Delta\leq 4$, and by assumption 
$\dim\Delta\geq 3$. Thus we may assume that $\Delta$ is connected and
of type $\aa_1$. Passing to the dual quadrangle, if necessary,  we may
assume moreover that $\pi_1\P=\ZZ_2$, and that $\pi_1\L=\ZZ$.

Let $\ell\in\L$ be a line. Since $\Delta$ cannot be transitive
on $\L$, the isotropy group $\Delta_\ell$ has positive dimension.
By \ref{L5_0_7} the isotropy group $\Delta_\ell$ is transitive on
the point row $L$. Since this is true for every point row,
$\Delta$ is transitive on $\P$, and thus $\Delta\cong{\bf SO}_3\RR$,
and $\P\cong\P_3\RR\cong{\bf SO}_3\RR$. Every point row in 
$\P$ is a translate of a $1$-parameter group in ${\bf SO}_3\RR$ and 
thus a classical point row in $\P_3\RR\cong{\bf SO}_3\RR$. By the result 
of Dienst \cite{Dienst} on quadrangles which are imbeddable in projective
spaces, $\frak P$ is the real symplectic quadrangle.\qed}

\Lemma{Let $\Gamma$ be a connected Lie group, and let $M=\Gamma/\Gamma_x$
be a homogeneous space of $\Gamma$, with fundamental group $\pi$.
Let $\Delta$ be a maximal compact subgroup of $\Gamma$. If $\Delta$
contains no almost-simple factor (i.e. if $\Delta$ is trivial or a torus),
then $M$ is an Eilenberg-MacLane space $K(\pi,1)$, and thus the manifold
$M$ has the same (co-) homology as the group $\pi$.

\proof Assume that $\Delta$ is a torus. Passing to the universal
covering group, we may assume that $\Gamma$,
as well as every closed connected 
subgroup of $\Gamma$, is contractible \cite[XV 3.1]{Hoch}.
It follows from the homotopy
sequence of $\Gamma_x\to\Gamma\to M$ that the homotopy groups of $M$
vanish in dimension $\geq 2$.\qed}

\Cor{Let $\frak P$ be a compact connected polygon. Let $\Gamma$ be
a Lie group that acts as a point transitive automorphism group on $\frak P$. Then
$\Gamma$ contains a compact semi-simple subgroup of 
positive dimension.

\proof Since $\P$ is homogeneous, $\frak P$ is finite-dimensional,
and $\pi_1\P$ is cyclic. The homology of the point space
is not isomorphic to the homology of a cyclic group, cp.
\cite[V 7.6]{Whitehead}.\qed}

\Cor{\label{4(11)Homogen}
Let $\frak P$ be a compact connected quadrangle with a 
point transitive automorphism group. If $\frak P$ has equal
parameters $m=m'=1$, then $\frak P$ is the real symplectic
quadrangle.\qed}

\Cor{\label{PHomClass}
Let $\frak P$ be a compact connected $n$-gon. Suppose that
a group of continuous automorphisms acts transitively on the point
space $\P$. Then $\frak P$ is finite dimensional. If $n\neq 4$, or
if $m=m'$, then $\frak P$ is (up to duality) one of the following
polygons:
\begin{itemize}
\item[$\bf3_\RR$] The real projective plane.
\item[$\bf3_\CC$] The complex projective plane.
\item[$\bf3_\HH$] The quaternion projective plane.
\item[$\bf3_\OO$] The octonion (Cayley) projective plane.
\item[$\bf4_\RR$] The real symplectic quadrangle.
\item[$\bf4_\CC$] The complex symplectic quadrangle.
\item[$\bf6_\RR$] The real Moufang hexagon (associated with ${\bf G}_{2(2)}$).
\item[$\bf6_\CC$] The complex Moufang hexagon (associated with ${\bf G}_2^\CC$).
\end{itemize}
\qed}

Ferus-Karcher-M\"unzner have constructed 
point homogeneous quadrangles (with $m\neq m'$) which are not Moufang \cite{FKM,GTH91}.

\chapter{Miscellanies}

\begin{flushright}
\em What I tell you three times is true. \\
L. Carroll, Hunting of the snark
\end{flushright}

\bigskip
In this chapter we collect some topological results which are
needed in various proofs, but which are not directly connected
with the geometry of a topological polygon.

\section{Cofibrations and ANRs}

\Def{Let $(X,A)$ be a topological 
pair. Recall that the inclusion $A\SUB X$ is
called a {\em cofibration}, if the extension problem
\[
\begin{array}{ccc}
A\times\II\cup X\times\{0\} & \too  & Y \\
\dto \\
X\times\II
\end{array}
\]
has a solution for every space $Y$
\cite[I]{tDKP}, \cite[I 5]{Whitehead},
\cite[II 3]{tDck}, \cite[VII 1]{Bredon93}. 
The cofibration is called {\em closed},
if $A$ is a closed subspace. If $X$ is a Hausdorff space, then $A$ is
necessarily closed \cite[1.17]{tDKP}. 

A pointed space $(X,*)$ is called {\em well-pointed}, if $\{*\}\SUB X$ is
a cofibration.

An inclusion $A\SUB X$ is a cofibration if and only if 
$A\times\II\cup X\times\{0\}$ is a retract of $X\times\II$
\cite[1.22]{tDKP}.}

\Prop{Consider the cocartesian diagram
\[
\begin{array}{ccc}
A    & \TO{f} & B            \\
\dto &        & \dto         \\
X    & \to    & X\cup_f B.
\end{array}
\]
If $A\SUB X$ is a cofibration, then the pushout $B\to X\cup_f B$ is also
a cofibration \cite[7.36]{tDKP}.\qed}

\Prop{\label{CollapseA}
If $A\SUB X$ is a closed cofibration, then $(X/A,*)$ is well-pointed.
If $A$ is contractible, then $(X,A)\to(X/A,*)$ is a homotopy equivalence
\cite[2.36]{tDKP}.\qed}

\Prop{\label{Collapse}
If $A\SUB X$ is a closed cofibration, then the collapsing map
\[
(X,A)\to(X/A,*)
\]
induces an isomorphism in homology and cohomology
\cite[VII 1.7]{Bredon93}.

In particular, if $(X,x)$ and $(Y,y)$ are well-pointed disjoint spaces, 
then the map $(X\cup Y,\{x,y\})\to(X\vee Y,*)$ induces an isomorphism
in (co-) homology, i.e. the (co-) homology of the one-point union
$(X\vee Y,*)$ is the direct sum of the (co-) homology of $(X,x)$ 
and $(Y,y)$.}

\Cor{\label{Wedge}
If $(X_1,x_1)$, $(X_2,x_2)$  are well-pointed homology and cohomology 
spheres of dimension $m$ and $m'$, respectively, then $(X\wedge Y,*)$ is
a homology and cohomology sphere of dimension $m+m'$.

\proof The subspaces $X_1\times\{x_2\},\{x_1\}\times X_2$
form an excisive couple, since the inclusion
$(X_1\times\{x_2\},\{(x_1,x_2)\})\SUB (X_1\vee X_2,\{x_1\}\times X_2)$
induces an isomorphism in singular homology by \ref{Collapse}, see
\cite[4.6.4]{Spanier}. By the K\"unneth Theorem \cite[5.3.10]{Spanier},
we get an isomorphism 
$[\H_\bullet(X_1,x_1)\otimes\H_\bullet(X_2,x_2)]_j\cong
\H_\bullet(X_1\times X_2,X_1\vee X_2)_j$.
The inclusion $X_1\vee X_2\SUB X_1\times X_2$ is a cofibration by
\cite[3.20]{tDKP}, hence we get an isomorphism
$\H_\bullet(X_1\times X_2,X_1\vee X_2)\cong\H_\bullet(X_1\wedge X_2,*)$.
Thus $X_1\wedge X_2$ is a homology $(m+m')$-sphere. It follows
from the universal coefficient theorem \cite[5.5.3]{Spanier} that
$X_1\wedge X_2$ is also a cohomology $(m+m')$-sphere.\qed}

\Thm{{\bf(Homotopy Excision Theorem)} \cite{Spa67} \label{HET}
Suppose $X=A\cup B$ is the union of two closed subspaces $A,B$. If
the intersection $A\cap B$ is a strong deformation retract of some
neighborhood in $A$, and if the pairs $(A,A\cap B)$ and $(B,A\cap B)$
are $n$- and $m$-connected, respectively, where $m,n\geq 0$,
then the map
\[
\pi_r(B,A\cap B)\to\pi_r(X,A)
\]
is an isomorphism for $0\leq r\leq n+m-1$, and an epimorphism for $r=m+n$.
\qed}

\Cor{\label{HETC}
Let $A\SUB X$ be a closed cofibration. If the pair $(X,A)$ is 
$m$-connected, and if $A$ is $n$-connected, for $m,n\geq 0$, then the
collapsing map induces isomorphisms 
\[
\pi_r(X,A)\to\pi_r(X/A)
\]
for $0\leq r\leq m+n$, and an epimorphism for $r=m+n+1$. In particular,
$X/A$ is $m$-connected.

\proof Let $C_A'$ denote the unreduced cone over $A$. We may apply
the Homotopy Excision Theorem \ref{HET} to the space $X\cup C_A'$. Since
$C_A'$ is contractible, the pair $(C_A',A)$ is $(n+1)$-connected. Thus
the inclusion $(X,A)\SUB(X\cup C_A',C_A')$ induces an isomorphism for
the homotopy groups in dimension $\leq n+m$, and an epimorphism
in dimension $m+n+1$. The collapsing map $(X\cup C_A',C_A')\to(X/A,*)$
is a homotopy equivalence by \ref{CollapseA}.\qed}

\Def{\label{D6_1_7}
A metrizable space $X$ is called an 
{\em absolute neighborhood retract} or ANR for short, if for
every imbedding $X\SUB M$ of $X$ as a closed
subspace into a metrizable space $M$,
there exists an open neighborhood $U$ of $X$ such that
$X$ is a retract of $U$ \cite[III 6]{Hu}.

Every locally contractible metrizable space of finite covering
dimension is an ANR \cite[V 7.1]{Hu}.

A closed subspace $A$ of an ANR $X$ is an ANR if and only if the
inclusion $A\SUB X$ is a cofibration \cite[IV 3.2]{Hu}.

Suppose that $X,A,B$ are ANR's, that $A\SUB X$ is closed, and that
$f:A\to B$ is continuous. If the pushout $X\cup_f B$ is metrizable,
then it is an ANR \cite[VI 1.2]{Hu}.

In particular, suppose that $Z=A\cup B$ is the union of two closed subsets.
If $Z$ is metrizable, and if $A$, $B$, and $A\cap B$ are ANR's,
then $Z$ is an ANR.

Any ANR is homotopy equivalent to a CW complex \cite[p.218]{Web}.\qed}

\section{Bundles and homotopy}

\Prop{Let $X,Y$ be pointed spaces, and let $[X;Y]^0$ denote the set
of all base-point preserving homotopy classes of maps from
$X$ to $Y$. If $X$ is well-pointed, then there is a natural
action of $\pi_1Y$ on $[X;Y]^0$ \cite[7 3.4]{Spanier}. If $\,Y$ is 
path-connected, then the space $[X;Y]$ of all free (i.e. not necessarily
base-point preserving) homotopy classes is the orbit space 
$[X;Y]^0/\pi_1Y$ \cite[III 1.11]{Whitehead} \cite[5.10.8]{tDck}.

The pointed space $Y$ is called {\em simple}, if it is path-connected, and 
if the $\pi_1Y$-action on $[X;Y]^0$ is trivial for all well-pointed
spaces $X$. Thus we have $[X;Y]^0\cong[X;Y]$ for simple spaces $Y$.
Every path-connected $H$-space (in particular every Eilenberg-MacLane
space) is simple \cite[7 3.5]{Spanier} \cite[III 4.18]{Whitehead}.\qed} 

\Prop{Suppose $Y,Y'$ are well-pointed spaces. Let $f:Y\to Y'$ be
a base-point preserving map. If $f$ induces an isomorphism
between all homotopy groups of $Y$ and $Y'$,
then $f$ induces an isomorphism $[X;Y]^0\cong[X;Y']^0$ for every 
well-pointed CW complex $X$ \cite[7 7.14,7.15]{Spanier}.\qed}

\Prop{Given a topological group $G$, let ${\bf B}G$ denote its
classifying space (unique up to homotopy equivalence).
If $X$ is a paracompact space, then the set of homotopy classes
$[X;{\bf B}G]$ is in one-to-one correspondence with the set of
all isomorphism classes of principal $G$-bundles over $X$. The total
space of the universal principal $G$-bundle is contractible, thus
$\pi_jG\cong\pi_{j+1}{\bf B}G$ \cite[I.4.12]{Hus}, \cite{Dold2}.\qed}

\Prop{\label{LineBundle}
Let $X$ be a paracompact space having the homotopy type of a
CW complex. The map that assigns to a line bundle over $X$ its first
Stiefel-Whitney class ${\bf w}_1$ is a bijection between the isomorphism
classes of line bundles over $X$ and $\H^1(X;\ZZ_2)$.
The map that assigns to an oriented $2$-plane bundle over $X$
its Euler class $\bf e$ is a bijection between the isomorphism classes of
oriented $2$-plane bundles over $X$ and $\H^2(X)$.

\proof The classifying spaces ${\bf BO}(1)$ and ${\bf BSO}(2)$
are Eilenberg-MacLane spaces of type $K(\ZZ_2,1)$ and $K(\ZZ,2)$,
respectively. We have to check that the maps
\[
[X;{\bf BO}(1)]\TO{{\bf w}_1}[X;K(\ZZ_2,1)]
\]
and 
\[
[X;{\bf BSO}(2)]\TO{\bf e}[X;K(\ZZ,2)]
\]
are bijections.

This is true if $X$ is a sphere, see eg. \cite[IV 4.10]{tDck}.
In case that $X$ is a CW complex, we get by the above remarks
a commutative diagram
\[
\begin{array}{ccc}
[X;{\bf BSO}(2)] & \OT\cong & [X;{\bf BSO}(2)]^0 \\
\DTO{\bf e}      &       & \DTO\cong \\
{[X;K(\ZZ,2)]}   & \OT\cong & [X;K(\ZZ,2)]^0,
\end{array}
\]
and a similar diagram for ${\bf BO}(1)$.\qed}

\Prop{\label{EulerClass}Let $\xi$ be an oriented
$2k$-plane bundle over a 
$2k$-sphere $\SS^{2k}$. Let $x$ denote a generator of $\H^{2k}(\SS^{2k})$, 
and let $\e_\xi=ax$ denote the Euler class of $\xi$. Then the 
cohomology ring of the Thom space $T_\xi$ of $\xi$ is given by 
\[
\H^\bullet(T_\xi)=\ZZ[u,v]/(v^2,u^2-av,uv),
\]
where $u$ and $v$ generate 
$\H^{2k}(T_\xi)$ and $\H^{4k}(T_\xi)$, respectively.

\proof Let $U_\xi\in\H^{2k}(E,E_0)$ be the orientation class 
of $\xi$, and let $\theta:\H^\bullet(\SS^{2k})\to\H^\bullet(E,E_0)$ 
denote the Thom isomorphism. We have the relation $\theta(1)=U_\xi$, and
$\theta(\e_\xi)=U_\xi^2$ \cite[p. 99]{MS}. On the other hand, 
$\H^{2k}(E,E_0)$ is generated by
$\theta(1)$, and $\H^{4k}(E,E_0)$ is generated by $\theta(x)$.
Finally, the reduced cohomology of $T_\xi$ is naturally isomorphic 
to that of $(E,E_0)$.\qed}

\Prop{\label{VectP2R}\cite{Levine} \cite[5.13 ]{tDck}
Let $\xi$ be a $2$-plane bundle over the real projective plane $\P_2\RR$.
Let $\eta$ denote the tautological line bundle over $\P_2\RR$, and let
$f:\SS^2\to\P_2\RR$ denote the two-fold covering of $\P_2\RR$.

If $\xi$ splits as a sum of two line bundles, then $\xi$ is isomorphic to
one of the bundles $\eps^2,\eps\oplus\eta,\eta\oplus\eta$.
(This follows from \ref{LineBundle} and the fact that 
$\H^1(\P_2\RR;\ZZ_2)=\ZZ_2$.)

Otherwise, $\xi$ is uniquely determined by its pullback $f^{*}\xi\in
{\bf Vect}^2(\SS^2)$, and hence by the fundamental group of the
sphere bundle $S(f^{*}\xi)$, cp. \cite[26.2]{Steenrod}.\qed}

\Thm{\label{Kneser}{\bf(Kneser)} \cite{Kne} \cite{Fri} \cite[V \S5]{KS}
Let ${\bf Top}(\SS^m)$ denote the group 
of all homeomorphisms of the $m$-sphere. If $m\leq 2$, then the inclusion 
${\bf O}(m+1)\to{\bf Top}(\SS^m)$ is a homotopy equivalence. Therefore
the structure group of any locally trivial $m$-sphere bundle can be 
reduced to ${\bf O}(m+1)$ in these dimensions, that is, there exists 
always an orthogonal $(m+1)$-plane bundle that has this sphere bundle 
as its unit sphere bundle.\qed}

\Thm{\label{Novikov}
Let $\xi,\xi'$ be two $k$-plane bundles over a paracompact space
$X$. If the underlying topological $\RR^k$-bundles are isomorphic, then
the rational Pontrjagin classes of $\xi$ and $\xi'$ are the same.

\proof The rational cohomology of ${\bf BO}(k)$ is a polynomial
ring generated by the Pontrjagin classes of the universal $k$-plane
bundle \cite[p.182]{MS}. Since the Pontrjagin classes are stable
\cite[15.2]{MS}, the cohomology of ${\bf BO}$ is a polynomial ring
generated by the universal Pontrjagin classes, and 
$\H^\bullet({\bf BO}(k);\QQ)\ot\H^\bullet({\bf BO};\QQ)$ is
an epimorphism. 

The homotopy fiber ${\bf TOP/O}$ is simply connected and
has finite homotopy groups in every dimension \cite[V \S5]{KS}. Thus
${\bf\tilde\H}{}^\bullet({\bf TOP/O};\QQ)=0$ by 
\cite[9 6.15, 5 2.8, 5 5.3]{Spanier}.
Hence the spectral sequence of the fibration ${\bf BO}\to{\bf BTOP}$
collapses, and we get an isomorphism
$\H^\bullet({\bf BO};\QQ)\OT\cong\H^\bullet({\bf BTOP};\QQ)$,
cp. also \cite[Epilogue]{MS}.

Let $f,f'$ be classifying maps
for $\xi$ and $\xi'$, respectively. The classifying maps of the
underlying $\RR^m$ bundles
$X\TO f  {\bf BO}\to{\bf BTOP}(k)$ and 
$X\TO{f'}{\bf BO}\to{\bf BTOP}(k)$ are homotopic,
and the claim follows from the homotopy-commutative diagram
\[
\begin{array}{ccccc}
X & \to & {\bf BO}(k)   & \to & {\bf BO} \\
  &     & \dto          &     & \dto \\
  &     & {\bf BTOP}(k) & \to & {\bf BTOP}.
\end{array}
\]
\qed}

\Cor{\label{NovikovCor} {\bf(Novikov's Theorem)}
Let $M$ be a smooth $k$-manifold. The rational Pontrjagin classes
of $M$ do not depend on the differentiable structure, but
only on the topological structure of $M$.

\proof By Kister's result \cite{Kis}, the underlying $\RR^k$ bundle
of the tangent bundle of $M$ is independent of the differentiable 
structure of $M$.\qed}

\Cor{\label{TopTrivial}
Let $\xi$ be a vector bundle over a smooth manifold $M$. If
the total space of $\xi$ is homeomorphic to an open subset
of some euclidean space, then
the total rational Pontrjagin classes of $M$ and of $\xi$ are inverse
to each other, i.e. 
\[
{\bf p}(M){\bf p}(\xi)=1.
\]

\proof We may as well assume that $\xi$ is smooth \cite[4 3.5]{Hir}.
The total Pontrjagin class of the total space $E(\xi)$ of $\xi$ vanishes
by \ref{NovikovCor}. Since $\xi$ is a normal bundle of $M\hookrightarrow
E(\xi)$, it follows that $1={\bf p}(TM\oplus\xi)={\bf p}(M)
{\bf p}(\xi)$, cp. \cite[15.3]{MS}.\qed}

\Prop{\label{OddBundle}
Let $E\to B$ be a locally trivial, orientable bundle. Assume that
the fiber $F$ is a rational homology $n$-sphere. If one of the 
maps $\H_n(F;\QQ)\to\H_n(E;\QQ)$ or $\H^n(F;\QQ)\ot\H^n(E;\QQ)$ is
trivial, then $n$ is odd.

\proof It follows from the universal coefficient theorem for field coefficients that $F$ is a rational cohomology $n$-sphere, and that the map
$\H^n(F;\QQ)\ot\H^n(E;\QQ)$, being the adjoint of $\H_n(F;\QQ)\to\H_n(E;\QQ)$,
is trivial, if the latter map is trivial \cite[\S 53]{Munkres}.

Suppose it is.
Consider the cohomology spectral
sequence of $E\to B$ \cite[9 5]{Spanier}. The map $\H^n(E;\QQ)\to\H^n(F;\QQ)$
is the composite
\begin{eqnarray*}
\lefteqn{\H^n(E;\QQ)=F^0\H^n(E;\QQ)\to F^0\H^n(E;\QQ)/F^1\H^n(E;\QQ)} \\
& = & E_\infty^{0,n}\to E_2^{0,n}=\H^0(B;\H^n(F;\QQ))\cong\H^n(F;\QQ)
\end{eqnarray*}
\cite[9 5]{Spanier}, hence the map $E_\infty^{0,n}\to E_2^{0,n}$ is trivial. 

Now $E_2^{s,t}=\H^s(B;\H^t(F;\QQ))=0$ unless $t\in\{0,n\}$, since the fiber
$F$ is a rational cohomology $n$-sphere. Therefore,
all differentials of the spectral sequence except for $E_2^{s,n}\TO{d_{n+1}}
E_2^{s+n+1,0}$ vanish. Since $E_{n+2}=E_\infty$ is the cohomology of this
cochain complex, there is an exact sequence
\[
0\to E_\infty^{0,n}\to E_2^{0,n}\to E_2^{n+1,0}\to E_\infty^{n+1,0}\to 0.
\]
It follows that the horizontal arrows in the diagram
\[
\begin{array}{ccc}
E_2^{0,n} & \to & E_2^{n+1,0} \\
\DTO\cong &     & \DTO\cong   \\
\H^0(B;\QQ) & \TO{\cup\Omega} & \H^{n+1}(B;\QQ)
\end{array}
\]
are monomorphisms, and thus $\Omega\neq 0\neq2\Omega$; hence
$n$ is odd \cite[9 5.2]{Spanier}.\qed}

\Lemma{\label{CWBundle}
Let $B$ be a finite CW complex, and let $p:E\to B$ be a locally 
trivial $m$-sphere bundle with a section $s:B\to E$. Then there is the
structure of a CW complex on $E$, with $s(B)$ as a subcomplex. For each
cell $\iota:e^k\to B$, the set $E_{\iota(e^k)}$ consists of an
$(m+k)$-cell.

\proof We proceed by induction on $k=\dim B$. If $B$ consists of a
finite number of points, then there is nothing to show. Now
assume that the claim holds for sphere bundles over
$k$-dimensional CW complexes. Let $B^{(k)}$ denote the $k$-skeleton of
$B$. We may assume that the restriction $E_{B^{(k)}}\to B^{(k)}$
has the claimed property. Let $(e^{k+1},\dot e^{k+1})\to (B,B^{(k)})$
be a cell, and consider the induced bundle
\[
\begin{array}{ccc}
E'  & \to & E \\
\dto\rlap{\scriptsize$\uto s'$} & & \dto\rlap{\scriptsize$\uto s$} \\
e^{k+1} & \to & B.
\end{array}
\]
Since $e^{k+1}$ is contractible, the bundle 
$E'\TO{\scriptstyle\OT{s'}}e^{k+1}$
is isomorphic to the trivial bundle $e^{k+1}\times\SS^m\to e^{k+1}$.
Thus we find a map 
\[
(e^{k+m+1},\dot e^{k+m+1})\to (E',E'_{\dot e^{k+1}}\cup s'(e^{k+1}))\to
(E,E_{ B^{(k)}} \cup s( B^{(k)} ) ).
\]
\qed}

\section{Generalized manifolds}
\label{GenMnf}

Let $R$ be a principal ideal domain.
We denote sheaf-theoretic cohomology by ${\bf\bar H}^\bullet$, and
Borel-Moore homology by ${\bf\bar H}_\bullet$. Compact supports
are indicated by the letter $c$. Borel-Moore homology
is defined only for locally compact Hausdorff spaces, and it
has some unusual features: it is not clear that
$\Hbar_{-1}(X;R)=0$, and the group $\Hbar_\bullet(X;R)$ depends
on whether $R$ is viewed as an $\ZZ$- or an $R$-module. 
However, these atrocities disappear for 'good' spaces
\cite[V]{Bredon65}.

The notation follows essentially Bredon's book \cite{Bredon65}.

\Lemma{\label{PreSheaf}
Let $S:U\mapsto S(U)$ be a presheaf on a space $X$, and 
let $\cal S$ denote the sheaf that it generates. Let $W\SUB X$
be open. The restriction ${\cal S}|W$ is canonically isomorphic
to the sheaf generated by the restriction of $S$ to $W$.

\proof This is clear from the construction of $\cal S$, see
\cite[I 1.3]{Bredon65}.\qed}

\Prop{\label{CompareHom}
Let $(X,A)$ be a pair. If $X$ and $A$ are paracompact and
$HLC$, then there is a natural isomorphism
\[
\H^\bullet(X,A;R)\cong\Hbar^\bullet(X,A;R)
\]
\cite[III 2.]{Bredon65}.
If $X$ and $A$ are locally compact and $HLC$, then there is a natural
isomorphism
\[
\H_\bullet(X,A;R)\cong\Hbar_\bullet^c(X,A;R)
\]
\cite[V 12.6]{Bredon65}.
Note that $\Hbar_\bullet^c(X,X-\{x\};R)\cong\Hbar_\bullet(X,X-\{x\};R)$
by \cite[V 5.8]{Bredon65}.\qed}

\Def{A locally compact Hausdorff space $X$ is called a
{\em generalized $n$-manifold over $R$}, if it has 
the following properties:
\begin{enumerate}
\item the sheaf-theoretic dimension $\dim_R X$ is 
finite \cite[II 15.]{Bredon65}.
\item $X$ is cohomological locally connected ($clc_R$), 
i.e. for every point
$x\in X$ and every open neighborhood $U$ of $x$, there exists an
open neighborhood $V\SUB U$ of $x$, such that the induced map
${\bf\tilde\Hbar}{}^\bullet(V;R)\ot{\bf\tilde\Hbar}{}^\bullet(U;R)$ is
trivial.
\item for each $x\in X$,
\[
\Hbar_i(X,X-\{x\};R)=
\left\{
\begin{array}{ccc}
R   & \mbox{for} & i=n \\
0   & \mbox{for} & i\neq n.
\end{array}
\right.
\]
\item The {\em orientation sheaf}, i.e. the sheaf generated by the presheaf 
$V\mapsto\Hbar_\bullet(X,X-V;R)$ is locally constant.
\end{enumerate}

Thus, a generalized $n$-manifold over $R$ is what Bredon \cite{Bredon65}
calls an $n$-$cm_R$.

We call $X$ a {\em generalized $n$-manifold}, if it is a generalized
$n$-manifold over every principal ideal domain $R$. We are mainly interested
in the cases $R=\ZZ,\ZZ_p,\QQ$.}

\Prop{\label{LocGenMnf}
Let $X$ be a locally compact Hausdorff space. If every
point $x\in X$ has an open neighborhood $U$ such that $U$ is 
a generalized $n$-manifold, then $X$ is a generalized $n$-manifold.
Conversely, every open subset of a generalized $n$-manifold is
again a generalized $n$-manifold.

A product of a generalized $n$-manifold and a generalized $m$-manifold
is a generalized $(m+n)$-manifold \cite[V 15.8]{Bredon65}.

\proof If every point has a neighborhood which has
finite sheaf-theoretic dimension, then $\dim_RX$ is finite
as well \cite[II 15.8]{Bredon65}. The orientation sheaf is locally
constant by \ref{PreSheaf}.\qed}

\Thm{\label{ANRMnf}
Let $X$ be a locally compact ANR of finite covering dimension.
Then the following are equivalent:
\begin{enumerate}
\item There exists a number $n$ such that for every $x\in X$
\[
\H_i(X,X-\{x\})=
\left\{
\begin{array}{ccc}
\ZZ & \mbox{for} & i=n \\
0   & \mbox{for} & i\neq n.
\end{array}
\right.
\]
\item For all $x,y\in X$,
the groups $\H_i(X,X-\{x\})$ and $\H_i(X,X-\{y\})$ are finitely
generated and isomorphic.
\item The space $X$ is a generalized $n$-manifold over $\ZZ$.
\item The space $X$ is a generalized $n$-manifold.
\end{enumerate}
If one of these equivalent conditions is satisfied, then
$\dim_R X=\dim_\ZZ X=\dim X=n$.

\proof Being a generalized manifold is a local property, hence
we may assume that $X$ is second countable by passing to a relatively
compact open subset. 

Since the covering dimension of $X$ is finite, its sheaf-theoretic
dimension is finite as well. The space $X$ is locally contractible
and hence $HLC$ and $clc_R$. Now $\Hbar_\bullet(X,X-\{x\})=
\H_\bullet(X,X-\{x\})$ by \ref{CompareHom}.
By \cite[3.2]{Bre69b}, $X$ is a generalized $n$-manifold over $\ZZ$.
This establishes the equivalence of (i), (ii), and (iii).

Now suppose that (i) holds, and let $R$ be a principal ideal
domain. It follows form the universal coefficient theorem that 
\[
\H_\bullet(X,X-\{x\})\otimes R\cong\H_\bullet(X,X-\{x\};R)
\cong\Hbar_\bullet(X,X-\{x\};R)
\]
and thus we infer from \cite{Bre69a} that $X$ is a generalized
manifold over $R$.\qed}

Now L\"owen's Theorem is an immediate consequence:

\Thm{\label{Lowen}{\bf (L\"owen)} \cite{Loew}
Let $X$ be a compact ANR of finite covering dimension $\dim X=n$. 
Suppose that for every $x\in X$, the complement $X-\{x\}$ is 
acyclic (i.e. ${\bf\tilde\H}{}_\bullet(X-\{x\})=0$).
Then $X$ is an (orientable) generalized $n$-manifold, and a
homology $n$-sphere. 

If $n\leq 2$, then $X$ is an $n$-sphere.

If $\pi_1(X-\{x\})=0=\pi_1(X-\{y\})$ for two elements $x,y\in X$, then
$X$ is homotopy equivalent to an $n$-sphere.

\proof Since $X$ is locally contractible and compact, its
homology is finitely generated in each dimension \cite[6 9.11]{Spanier}.
The complement of each $x\in X$ is acyclic, and therefore
all local homology groups in dimension $i$ are isomorphic to $\H_i(X)$.
By \ref{LocGenMnf}, $X$ is a generalized $n$-manifold.
Thus
${\bf\tilde\H}{}_\bullet(X;R)\cong\H_{\bullet}(X,X-\{x\};R)\cong
{\bf\tilde\H}{}_\bullet(\SS^n;R)$,
and therefore $X$ is $R$-orientable.

If $\dim X\leq 2$, then $X$ is locally euclidean \cite[p.272]{Wil},
and the claim follows from the classification of the compact
orientable $1$- and $2$-manifolds.

If $n>1$, then the Mayer-Vietoris
sequence yields $\H_0(X-\{x,y\})=\ZZ$, and from the Seifert-Van Kampen
Theorem we conclude that $\pi_1X=0$. Thus $X$ is homotopy equivalent
to a sphere by \ref{CWSphere}.\qed}

\Lemma{Let $X$ be a locally compact, finite-dimensional ANR.
The (unreduced) open cone $M=C'_X-X$ is a generalized $(n+1)$-manifold
if and only if $X$ is a generalized $n$-manifold with
$\H_\bullet(X)\cong\H_\bullet(\SS^n)$.

\proof Clearly, $C'_X-X$ is a finite-dimensional ANR.
Let $x$ denote the tip of the cone. Now $\H_\bullet(C'_X-(X\cup\{x\}))=
\H_\bullet(X\times (0,1))\cong\H_\bullet(X)$, and the claim
follows from \ref{ANRMnf}.\qed}

\Cor{\label{ManBundle}
Let $E\to B$ be a locally trivial bundle with fiber $F$. Suppose that
$E,B,F$ are finite-dimensional ANRs.
If $E$ is a generalized manifold, then the fiber $F$ and the base $B$ are
generalized manifolds. The (unreduced) open mapping cylinder is
a generalized manifold if and only if
$\H_\bullet(F)=\H_\bullet(\SS^{\dim F})$.

\proof This follows from the previous lemma and
\cite[V 15.8]{Bredon65}.\qed}

\section{M\"unzner's Theorem}
\label{MueThm}

\Thm{\label{Muenzner}{\bf(M\"unzner)}
Let $\F\TO{\pr_1}\P$ and $\F\TO{\pr_2}\L$ be locally trivial bundles
with homotopy $m'$- and $m$-spheres as fibers.

Suppose that the fibers and the bases of these bundles are locally compact,
finite dimensional ANRs. 

Assume moreover that $\P,\L$ and $\F$ are generalized manifolds 
of dimension $r-m',r-m$ and $r$, respectively. 

Consider the double mapping cylinder $D\F$ over $\pr_1,\pr_2$.
If $D\F$ is a homology $(r+1)$-sphere, then
\begin{eqnarray*}
\H^\bullet(\F;\ZZ_2) & \cong & \ZZ_2^{2n} \\
\H^\bullet(\P;\ZZ_2) & \cong & \ZZ_2^n \\
\H^\bullet(\L;\ZZ_2) & \cong & \ZZ_2^n
\end{eqnarray*}
where $n\in\{1,2,3,4,6\}$. 
If $m,m'>1$, then
\begin{eqnarray*}
\H^\bullet(\F) & \cong & \ZZ^{2n} \\
\H^\bullet(\P) & \cong & \ZZ^n \\
\H^\bullet(\L) & \cong & \ZZ^n.
\end{eqnarray*}

More specifically, the cohomology rings of $\P,\L$ and $\F$ is as
follows. Here, the natural inclusions 
between the listed rings correspond to the monomorphisms
\[
\H^\bullet(\P;R)\TO{\pr_1^\bullet}\H^\bullet(\F;R)
\OT{\pr_2^\bullet}\H^\bullet(\L;R).
\]
Compare also Hebda \cite{Hebda} for the $\ZZ_2$-cohomology.
The subscripts indicate the degree of the cohomology classes.

\begin{enumerate} 
\item[$\bf1$] $n=1$ and $m=m'$.
The cohomology of $\P\ot\F\to\L$ is the same as that of
\[
{*}\ot\SS^m\to{*}
\]

\item[$\bf2$] $n=2$, $m,m'$ arbitrary.
The cohomology of $\P\ot\F\to\L$ is the same as that of
\[
\SS^m\ot\SS^m\times\SS^{m'}\to\SS^{m'}
\]

\item[$\bf3_1$] $n=3$ and $m=m'=1$.
\begin{eqnarray*}
\H^\bullet(\P;\ZZ_2) & = & \ZZ_2[x_1]/(x_1^3) \\
\H^\bullet(\L;\ZZ_2) & = & \ZZ_2[y_1]/(y_1^3) \\
\H^\bullet(\F;\ZZ_2) & = & \ZZ_2[x_1,y_1]/(x_1^3,y_1^3,x_1^2+x_1y_1+y_1^2).
\end{eqnarray*}

\item[$\bf3_2$] $n=3$ and $m=m'\in\{2,4,8\}$
\begin{eqnarray*}
\H^\bullet(\P) & = & \ZZ[x_m]/(x_m^3) \\
\H^\bullet(\L) & = & \ZZ[y_m]/(y_m^3) \\
\H^\bullet(\F) & = & \ZZ[x_m,y_m]/(x_m^3,y_m^3,x_m^2-x_my_m+y_m^2).
\end{eqnarray*}

\item[$\bf4_1$] $n=4$ and $m=m'=1$
\begin{eqnarray*}
\H^\bullet(\P;\ZZ_2) & = & \ZZ_2[x_1]/(x_1^4)  \\
\H^\bullet(\L;\ZZ_2) & = & \ZZ_2[y_1,y_2]/(y_1^2,y_2^2) \\
\H^\bullet(\F;\ZZ_2) & = & \ZZ_2[x_1,y_1,y_2]/(x_1^4,y_1^2,y_2^2,
y_2+x_1^2+x_1y_1) 
\end{eqnarray*}

\item[$\bf4_2$] $n=4$, $m=1$ and $m'>1$.

The $\ZZ_2$-cohomology rings of $\P$, $\L$, and $\F$ are
graded $\ZZ_2$-algebras generated by homogeneous elements $(x_m,x_{m+m'})$,
$(y_{m'},y_{m'+m})$, and $(x_m,y_{m'},x_{m+m'},y_{m'+m})$, respectively,
subject to the relations $x_m^2=x_{m+m'}^2=0$,
$y_{m'}^2=y_{m'+m}^2=0$ and
$x_m^2=y_{m'}^2=x_{m+m'}^2=y_{m'+m}^2=y_{m'+m}+x_{m+m'}+x_my_{m'}=0$
(the subscripts indicate the degrees).

\item[$\bf4_3$] $n=4$, $m,m'>1$, and $m+m'$ odd.

The integral cohomology rings of $\P$, $\L$, and $\F$ are anticommutative
graded $\ZZ$-algebras generated by homogeneous elements $(x_m,x_{m+m'})$,
$(y_{m'},y_{m'+m})$, and $(x_m,y_{m'},x_{m+m'},y_{m'+m})$, respectively,
subject to the relations $x_m^2=x_{m+m'}^2=0$,
$y_{m'}^2=y_{m'+m}^2=0$ and
$x_m^2=y_{m'}^2=x_{m+m'}^2=y_{m'+m}^2=y_{m'+m}+x_{m+m'}-x_my_{m'}=0$
(the subscripts indicate the degrees).

The restriction $m,m'>1$ is not really necessary, cp. \ref{OddCor}.

\item[$\bf4_4$] $n=4$ and $m=m'\in\{2,4\}$. 
\begin{eqnarray*}
\H^\bullet(\P) & = & \ZZ[x_m]/(x_m^4) \\
\H^\bullet(\L) & = & \ZZ[y_m,y_{2m}]/(y_{2m}^2,y_m^2-2y_{2m}) \\
\H^\bullet(\F) & = & \ZZ[x_m,y_m,y_{2m}]/(x_m^4,y_{2m}^2,y_m^2-2y_{2m},
y_{2m}+x_m^2-y_mx_m)
\end{eqnarray*}

\item[$\bf6_1$] $n=6$ and $m=m'=1$.
\begin{eqnarray*}
\H^\bullet(\P;\ZZ_2) & = & \ZZ_2[x_1,x_3]/(x_1^3,x_3^2) \\
\H^\bullet(\L;\ZZ_2) & = & \ZZ_2[y_1,y_3]/(y_1^3,y_3^2) \\
\H^\bullet(\F;\ZZ_2) & = & \ZZ_2[x_1,x_3,y_1,y_3]/
(x_1^3,y_1^3,x_3^2,y_3^2,x_1^2+y_1^2+x_1y_1,x_3+y_3+x_1^2y_1)
\end{eqnarray*}

\item[$\bf6_2$] $n=6$ and $m=m'\in\{2,4\}$.
\[
\H^\bullet(\P) = \ZZ[x_m,x_{3m}]/(x_{3m}^2,x_m^3-2x_{3m})
\]
The rings $\H^\bullet(\L)$ and $\H^\bullet(\F)$ have bases
$\{y_m,y_{2m},y_{3m},y_{4m},y_my_{4m}\}$ and 
$\{x_m,y_m,
\linebreak
x_m^2,y_{2m},x_{3m},y_{3m},x_mx_{3m},y_{4m},x_m^2x_{3m},
y_my_{4m},x_{3m}y_{3m}\}$ respectively. The missing products can easily
be calculated from the equations $y_{m}^2=3y_{2m}$, $y_my_{2m}=2y_{4m}$,
$y_my_{3m}=3y_{4m}$, $y_m^2y_{4m}=0$, $x_my_m=x_m^2+y_{2m}$, 
and the fact that the integral
cohomology modules are torsion-free. For example $3y_{2m}y_{3m}=
y_m^2y_{3m}=3y_my_{4m}$, hence $y_{2m}y_{3m}=y_my_{4m}$.
\end{enumerate}

\proof M\"unzner uses only standard properties of manifolds, like
Poincar\'e and Alexander duality, and the coefficient sequence
$0\to\ZZ_2\to\ZZ_4\to\ZZ_2\to 0$. His proof is purely algebraic,
hence his arguments are also valid for generalized manifolds.
The restrictions on the dimension in the cases $3_2,4_4,6_2$ are
not given in \cite{Mue}, but they follow from the structure of
the cohomology rings and \ref{Adams}.\qed}

\section{Odds and ends}

\Prop{\label{CWSphere}
Suppose $X$ is a simply connected CW complex and a homology 
$m$-sphere. Then $X$ is homotopy equivalent to an $m$-sphere.

\proof Consider the Hurewicz isomorphism $\pi_m(X)\to\H_m(X)$,
$h\mapsto h_\bullet[\SS^m]$. Let $f:\SS^m\to X$ denote the generator 
of $\pi_m(X)$. The map $f$ induces an isomorphism 
$\H_\bullet(\SS^m)\to\H_\bullet(X)$ (because $f_\bullet[\SS^m]$
is a generator of $\H_m(X)$), and therefore $f$ is a homotopy 
equivalence \cite[7.6.24]{Spanier} \cite[VII 11.15]{Bredon93}.\qed}

\Prop{\label{LocSur}
Let $M$, $N$ be compact connected $k$-manifolds, and let $f:M\to N$
be a continuous map. Suppose that there exists an element $x\in M$
with $\{x\}=f^{-1}(f(x))$, and a neighborhood $U$ of $x$ such that
the restriction $f|U$ is injective. Then $f$ is surjective.

\proof By invariance of domain, $f|U:U\to f(U)$ is a homeomorphism.
Consider the diagram
\[
\begin{array}{ccccc}
\H_k(M;\ZZ_2) & \to & \H_k(M,M-\{x\};\ZZ_2) & \ot & \H_k(U,U-\{x\};\ZZ_2) \\
\dto          &     & \dto                  &     & \dto \\
\H_k(N;\ZZ_2) & \to & \H_k(N,N-\{f(x)\};\ZZ_2) & \ot & 
\H_k(f(U),f(U)-\{f(x)\};\ZZ_2) .
\end{array}
\]
The vertical arrow at the right is an isomorphism. The horizontal arrows
on the left-hand side are isomorphisms, because $M,N$ are 
$\ZZ_2$-orientable, and the horizontal
arrows on the right-hand side are excision
isomorphisms. Thus $f_\bullet:\H_k(M;\ZZ_2)\to\H_k(N;\ZZ_2)$ is an
isomorphism. Suppose there is an element $y\in N-f(M)$. Then $f$ factors
as $M\to N-\{y\}\to N$. But $\H_k(N-\{y\};\ZZ_2)=0$, a contradiction.\qed}

\Thm{{\bf(Generalized Poincar\'e conjecture)}\label{Poincare}
Let $M$ be a simply connected $n$-manifold. If $M$ has the same
homology as an $n$-sphere, and if $n\geq 4$, then $M$ is homeomorphic
to a sphere.

\proof For $n=4$, this is included in Freedman's classification
of $4$-manifolds \cite{Freedman}.

For $n\geq 5$, this follows from Smale's proof
\cite{Smale} of the combinatorial generalized Poincar\'e conjecture
and the fact that the Kirby-Siebenmann obstruction \cite{KS} vanishes,
since $\H^4(M;\ZZ_2)=0$. There is also a direct proof by Newman
\cite{Newman}.\qed}

\Thm{{\bf(Adams-Atiyah)}\label{Adams} \cite{AA}
Let $X$ be a space having the homotopy type of a finite CW complex.
Assume that the integral cohomology of $X$ has no torsion, and that the
integral cohomology groups of $X$ vanish in all dimensions different 
from $0,2m,4m,8m,\ldots$.

If there is an $x\in\H^{2m}(X;\ZZ_2)$ with $x^2\neq 0$, then
$2m\in\{2,4,8\}$. 

If there is an $x\in\H^{2m}(X;\ZZ_3)$ with $x^3\neq 0$, then
$2m\in\{2,4\}$.\qed}

\end{document}